\documentclass[12pt]{article}
\usepackage{amsmath,amsfonts,amssymb}
\usepackage{alltt}

\usepackage{amsmath,amsfonts,ifthen,fullpage,enumerate}  
\usepackage{mathrsfs}
\usepackage{hyperref}
\usepackage{mathtools}
\usepackage{tikz-cd}
\usetikzlibrary{arrows.meta}
\tikzset{>={Latex[width=2mm,length=2mm]}}

\usepackage{colortbl}

\usepackage{enumerate}

\usepackage{amsmath,amsfonts,amssymb}
\usepackage{float}

\usepackage{makecell}  

\usepackage{graphicx}        

\usepackage{tikz-cd}

\usepackage{tcolorbox}
\def\greyboxvar#1#2
{\begin{tcolorbox} [boxrule=-5pt, boxsep=-0.1cm, colback=white!95!black, 
	width=#2] #1\end{tcolorbox}}

\def\greybox#1
{
	\begingroup
\setlength{\tabcolsep}{2pt} 
\renewcommand{\arraystretch}{1.3} 
	\begin{tabular}{c}
\cellcolor{gray!25}#1\end{tabular}
\endgroup
}

\def\greyeq#1{\colorbox{black!8!white}{$\displaystyle #1$}}

\usepackage[all,pdf]{xy}
\usepackage{lmodern,amssymb}   

\def\gcd{\mathop{\rm gcd}\nolimits}
\def\lcm{\mathop{\rm lcm}\nolimits}

\def\divides{\mid}
\def\ndivides{\nmid}

\def\Re{\mathop{\rm Re}\nolimits}

\def\mult{\mathop{\rm mult}\nolimits}
\def\spec{\mathop{\rm spec}\nolimits}

\def\trace{\mathop{\rm trace}\nolimits}

\def\rank{\mathop{\rm rank}\nolimits}

\def\pmat#1{\begin{pmatrix}#1\end{pmatrix}}
\def\mat#1{\begin{matrix}#1\end{matrix}}
\def\question#1{{\bf Question: }#1}
\def\question#1{}

\def\cG{{\cal G}}

\def\cO{{\cal O}}

\def\cT{{\cal T}}

\def\CC{\mathbb{C}}

\def\HH{\mathbb{H}}

\def\Cd{\C^d}

\def\Hd{\HH^d}

\def\C{\mathbb{C}}

\newcommand{\RR}{\mathbb{R}}

\newtheorem{theorem}{Theorem}[section]
\newtheorem{corollary}{Corollary}[section]
\newtheorem{lemma}{Lemma}[section]
\newtheorem{example}{Example}[section]
\newtheorem{remark}{Remark}[section]
\newtheorem{proposition}{Proposition}[section]

\newenvironment{proof}{{\noindent \it
Proof.}}{\hfill$\Box$\medskip}
\newif\ifdraft\def\draft{\drafttrue\hoffset=.8truecm\showlabeltrue
\def\comment##1{{\bf comment: ##1}}
\headline={\sevenrm \hfill \ifx\filenamed\undefined\jobname\else\filenamed\fi%
(.tex) (as of \ifx\updated\undefined???\else\updated\fi)
 \TeX'ed at {\hour\time\divide\hour by 60{}%
\minutes\hour\multiply\minutes by 60{}%
\advance\time by -\minutes
\the\hour:\ifnum\time<10{}0\fi\the\time\  on \today\hfill}}
}

\def\inpro#1{\langle#1\rangle}
\def\ip#1{\langle\kern-.28em\langle#1\rangle\kern-.28em\rangle_\nu}

\def\cC{{\cal C}}
\def\cD{{\cal D}}
\def\cO{{\cal O}}
\def\cT{{\cal T}}
\def\cI{{\cal I}}
\def\cH{{\cal H}}

\def\openR{{{\rm I}\kern-.16em {\rm R}}}

\let\ga\alpha

\let\gb\beta

\let\gTh\Theta

\let\gl\lambda

\let\gL\Lambda

\let\go\omega
\let\gO\Omega
\let\ga\alpha

\let\gb\beta

\def\mC{\mathsf{C}}
\def\mD{\mathsf{D}}
\def\mG{\mathsf{G}}
\def\mH{\mathsf{H}}
\def\mT{\mathsf{T}}
\def\mO{\mathsf{O}}
\def\mI{\mathsf{I}}
\def\inpro#1{\langle#1\rangle}

\def\ker{\mathop{\rm ker}\nolimits}
\def\GL{\mathop{\it GL}\nolimits}

\def\Iff{\hskip1em\Longleftrightarrow\hskip1em}
\def\Implies{\hskip1em\Longrightarrow\hskip1em}

\def\formeq{\the\sectionno.\the\equationno}  
\def\elabel#1/#2/#3/{\global\advance\equationno by 1 %
\ifx#1\empty\else\emember#1%
\ifshowlabel\marginal{\string#1}\fi\fi%
\ifmmode\eqno{#3(\formeq#2)}\else#3\formeq#2\fi} 

\def\makeblanksquare#1#2{
\dimen0=#1pt\advance\dimen0 by -#2pt
      \vrule height#1pt width#2pt depth0pt\kern-#2pt
      \vrule height#1pt width#1pt depth-\dimen0 \kern-#1pt
      \vrule height#2pt width#1pt depth0pt \kern-#2pt
      \vrule height#1pt width#2pt depth0pt
}

\title{\bf 
The quaternionic reflection groups with imprimitive complexifications
and the finite subgroups of $GL(\CC^2)$
}

\author{Shayne Waldron\\ 
 \\
Department of Mathematics \\ University of Auckland\\
Private
Bag 92019, Auckland, New Zealand\\
e--mail: waldron@math.auckland.ac.nz}

\begin{document}

\maketitle 

\begin{abstract}

We give a transparent development of the quaternionic reflection groups 
which have 
imprimitive complexifications, including a description of their 
reflections and the corresponding normal reflection subgroups.
They are seen to be of rank one or rank two, and are obtained from minimal 
complex reflection groups by 
the addition of one quaternionic reflection and a 
	complex scalar matrix of finite order.
In the case of the rank two groups, we use our techniques to consider 
the broader question of what are
the minimal irreducible subgroups of $GL(\CC^2)$ with a given collineation 
	group.
This leads to a complete classification of the finite irreducible subgroups of $GL(\CC^2)$ in terms of the minimal groups,
with the groups giving rise to quaternionic reflection groups located within it.

%
\end{abstract}

\bigskip
\vfill


 \ifdraft (Other designs, configurations) \else\fi
15B33, \ifdraft (Matrices over special rings (quaternions, finite fields, etc.) \else\fi
20C25, \ifdraft Projective representations and multipliers \else\fi
20G20, \ifdraft Linear algebraic groups over the reals, the complexes, the quaternions \else\fi
51M05, \ifdraft (Euclidean geometries (general) and generalizations) \else\fi
51M20, \ifdraft (Polyhedra and polytopes; regular figures, division of spaces [See also 51F15]) \else\fi
\quad
secondary
15B57, \ifdraft (Hermitian, skew-Hermitian, and related matrices) \else\fi
51E99, \ifdraft Geometry,  None of the above, but in this section \else\fi
51M15, \ifdraft (Geometric constructions in real or complex geometry) \else\fi
65D30. \ifdraft (Numerical integration) \else\fi

\vskip .5 truecm
\hrule
\newpage

\section{Introduction}

A group $G$ of complex matrices in $M_{2d}(\CC)$ gives rise to
a group of quaternionic matrices in $M_d(\HH)$ if it
is conjugate to a group of matrices of the 
{\bf symplectic  form}
\begin{equation}
\label{C2dHdcorres}
[A+Bj]_\CC:=\pmat{A&-B\cr\overline{B}&\overline{A}}\in M_{2d}(\CC)
\Iff A+Bj  \in M_d(\HH).
\end{equation}
A finite matrix group over 
$\RR$, $\CC$ or $\HH$ 
is 
a {\bf reflection group} if it is generated by reflections
(matrices which pointwise fix a hyperplane).
Here we consider the quaternionic reflection groups whose symplectic matrices have 
the block monomial symplectic forms
\begin{equation}
\label{ABmonomialform}
\pmat{A&0\cr0&\overline{A}} \iff A, \qquad
\pmat{0&-B\cr\overline{B}&0} \iff Bj , \qquad A,B\in M_d(\CC),
\end{equation}
i.e., the complexification has a system of imprimitivity consisting
of two $d$-dimensional complex subspaces. 
Cohen \cite{C80} shows that these are the only possible
complex system of imprimitivity for a primitive quaternionic reflection group,
and that the block monomial form of (\ref{ABmonomialform}) can 
be obtained after a suitable change of variables. 
In the classification of Schmitt \cite{S23} Fig.\ 1.2.1,  the groups we consider are
the ``symplectically irreducible complex irreducible symplectically primitive
complex imprimitive symplectic reflection groups''.


We will see from an elementary argument (Lemma \ref{BjReflects})
that reflections of the type $Bj$ are highly restricted. As a consequence
(the reflections of the type $A$ being well understood),
all of the quaternionic reflection groups 
with matrices of the form (\ref{ABmonomialform}) 
can be written
$$ \inpro{G,\zeta,s}, \qquad \rank(G)=1,2, \qquad 
s=\pmat{j},\pmat{0&j\cr-j&0} \quad\hbox{(respectively)}, $$
where $G$ is a complex reflection group, $\zeta$ is a complex 
scalar matrix, and the $s$ above is a quaternionic reflection
of the form $Bj$.
From this observation, we proceed as follows:
\begin{itemize}
\item We determine the minimal collineation preserving reflection subgroups 
$G$ in $\GL(\CC^2)$.
\item We consider the structure of the reflection groups $\cH_a(G)$ obtained 
from them by adding $s$ and a scalar matrix $\zeta$ of finite order
(Lemma \ref{internaldirectsum}), and determine when a given group
can be obtained from different pairs $(G,\zeta)$, thereby 
		classifying all such reflection groups
(Theorems \ref{TOIfamilyclass} and \ref{G(n,a,b,r)familyclass}).
\item We undertake the corresponding classification of the finite 
subgroups of $\GL(\CC^2)$ in terms of the minimal subgroups with
		a given collineation group,
choosing the parameters so that no groups are repeated
(Theorems \ref{GL2Csubgroups-primitive}, \ref{GL2Csubgroups-imprimitive},
and \ref{GL2Csubgroups-all}). 
\item We then investigate the lattice of normal reflection subgroups of the
	quaternionic reflection groups $\cH_a(G)$ in terms
of the indexing lattice of \cite{W26c}. By way of comparison with the 
complex reflection groups where the indexing is $1$--$1$,
it is seen that this is often not $1$--$1$ (Theorem \ref{HaGnormalrefgroups}
). 
Indeed, there are quaternionic reflection groups which do not require at 
least one reflection of each order to generate it.
\end{itemize}
We give numerous examples, tables and figures which summarise the 
algebraic and geometric structures of the complex and quaternionic groups that we consider.


\section{Quaternionic reflections which have imprimitive complexifications}

We assume, without loss of generality (since finite subgroups
of $GL(\CC^2)$ are conjugate to unitary ones), that the matrices
$A$ and $Bj$ of (\ref{ABmonomialform}) are unitary, which by the calculation
$$ (Bj)(Bj)^*=Bjj^*B^*=BB^*, $$ 
is equivalent to the complex matrices $A$ and $B$ being unitary.
Reflections of the type $A$ are simply complex reflections. 
We assume basic familiarity with complex reflection groups (e.g.,
\cite{ST54}, \cite{C76}, \cite{OT92}, \cite{BMM99}, \cite{LT09}, 
\cite{B10},
also \cite{BST23}, \cite{RS26}),

We now determine the reflections of the type $Bj$,
by using
\begin{equation}
\label{jzcommute}
jz=\overline{z}j, \qquad\forall z\in\CC. 
\end{equation}

\begin{lemma}
\label{BjReflects}
A matrix $Bj$, $B\in M_d(\CC)$, is a unitary reflection in precisely 
the 
	cases
\begin{enumerate}[\rm(i)]
	\item $d=1$ and $Bj$ is of order $4$ given by $Bj=\pmat{\xi j}$, $\xi\in\CC$, $|\xi|=1$.
\item $d=2$ and $Bj$ is of order $2$ given by the monomial matrix
\end{enumerate}
\begin{equation}
\label{Bform}
Bj=\pmat{\xi&0\cr0&\xi}\pmat{0&j\cr-j&0}=\pmat{0&\xi j\cr-\xi j&0}, \qquad
\xi\in\CC,\quad |\xi|=1.
\end{equation}
\end{lemma}

\begin{proof} We first observe that 
$$ (Bj)^{2m+1}=(BjBj)^mBj=(-B\overline{B})^mBj\ne I, $$
so that no matrix $Bj$ can have an odd order.
The reflection $r$ in $U(\Hd)$ with root $a\in\Hd$ and $ra=a\gl$, $\gl\in\HH$,
$|\gl|=1$, is given by the formula
$$ r_{a,\gl}=I-a(1-\gl)a^*=Bj, \qquad B=(I-a(1-\gl)a^*)(-j), \qquad\inpro{a,a}=1. $$
The condition for this $B\in M_d(\HH)$ to be complex is that its $j$ and $k$ parts are zero, i.e.,
$$ B+iBi-jBj+kBk=0, \qquad B+iBi+jBj-kBk=0. $$
\noindent
Now  
$$ B+iBi  
	= 2I(-j)+a(1-\gl)a^*j-2ia(1-\gl)a^*k=0 \Iff 2I=a(1-\gl)a^*-ia(1-\gl)a^*i, $$
$$ jBj-kBk 
	= 2I(j)-ja(1-\gl)a^*-ka(1-\gl)a^*i=0 \Iff  2I=a(1-\gl)a^*-ia(1-\gl)a^*i, $$
so that $B$
is complex if and only if
\begin{equation}
\label{2Iformula}
2I=a(1-\gl)a^*-ia(1-\gl)a^*i. 
\end{equation}
Let $a=u+vj$, $u,v\in\Cd$, so that $\inpro{a,a}=\inpro{u,u}+\inpro{v,v}=1$,
and $\gl=\ga+\gb j$, $\ga,\gb\in\CC$. 
Then repeatedly using (\ref{jzcommute}) gives
\begin{align}
\label{aglaformula}
a(1-\gl)a^*
&= \{(1-\ga)uu^*+\overline{(1-\ga)}vv^*-\gb uv^*+\overline{\gb}vu^*\} \cr
&\quad +\{-(1-\ga)uv^T +\overline{(1-\ga)}vu^T -\gb uu^T-\overline{\gb}vv^T\}j, 
\end{align}
so that   
	(\ref{2Iformula}) divided by $2$ becomes
\begin{align}
\label{basicformula}
I &=(1-\ga)uu^* +\overline{(1-\ga)}v v^* -\gb u v^* +\overline{\gb}v u^* \cr
&=\{(1-\ga)u+\overline{\gb}v\} u^* +\{\overline{(1-\ga)}v-\gb u\} v^* .
\end{align}
We observe the right hand side of (\ref{basicformula}) has
rank at most $2$, so we can only have $d=1,2$.
Taking the trace of (\ref{basicformula}) gives
$$ d=(1-\ga)\inpro{u,u}+\overline{(1-\ga)}\inpro{v,v}
-\gb\inpro{v,u}+\overline{\gb}\inpro{u,v}, $$
and since $-\gb\inpro{v,u}+\overline{\gb}\inpro{u,v}$ is pure imaginary,
taking the real part of the above gives
$$ d=\Re(1-\ga)(\inpro{u,u}+\inpro{v,v}) =1-\Re(\ga) \Implies
	\Re(\ga)=1-d. $$

For $d=1$, $\xi$ is unitary if and only if $|\xi|=1$, and the calculation
$$ (\xi j)^2=\xi j\xi j=\xi\overline{\xi}j^2=-1, $$
shows that $Bj=r_{1,\xi j}=\xi j$ is a reflection of order $4$.

For $d=2$, we have $\Re(\ga)=1-d=-1$, so that $\gl=\ga=-1$ and $\gb=0$.
By (\ref{basicformula}), we have $I=2(uu^*+vv^*)$, so that $u$ and $v$ are orthogonal and of
length ${1\over\sqrt{2}}$.
Hence by (\ref{aglaformula}), 
the reflection is given by
$$ r_{a,-1}= I-2aa^*=2(uv^T-vu^T)j = 2\pmat{0&u_1v_2-v_1u_2\cr
-(u_1v_2-v_1u_2)&0}j, $$
and so is monomial. The following choice for the root vector $a$ gives
$$ u={1\over\sqrt{2}}\pmat{1\cr0}, \quad
v={1\over\sqrt{2}}\pmat{0\cr\xi} \Implies a={1\over\sqrt{2}}\pmat{1\cr\xi j}
\Implies r_{a,-1}=\pmat{0&\xi j\cr-\xi j&0}. $$
and so we obtain all the reflections asserted.
\end{proof}

\begin{remark}
A unitary matrix $Bj$ has order $2$ if and only if
$$ (Bj)^2=Bj Bj=B\overline{B}j^2= - B(B^T)^*=I \Iff
-B=B^T, $$
i.e., $B$ is skew symmetric, and so for $d=2$ all
elements $Bj$ of order $2$ are reflections.

It is possible to construct matrices $Bj$ of order $2m$.
For example, one may take a unitary $B$ of the form
$$ B=\xi \pmat{a&\overline{b}\cr -b&\overline{a}}, \qquad |\xi|=1. $$
A simple calculation shows that $\det((Bj)^2)=\det(-B\overline{B})=1$,
so the eigenvalues of $-B\overline{B}$ (of order $m$) must be
	$\zeta$, $\overline{\zeta}$, with $\zeta$ a primitive $m$-th root of unity.
By considering 
$$ \trace(-B\overline{b})=\zeta+\overline{\zeta}, $$
we obtain the following sufficient condition for $Bj$ to have order $2m$
$$ \Re(b)=\cos{\pi\over m}, $$
and as a particular case, we may take
$$ B_m=\xi\pmat{\sin{\pi\over m}&\cos{\pi\over m}\cr-\cos{\pi\over m}&\sin{\pi\over m}}, \qquad m=1,2,\ldots. $$
\vskip-0.7truecm
\noindent
This gives
$$ B_1j=\pmat{0&\xi j\cr-\xi j&0}, \quad
B_2j=\xi\pmat{1&0\cr0&1}j, \quad
B_3j=\xi\pmat{{\sqrt{3}\over2}&{1\over2}\cr-{1\over2}&{\sqrt{3}\over2}}j, \quad
	B_4j={\xi\over\sqrt{2}}\pmat{1&1\cr-1&1}j, 
	$$
of orders $2$, $4$, $6$ and $8$. There are many more such matrices that can obtained by 
starting with $B$ of a more general form, e.g.,
$$ 
Bj={\xi\over\sqrt{2}}\pmat{1&1\cr1&-1} \ \hbox{(order $4$)}, \qquad
Bj={\xi\over\sqrt{2}}\pmat{1&-i\cr1&i} \ \hbox{(order $6$)}.
$$
\end{remark}

The $\xi j$ for the reflections of the type (\ref{Bform}) in a finite monomial group
form what is called a {\it reflection system} in \cite{W25}. In particular, they 
are closed under the operation $a\circ b:=ab^{-1}a$, which in this case is
$$ \xi j\circ \go j= \xi j(\go j)^{-1}\xi j
=\xi j j^*\overline{\go} \xi j =\xi \overline{\go}\xi
=(\xi\circ\go) j, $$
i.e., the complex factors $\xi$ are closed under $a\circ b$.
It is always possible to apply a change of variables to a monomial group 
so that a reflection system 
transforms to one containing any given unit quaternion (in the formal definition 
this is taken to be $1$). In our situation, this implies that if there is 
reflection of type $Bj$ in a group, then we can take it to be $s$,
e.g., by the explicit change of variables
$$ \pmat{j&0\cr0&\xi j}\pmat{0&\xi j\cr-\xi j&0}\pmat{j&0\cr0&\xi j}^{-1}
=\pmat{0&-\overline{\xi}\cr1&0}\pmat{-j&0\cr0&-\xi j}
=\pmat{0&j\cr-j&0}. $$
It then follows that all rank two quaternionic reflection groups with 
with imprimitive complexifications have the form:
\begin{equation}
\label{generalform}
	\inpro{G,\zeta_m,s}, \qquad \zeta_m:=e^{2\pi i\over m}, \quad s:=\pmat{0&j\cr-j&0},
\end{equation}
where $G$ is a complex reflection group, and the reflections of type (\ref{Bform}) 
in the group are given by $\xi s$ with $\xi\in\inpro{Z(G),\zeta_m}$ 
(we recall by Schur's lemma that if $G$ is irreducible then its
centre consists of scalar matrices).
Before classifying all such groups and determining their reflections, there
are two technical points to consider:
\begin{enumerate}
\item If a complex unitary change of variables $UGU^{-1}$ is applied to $G$,
then does one still get the same group?
\item Does $\inpro{G,\zeta_m,s}$ give a finite group, and what is its structure?
\end{enumerate}
For the first point, we calculate that
$$ UsU^{-1}=U\pmat{0&1\cr-1&0}jU^*
=U\pmat{0&1\cr-1&0}U^Tj=\det(U)\pmat{0&1\cr-1&0}j=\det(U)s, $$
so that $s$ is fixed by the change of variables $U$ if and only if
$\det(U)=1$. Therefore, with $\hat U:=\ga U$,
$\ga:=1/\sqrt{\det(U)}$, we have $\det(\hat U)=1$, and hence
$$ \hat U\inpro{G,\zeta_m,s} \hat U^{-1} = 
\inpro{ (\ga U)G(\ga U)^{-1},\hat U \zeta_m \hat U^{-1}, \hat U s \hat U^{-1}}
= \inpro{ UGU^{-1},\zeta_m,s}, $$
i.e., the same quaternionic reflection group is obtained.

For the second point, we first consider how 
a matrix $A\in GL(\CC^2)$ is conjugated by $s$  
$$ A^s
=\pmat{a&b\cr c&d}^s
=\pmat{\overline{d}&-\overline{c}\cr-\overline{b}&\overline{a}}
=\Bigl\{\det(\overline{A}){1\over \det(\overline{A})}
\pmat{\overline{d}&-\overline{b}\cr-\overline{c}&\overline{a}}\Bigr\}^T
        =\det(\overline{A}) (\overline{A}^{-1})^T, $$
with $\det(A^s)=\overline{\det(A)}$.
Hence, 
for $A=g$ a unitary matrix, we have
\begin{equation}
\label{gsconjrelation}
g^s=sgs^{-1}=\det(g^*) g, \qquad g\in U(\CC^2).
\end{equation}
Thus $\inpro{G,\zeta_m,s}$ gives a finite group 
\begin{equation}
\label{hatGdefn}
\inpro{G,\zeta_m,s}=\inpro{\hat G,\zeta_m,s}, \qquad \hat G:=\inpro{G,G^s}
	=\inpro{G,\det(G)},
\end{equation}
with its structure
(multiplication) determined by (\ref{gsconjrelation}).
By (\ref{gsconjrelation}), 
$\hat G=\inpro{G,G^s}$
is a finite complex reflection group (generated by the
reflections in $G$ and $G^s$), which is 
closed under conjugation by $s$. 
Thus, it suffices to consider only the complex reflection groups 
that are closed under conjugation by $s$. This is (for better or worse)
most of them.

\begin{figure}[H]
\begin{center}
	\begin{tabular}[t]{ p{7.5truecm} p{2.8truecm} p{2.9truecm} }
\begin{tikzpicture}
    \matrix (A) [matrix of nodes, row sep=1.0cm, column sep = 0.3 cm]
    {
	    && $G_{7},\inpro{\zeta_{12}}$ && \\
	 $G_{6},\inpro{\zeta_{4}}$ && $G_{5},\inpro{\zeta_{6}}$ && $G_{6},\inpro{\zeta_{4}}$  \\
	 $G_{4},\inpro{\zeta_{2}}$ && $\mT,\inpro{\zeta_2}$ && $G_{4},\inpro{\zeta_{2}}$ \\
    };
    \draw (A-1-3)--(A-2-1);
    \draw (A-1-3)--(A-2-3);
    \draw (A-1-3)--(A-2-5);
    \draw (A-2-1)--(A-3-1);
    \draw (A-2-3)--(A-3-1);
    \draw (A-2-3)--(A-3-3);
    \draw (A-2-3)--(A-3-5); 
    \draw (A-2-5)--(A-3-5); 
\end{tikzpicture} 
		& 
		  \begin{tikzpicture}
    \matrix (A) [matrix of nodes, row sep=1.0cm, column sep = -0.2 cm]
			  {  \\ $\mat{G\mapsto\hat G\cr\longrightarrow }$  \\  \\ };
\end{tikzpicture}
		&
	\begin{tikzpicture}
    \matrix (A) [matrix of nodes, row sep=1.0cm, column sep = -0.2 cm]
    {
            $G_{7},\inpro{\zeta_{12}}$ \\
            $G_{5},\inpro{\zeta_{6}}$  \\
		$\mT,\inpro{\zeta_2}$ \\
    };
    \draw (A-1-1)--(A-2-1);
    \draw (A-2-1)--(A-3-1);
\end{tikzpicture} 
       \end{tabular} 
\end{center}
\vskip-0.5truecm
\caption{\label{G7-Zlat} The lattice of the collineation preserving normal reflection 
	subgroups of $G_7$ and their centres, and its image under $G\mapsto\hat G$.  
	The remaining collineation preserving normal subgroups are the 
	reflection free subgroup $\mT$, also depicted, and $\inpro{\zeta_4}\mT$.
	}
\end{figure}
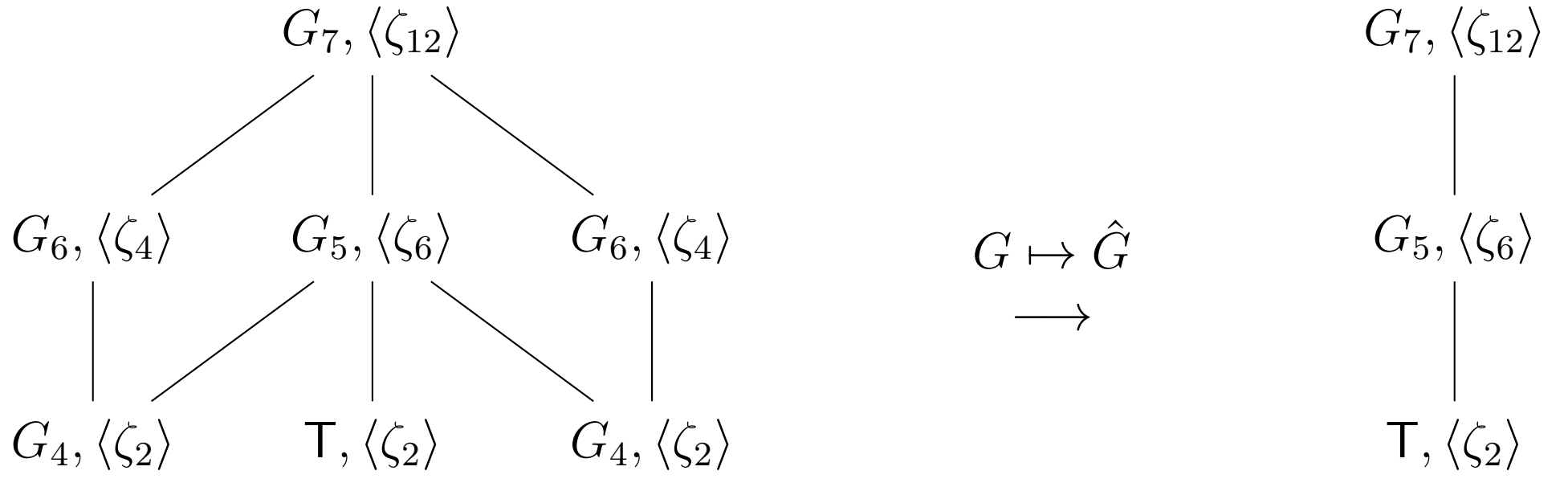

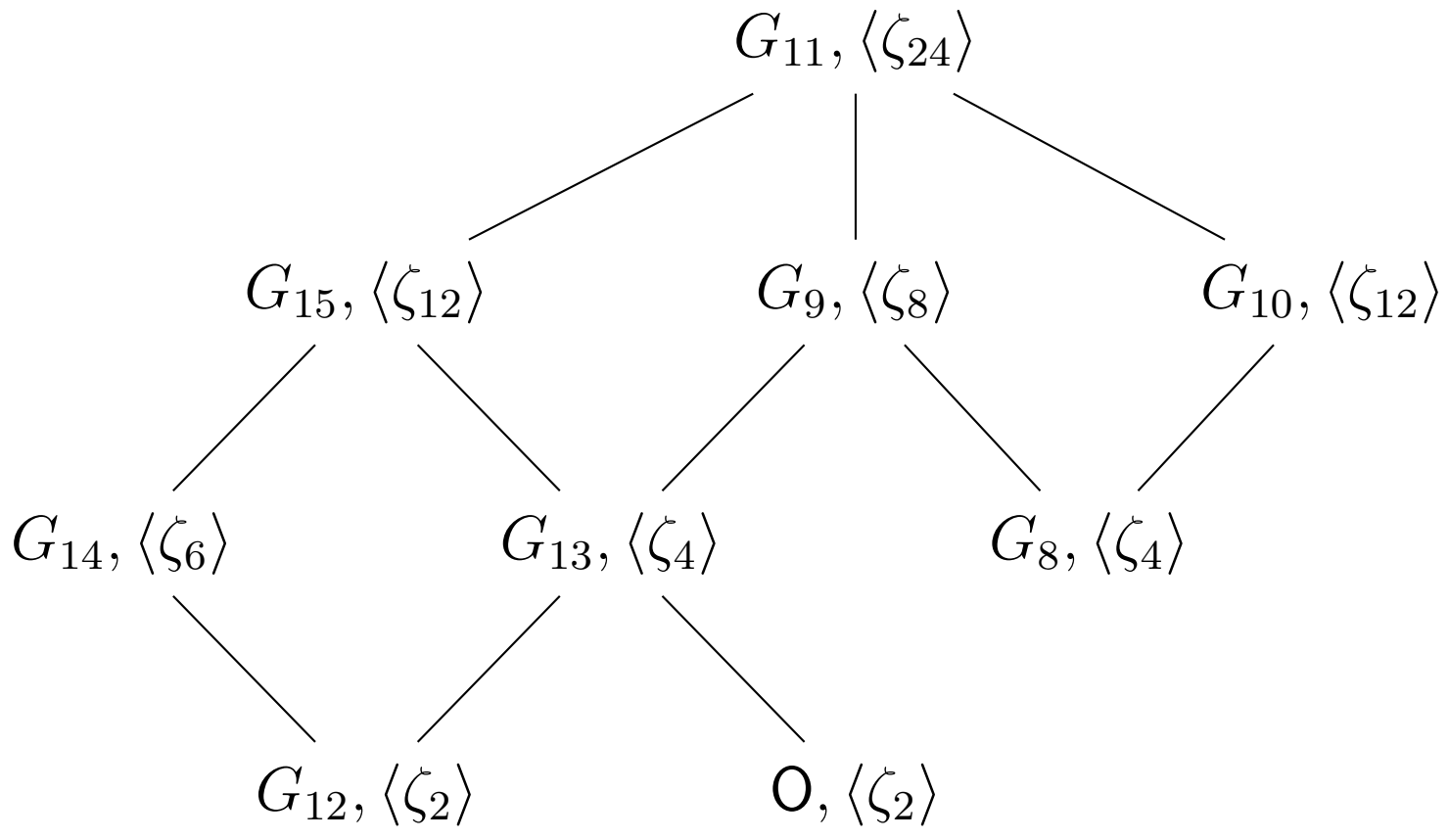
\begin{figure}[h]
\begin{center}
\begin{tikzpicture}
    \matrix (A) [matrix of nodes, row sep=1.0cm, column sep = -0.2 cm]
    {
	    &&& $G_{11},\inpro{\zeta_{24}}$ && \\
	& $G_{15},\inpro{\zeta_{12}}$ && $G_{9},\inpro{\zeta_{8}}$ && $G_{10},\inpro{\zeta_{12}}$ \\
	$G_{14},\inpro{\zeta_{6}}$ && $G_{13},\inpro{\zeta_{4}}$ && $G_{8},\inpro{\zeta_{4}}$ & \\
	& $G_{12},\inpro{\zeta_{2}}$ && $\mO,\inpro{\zeta_2}$&& \\
    };
    \draw (A-1-4)--(A-2-2); \draw (A-1-4)--(A-2-4);
    \draw (A-1-4)--(A-2-6); \draw (A-2-2)--(A-3-1);
    \draw (A-2-2)--(A-3-3); \draw (A-2-4)--(A-3-3);
    \draw (A-2-4)--(A-3-5); \draw (A-2-6)--(A-3-5);
    \draw (A-3-1)--(A-4-2); 
    \draw (A-3-3)--(A-4-2); 
    \draw (A-3-3)--(A-4-4); 
\end{tikzpicture}
\end{center}
\vskip-0.5truecm
\caption{\label{G11-Zlat}
The lattice of the eight collineation preserving normal reflection subgroups of 
	$G_{11}$ and their centres. The remaining normal collineation preserving
	subgroups are the reflection free subgroup $\mO$, also depicted,
	and $\inpro{\zeta_6}\mO$, for which $\mO\subset\inpro{\zeta_6}\mO
	\subset G_{15}$. 
}
\end{figure}

\begin{figure}[H]
\begin{center}
\begin{tikzpicture}
    \matrix (A) [matrix of nodes, row sep=1.0cm, column sep = 0.4 cm]
    {
	    & $G_{19},\inpro{\zeta_{60}}$ & \\
$G_{21},\inpro{\zeta_{12}}$ & $G_{17},\inpro{\zeta_{20}}$ & $G_{18},\inpro{\zeta_{30}}$ \\
$G_{22},\inpro{\zeta_{4}}$ & $G_{20},\inpro{\zeta_{6}}$ & $G_{16},\inpro{\zeta_{10}}$ \\
	&$\mI,\inpro{\zeta_2}$ & \\
    };
    \draw (A-1-2)--(A-2-1);
    \draw (A-1-2)--(A-2-2);
    \draw (A-1-2)--(A-2-3);

    \draw (A-2-1)--(A-3-1);
    \draw (A-2-1)--(A-3-2);
    \draw (A-2-2)--(A-3-1);
    \draw (A-2-2)--(A-3-3);
    \draw (A-2-3)--(A-3-2);
    \draw (A-2-3)--(A-3-3);
    \draw (A-3-1)--(A-4-2); \draw (A-3-2)--(A-4-2); \draw (A-3-3)--(A-4-2); 
\end{tikzpicture}
\end{center}
\vskip-0.5truecm
\caption{\label{G19-Zlat} The lattice of the collineation preserving normal reflection subgroups of $G_{19}$ and their centres,
	and 
the collineation preserving normal subgroup 
	$\mI=G_{22}\cap G_{20}\cap G_{16}$. }
\end{figure}
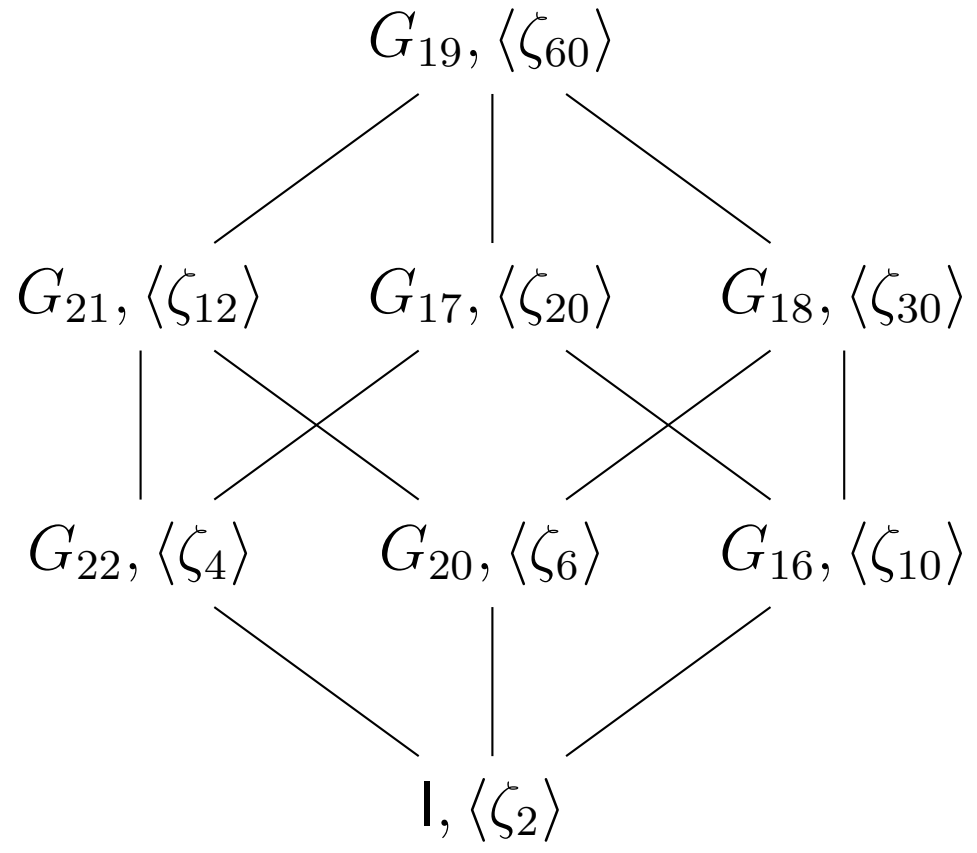

The {\bf collineation group} associated with a matrix group $G$ 
is the group considered up to multiplication by a scalars \cite{B17}, 
i.e., for irreducible groups $G_{\rm coll}=G/Z(G)$.
We recall from \cite{W26c} that each primitive reflection group in $U(\CC^2)$
is a collineation preserving subgroup of the maximal primitive reflection groups
$$ G_7=\inpro{R^2,Z,Z^R}, \qquad
G_{11}=\inpro{F,R,Z}, \qquad G_{19}=\inpro{R^2,Z,M}, $$
as given in Figures \ref{G7-Zlat}, \ref{G11-Zlat}, \ref{G19-Zlat},
where
\begin{equation}
\label{FRZdefn}
F:={1\over\sqrt{2}}\pmat{1&1\cr1&-1}, \qquad
R:=\pmat{1&0\cr0&i^3}, \qquad
Z:= e^{{2\pi i\over 24}} RF,
\end{equation}
\begin{equation}
\label{Mdefn}
M:={e^{-{2\pi i\over 5}}\over2} \pmat{\tau^{-1}-\tau i&-1\cr 1& \tau^{-1}+\tau i},
\qquad
\tau={1+\sqrt{5}\over2},
\end{equation}
are reflections of orders $2,4,3,5$ (respectively).

\begin{lemma} 
\label{hatGformulas}
Let $G\subset U(\CC^2)$ be an irreducible complex reflection group. 
Then $G$ is closed under conjugation by $s$, i.e., $\hat G=G=G^s$,
	equivalently $\det(G)\subset G$, 
in all cases except 
$$ \hat G_4=G_5, \qquad \hat G_6=G_7, $$
$$  \hat G(m,p,2)=G(2m,2p,2), \quad \hbox{$m$ odd}. $$
\end{lemma}

\begin{proof}
In view of (\ref{gsconjrelation}), $G$ is closed under conjugation by $s$
if and only if 
$$ \det(g^*)=\det(g^{-1})=\det(g)^{-1}\in G, \qquad \forall g\in G. $$
Since $g\mapsto\det(g)$ is a homomorphism, it suffices
to check whether (or not) $\det(g)\in G$ for all $g$ from a generating set of 
reflections. If this is not case, then 
$$ \det(g)=gsg^{-1}s \not\in Z(G), $$
for some reflection $g\in G$.


	We first consider the primitive reflection groups, each of which 
is a collineation preserving subgroup of the maximal primitive reflection groups
$$ G_7=\inpro{R^2,Z,Z^R}, \qquad
G_{11}=\inpro{F,R,Z}, \qquad G_{19}=\inpro{R^2,Z,M}, $$
	as given in Figures \ref{G7-Zlat}, \ref{G11-Zlat}, \ref{G19-Zlat}.
It follows from 
$$ \det(F)=-1, \qquad \det(R)=-i, \qquad \det(Z)=\zeta_3, \qquad \det(M)=\zeta_5^{-1}, $$
and the scalar matrices in these subgroups 
	(the centre of these subgroups, see Table \ref{GjTablerank2}), 
that they are all
closed under conjugation by $s$, except for $G_4$ and $G_6$, which each appear twice 
as normal subgroups of $G_7$. i.e.,
	$$ G_4=\inpro{Z,Z^S},\inpro{Z^R,Z^{SR}}, \qquad
	G_6=\inpro{R^2,Z},\inpro{R^2,Z^R}. $$
In this case
	$$  \zeta_3=\det(Z)=\det(Z^R)=ZsZ^{-1}s\not\in Z(G_4),Z(G_6), \qquad Z(G_4)=\inpro{-1},\quad
Z(G_6)=\inpro{i}, $$
and since $\det(R^2)=-1$, we have
$$ \inpro{G_4,s}=\inpro{G_4,\zeta_3,s}=\inpro{G_5,s}=\inpro{\hat G_5,s}, \quad
\inpro{G_6,s}=\inpro{G_6,\zeta_3,s}=\inpro{G_7,s}=\inpro{\hat G_7,s}. $$
Hence $\hat G_4=G_5$, $\hat G_6=G_7$ (in each case).

The calculations for the imprimitive reflection groups
\begin{equation}
\label{G(m,p,2)form}
G(m,p,2) = \inpro{\pmat{0&1\cr1&0},\pmat{0&\zeta_m\cr\overline{\zeta_m}&0},
\pmat{\zeta_{m\over p}&0\cr0&1}}, \qquad Z(G(m,p,2))=
\inpro{\zeta_{\gcd(p,2){m\over p}}},
\end{equation}
are similar. The determinants of the above generators are $-1$, $-1$ and $\zeta_{m\over p}
\in \inpro{\zeta_{\gcd(p,2){m\over p}}}$. Thus $G(m,p,2)$ is closed under 
conjugation except for when $m$ is odd, in which case
\begin{align*}
\inpro{G(m,p,2),s}
&=\inpro{\pmat{0&1\cr1&0},\pmat{0&\zeta_m\cr\overline{\zeta_m}&0},
\pmat{\zeta_{m\over p}&0\cr0&1},-1,s} \cr
&=\inpro{\pmat{0&1\cr1&0},\pmat{0&-\zeta_m\cr-\overline{\zeta_m}&0},
\pmat{\zeta_{m\over p}&0\cr0&1},s} \cr
        &=\inpro{\pmat{0&1\cr1&0},\pmat{0&\zeta_{2m}\cr\overline{\zeta_{2m}}&0},
\pmat{\zeta_{2m\over 2p}&0\cr0&1},s} \cr
        &=\inpro{G(2m,2p,2),s},
\end{align*}
which gives $\hat G(m,p,2)=G(2m,2p,2)$.
\end{proof}

\begin{table}[H]
\caption{
\label{GjTablerank2}
The reflection types and collineation groups of the primitive complex reflection groups of rank $2$ (with generators) as given in 
	\cite{W26c} (also see Tables C.1, C.2 \cite{OT92}).  }
\smallskip
\begin{center}
\begin{tabular}{ |  >{$}l<{$} | >{$}l<{$} | >{$}l<{$} | >{$}l<{$} | >{$}l<{$} | >{$}l<{$} |  >{$}l<{$} | }
\hline
        &&&&&&\\[-0.3cm]
        \hbox{ST} & \hbox{Order} & \hbox{Reflection type} & |Z(G)| & 
        \hbox{Reflections} & G_{\rm coll} & \hbox{Generators} \\[0.1cm]
\hline
  &&&&&&\\[-0.3cm]
G_4 &  24 & 4 C_3 & 2 & 3^{8} & A_4 & Z,Z^S \hbox{ or } Z^R, Z^{SR} \\
G_5 &  72 & 4 C_3, 4 C_3 & 6 & 3^{16} & A_4 & Z, Z^R \\
G_6 &  48 & 6 C_2, 4 C_3 & 4 & 2^{6} 3^{8} & A_4 & R^2, Z \hbox{ or } R^2, Z^R \\
G_7 &  144 & 6 C_2, 4 C_3, 4 C_3 & 12 & 2^6 3^{16} & A_4 & R^2, Z, Z^R \\
  &&&&&&\\[-0.4cm]
G_8 &  96 & 6 C_4 & 4 & 2^6 4^{12} & S_4 & R, R^F \\
G_9 &  192 & 12 C_2, 6 C_4 & 8 & 2^{18} 4^{12} & S_4 & F, R \\
G_{10} & 288 & 6 C_4, 8 C_3 & 12 & 2^{6} 3^{16} 4^{12} & S_4 & R, Z \\
G_{11} & 576 & 12 C_2, 6 C_4, 8 C_3 & 24 & 2^{18} 3^{16} 4^{12} & S_4 & F, R, Z \\
G_{12} & 48 & 12 C_2 & 2 & 2^{12} & S_4 &  F, F^Z, F^{Z^2} \\
G_{13} & 96 & 12 C_2, 6 C_2 & 4 & 2^{18} & S_4 & F, F^Z, R^2 \\
G_{14} & 144 & 12 C_2, 8 C_3 & 6 & 2^{12}3^{16} & S_4 & F, Z \\
G_{15} & 288 & 12 C_2, 6C_2, 8 C_3 & 12 & 2^{18}3^{16} & S_4 & F, R^2, Z \\
        &&&&&&\\[-0.4cm]
G_{16} & 600 & 12C_5 & 10 & 5^{48} & A_5 & M, M^{R^2} \\
G_{17} & 1200 & 30C_2, 12C_5 & 20 & 2^{30}5^{48} & A_5 & R^2, M \\
G_{18} & 1800 & 20C_3, 12C_5 & 30 & 3^{40}5^{48} & A_5 & Z, M \\
G_{19} & 3600 & 30C_2, 20C_3, 12C_5 & 60 & 2^{30}3^{40}5^{48} & A_5 & R^2, Z, M \\
G_{20} & 360 & 20C_3 & 6 & 3^{40} & A_5 & Z, Z^M \\
G_{21} & 720 & 30C_2, 20C_3 & 12 & 2^{30}3^{40} & A_5 & R^2, Z^M \\
G_{22} & 240 & 30C_2 & 4 & 2^{30} & A_5 & R^2, (R^2)^Z, (R^2)^M \\ [0.1cm]

\hline
\end{tabular}
\end{center}
\end{table}

\begin{remark}
In the notation of \cite{C80}, the group $\inpro{G,s}$ for a finite
group of unitary matrices closed under conjugation, i.e., with 
$\hat G=G$ or equivalently $\det(g)\in G$, $\forall g\in G$,
is given by
$$ E(G)=\hbox{the complexification of $\inpro{G,s}$}. $$ 
The $E(H)$ given for $H\subset\GL(\CC^2)$ with
$\det(h)\in H$, $\forall h\in H$, consists of matrices in the
symplectic form (\ref{C2dHdcorres}) if and only if the matrices in $H$ are unitary.
Indeed
$$ E(H):=\inpro{\hat h,[s]_\CC:h\in H}, \qquad
\hat h:=\pmat{h&0\cr 0&(h^T)^{-1}}. $$
To put this in the symplectic form, let $U\in\GL(\CC^2)$ be a 
change of variables for which $G=UHU^{-1}$ is a group of unitary matrices
and  
	$\det(U)=1$ (this preserves the condition $\det(h)\in H$, $\forall h\in H$). Then
$$ \hat U \hat h \hat U^{-1}
=\pmat{UhU^{-1}&0\cr0&(U^T)^{-1}(h^T)^{-1}U^T}=[UhU^{-1}]_\CC,  $$
since
$$ \overline{UhU^{-1}} =((UhU^{-1})^T)^* =((U^{-1})^Th^TU^T)^{-1} 
=(U^T)^{-1}(h^T)^{-1} U^T. $$
Similarly,
	$$  \hat U [s]_\CC \hat U^{-1}
=\pmat{0&-USU^T\cr(U^T)^{-1}SU^{-1}&0}=\pmat{0&-S\cr S&0}=[s]_\CC, 
\qquad S:=\pmat{0&1\cr-1&0}, $$
since
$$ U^TSU=\det(U)S=S, \qquad (U^{-1})^TSU^{-1}=\det(U^{-1})S=\det(U)^{-1}S=S, $$
and we have
$$ \hat U (E(H)) \hat U^{-1}
= \inpro{\hat U \hat H \hat U^{-1},\hat U [s]_\CC \hat U^{-1}}
	= \inpro{[G]_\CC,[s]_\CC}
	=[ \inpro{G,s}]_\CC, $$
	where $G=UHU^{-1}$ is unitary.
\end{remark}

Let $G\subset U(\CC^2)$ be an irreducible finite group. 
Then by Schur's lemma, the group of scalar matrices in $G$ is its centre $Z(G)$.
We may add to 
$G$ the reflection $s$ and a complex scalar matrix $\zeta_m$ of order $m$,
which is equivalent to adding the reflections $s$ and $\zeta s$,
to obtain the following class of finite irreducible quaternionic groups
\begin{equation}
\label{HaGdefn}
\cH_a(G):=\inpro{G,\zeta_{|Z(\hat G)|a},s}, \qquad a\ge 1.
\end{equation}
Here the scalar added has the form $\zeta_{|Z(\hat G)|a}$ since
$\hat G\subset\inpro{G,s}$ already contains $\zeta_{|Z(\hat G)|}$, and so we have
$$ \cH_1(G)=\inpro{G,s}, \qquad
\cH_a(G)=\cH_a(\hat G), $$
with
\begin{equation}
\label{HaGinclusions}
\cH_b(H)\subset\cH_a(G), \qquad b\divides a,\quad H\subset G.
\end{equation}
In view of 
the observation (\ref{generalform}), all of the rank two quaternionic reflection 
groups with imprimitive complexifications have the form (\ref{HaGdefn}).
They have the following structure.


\begin{lemma} 
\label{internaldirectsum}
Let $G\subset U(\CC^2)$ be an irreducible finite group.
Then $\cH_a(G)$ is the (internal) semidirect product
\begin{equation}
\label{HaGsemidirectprod}
\cH_a(G) = \inpro{\hat G,\zeta_{|Z(\hat G)|a}} \rtimes \inpro{s}, 
\qquad |\cH_a(G)|=2|\hat G|a,
\end{equation}
where the multiplication is given by
\begin{equation}
\label{HaGmultiplication}
s^2=1, \qquad sgs=\det(g^{-1}) g, \qquad\forall g\in\inpro{\hat G,\zeta_{|Z(\hat G)|a}}.
\end{equation}
Moreover, $\cH_a(G)$ is a quaternionic reflection group if and only 
if $\inpro{\zeta_{|Z(\hat G)|a}}G$ contains a set of reflections
which generate a collineation preserving subgroup of $G$.
\end{lemma}

\begin{proof}
Since $\hat G$ is closed under conjugation by $s$, and 
$$ s\xi s^{-1}=\overline{\xi}=\xi^{-1}, \qquad\xi\in\CC,\ |\xi|=1, $$
$N=\inpro{\hat G,\zeta_{|Z(\hat G)|a}}$ is a normal reflection
subgroup of 
$\cH_a(G)$. 
Since $N$ consists of complex matrices, $\cH_a(G)=N\inpro{s}$
with $N\cap\inpro{s}=1$, i.e., we have the asserted semidirect product.
Since the scalar matrices in $\hat G$ and $\cH_a(G)$ are $Z(\hat G)$ and $\inpro{\zeta_{|Z(\hat G)|a}}$, we have
$$ |\cH_a(G)|= |\inpro{s}| \cdot |N|
=2 \cdot {|\hat G|\over|Z(\hat G)|} |Z(\hat G)|a=2|\hat G|a. $$

	Since $s^2=1$, the conjugation formula (\ref{gsconjrelation}) 
establishes (\ref{HaGmultiplication})
for $g\in\hat G$, and for $\xi g$, we use this
	and the multiplicativity of the determinant function to obtain (\ref{HaGmultiplication}), 
	i.e.,  
\begin{align*} s(\xi g)s=(s\xi s)(sgs)
&=\pmat{\xi&0\cr0&\xi}^{-1}\det(g^{-1})g \cr
	&=\det\bigl(\pmat{\xi&0\cr0&\xi}^{-1}\bigr) \xi\det(g^{-1})g
=\det((\xi g)^{-1})(\xi g). 
\end{align*}
The largest reflection subgroup contained in $\cH_a(G)$ is 
generated by the reflections in $\inpro{\zeta_{|Z(\hat G)|a}}G$ and $s$, 
and this can equal $\cH_a(G)$ if only if the complex reflections
	generate a collineation preserving subgroup of $G$.
\end{proof}

A finite subgroup $G\subset GL(\CC^2)$ is said to be {\bf minimal}
(for its collineation group) if there is no proper subgroup $H$ of $G$
for which $G=\inpro{\zeta_m}H$ for same scalar $\zeta_m$.
Clearly, all the finite irreducible subgroups of $GL(\CC^2)$ are of 
the form $\inpro{\zeta_m}H$, for some minimal group $H$ (see \cite{NVT08}),
and in particular the groups $\inpro{G,\zeta_m}$ giving rise to 
quaternionic reflection groups $\inpro{G,\zeta_m,s}$ have this form.
The basic mechanics of our development relies on the following fact:
\begin{itemize}
\item For every $g\in U(\CC^2)$ of finite order which is not a scalar matrix,
	two scalar multiples of it give a reflection (a matrix with
		eigenvalue $1$) of finite order.
\end{itemize}
If such a matrix $g$ is not a reflection, then we call it a {\bf hidden reflection}
(see \cite{W26c}), since it may be multiplied by a scalar to obtain/reveal
a reflection. Thus for any finite irreducible group $H$ (minimal or not), 
adding certain scalars to it will give a group of the form $\inpro{G,\zeta}$,
$G$ a reflection group, with the same collineation group as $H$,
i.e., $G/Z(G)=H/Z(H)$. In this way, every quaternionic reflection group
of the type (\ref{HaGdefn})
can be written
$$ \cH_a(G)=\inpro{G,\zeta_{|Z(\hat G)|a},s},  $$
where either
\begin{enumerate}
\item $G$ is a reflection group which is minimal
(has no proper collineation preserving reflection subgroups),
and there is no restriction on $a$.
\item $G$ is a minimal group, and $a$ is such that $\cH_a(G)$ is a reflection group.
\end{enumerate}

We begin with the first, 
our initial approach, ultimately adopting the second in some cases where
it is more convenient. The 
collineation groups $G/Z(G)$, i.e., the irreducible finite subgroups of $PGL(\CC^2)$ were
classified by Klein \cite{K56}. As abstract groups they are
$$ D_n, \quad n\ge 2\quad \hbox{(dihedral group of order $2n$)}, \qquad 
A_4,\qquad S_4, \qquad A_5. $$
A way to give a standardised/canonical group of matrices 
with a given collineation group, which is akin to revealing reflections, is to
take any such matrix group and replace each of its elements by the two scalar multiples
which have determinant $1$, thereby obtaining the finite 
irreducible subgroups of $SU(\CC^2)$ (of twice the order), which we denote by 
$$ \mD_n, \quad n\ge 2\quad \hbox{(binary dihedral group of order $4n$)}, \qquad 
\mT,\qquad \mO, \qquad \mI. $$

We now consider the primitive groups $\mT,\mO,\mI$ which are called the
binary tetrahedral, octahedral and icosahedral groups (of orders $24,48,120$).
In \cite{W26c}, it was observed that any reflection group is a normal subgroup
of the reflection group 
that is obtained
by revealing all its hidden reflections,
the so called {\bf maximal reflection group} for the underlying collineation group.
This is 
the reflection group obtained by revealing all the hidden reflections
of the reflection-free groups $\mD_n,\mT,\mO,\mI$, 
respectively. It is convenient to denote all the groups that can be obtained
from a group $G$ by adding scalars by
\begin{equation}
\label{MaGdefn}
M_a(G):=\inpro{G,\zeta_{|Z(G)|a}}, \qquad a\ge1, \qquad 
M_1(G)=G, \quad|M_a(G)|= |G|\cdot a.
\end{equation}
Each group is different, with
$$ M_b(H)\subset M_a(G),\quad b\divides a, \ H\subset G, \qquad
\cH_a(G)=\inpro{M_a(G),s}, \quad\hbox{when $\hat G=G$}.
$$
The lattices of the collineation preserving normal reflection subgroups of the
maximal reflection groups were described in \cite{W26c} by taking subsets
of the reflection orbits to move down the lattice, i.e., hiding reflections
without changing the collineation groups. Conversely, one can move up the
lattice by adding scalars to reveal reflections. The scalars that must be
added to move to a reflection group above is simply its centre.

We now consider $A_4$ (equivalently $\mT$), which has maximal reflection group $G_7$.
This case is slightly more technical since $\hat G_4=G_5$, $\hat G_6=G_7$.
The minimal reflection subgroups in the reflection subgroup lattice
of Figure \ref{G7-Zlat} are $G_4$ which appears twice. This is also a minimal group.
Applying $G\mapsto\hat G$ to the lattice gives the relevant subgroups
for the groups $\cH_a(G)$. There is one minimal reflection subgroup $G_5$, which 
is not a minimal group, since $\mT$ is a subgroup with $G_5=M_3(\mT)=\inpro{\zeta_6}\mT$.
Starting with $G_5$, we therefore have the quaternionic reflection groups 
\begin{equation}
\label{Ta-family}
\cT_1(G_5):=\cH_a(G_5), \qquad a\ge1.
\end{equation}

We now consider $S_4$ (equivalently $\mO$), which has maximal 
reflection group $G_{11}$. Here there are two minimal
reflection subgroups $G_{12}$ and $G_8$, see Figure \ref{G11-Zlat},
which are also minimal groups.
The corresponding families of quaternionic reflection groups are
\begin{equation}
\label{Oa-family}
\cO_a(G_{12}):=\cH_a(G_{12}), \quad a\ge1, \qquad
 \cO_b(G_{8}):=\cH_b(G_{8}), \quad b\ge1.
\end{equation}
These cross over (first intersect) at
\begin{equation}
\label{OaOb-crossover}
\cO_4(G_{12}) = \inpro{\zeta_8} G_{12}= G_9 = \inpro{\zeta_4}G_8=\cO_2(G_8).
\end{equation}
The presentation in \cite{C80} consists of the following subfamily of $\cO_a(G_{12})$
$$ \cH_{2a}(\mO)=\inpro{\zeta_{4a}}\mO
= \cH_a(G_{13})
= \cH_{2a}(G_{12})
=\inpro{\zeta_{4a}} G_{12}=\cO_{2a}(G_{12}), \qquad a\ge 1, $$
and a family denoted by \cite{C80}, \cite{T25}, \cite{S23} (respectively) as
\begin{align}
\label{OTgrouplabels}
	& E\bigl((\mu_{4d}|\mu_{2d};\mO|\mT)\bigr)&&\quad
	\hbox{($d$ an odd multiple of $1$ or $2$)}, \cr
	& E(\cC_d\boxdot_2\cO)&&\quad
	\hbox{($d$ a multiple of $4$ but not divisible by $16$)},  \cr
	& E(\mO\mT_{2d})&&\quad 
	\hbox{($d$ not divisible by $4$)}. &
\end{align}
The family (\ref{OTgrouplabels})
will be seen to be the reflection groups for which the intersection
of the underlying complex group with $SU(\CC^2)$ is $\mT$ (as opposed to $\mO$).

For $A_5$ (equivalently $\mI$), the maximal
reflection group is $G_{19}$. There are three minimal
reflection subgroups $G_{22}$, $G_{20}$ and $G_{16}$, see Figure \ref{G19-Zlat},
none of which is a minimal group since
\begin{equation}
\label{Igrpsnotminimal}
G_{22}=\inpro{\zeta_4}\mI, \qquad
G_{20}=\inpro{\zeta_6}\mI, \qquad
G_{16}=\inpro{\zeta_{10}}\mI.
\end{equation}
This leads to three infinite (and intersecting) families
\begin{equation}
\label{Iabc-families}
\cH_a(G_{22}), \quad a\ge1, \qquad
\cH_b(G_{20}), \quad b\ge1, \qquad
\cH_c(G_{16}), \quad c\ge1,
\end{equation}
which, in view of (\ref{Igrpsnotminimal}), are most easily described by
$$ \cI_a(\mI):=\cH_a(\mI), \qquad
\hbox{for appropriate choices of $a\ge1$}. $$

Recently,  
it has been shown that a primitive complex reflection group can have
a quaternionic system of imprimitivity, and therefore be 
conjugate to an imprimitive quaternionic reflection group.
The only cases are $G_{12}$, $G_{13}$ and $G_{22}$ 
(see Theorem 4.2 of \cite{W26}
and Theorem 7.3 of \cite{T25}).
It is possible that adding the quaternionic reflection $s$ to these groups might 
give a quaternionic reflection group which shares 
this system of imprimitivity (and hence is imprimitive).
Indeed, this is exactly what happens.

\begin{lemma}
\label{HaGquatimprimitiveid}
Let $G$ be a rank two primitive complex reflection group.
	The quaterionic
	reflection group $\cH_a(G)$ is imprimitive
in precisely the three cases
\begin{align}
\label{HaGnotprimitive}
\cO_1(G_{12})&=\inpro{G_{12},s}=G_\cO(L_{14}^\cO,1), \cr
\cO_2(G_{12})&=\inpro{G_{13},s}=G_\cO(L_{20}^\cO,C_2), \cr
\cI_2(\mI)&=\inpro{G_{22},s}=G_\cI(L_{32}^\cI,C_2),
\end{align}
where the latter notation in terms of reflection systems is from \cite{W25}.
\end{lemma}

\begin{proof}
The quaternionic system of imprimitivity for each of 
$G=G_{12},G_{13},G_{22}$ 
	(as presented here)  is given by (see \cite{W26})
\begin{equation}
\label{Uimpsysform}
U={1\over\sqrt{2}}\pmat{1&-zj\cr zj&-1}, \qquad
z\in\CC, \quad |z|=1, 
\end{equation}
where the normalisation of $U$ ensures that $U^2=1$.
For $\inpro{G,s}$ to remain imprimitive, $U$ must 
conjugate $s$ to a monomial matrix. We calculate
$$ UsU^{-1}={1\over 2} \pmat{z-\overline{z} & -z^2j-j\cr j-z^2j & z-\overline{z}}, $$
which is monomial if and only if $z=\pm1$, and so
$\inpro{G_{12},s},\inpro{G_{13},s},\inpro{G_{22},s}$ are imprimitive.

To see if there are additional imprimitive groups $\cH_a(G)$,
	we consider which complex scalar matrices $\xi$ are conjugated to monomial 
matrices by $U$. We have 
$$ U\xi U^{-1}=\pmat{\xi+\overline{\xi} & \overline{\xi}zj -\xi zj \cr
	\overline{\xi}zj-\xi zj & \xi+\overline{\xi} }, $$
which is monomial if and only if $\overline{\xi}=\pm \xi$, i.e., 
$\xi$ is a $4$-th root of unity.
Since $G_{12}\subset G_{13}$ and $Z(G_{13})=Z(G_{22})=\inpro{\zeta_4}$,
adding any additional complex scalar matrices to the above imprimitive groups
gives a primitive group, and so there are not additional imprimitive groups
	of the form $\cH_a(G)$.
\end{proof}


We are now in a position to describe all the primitive quaternionic 
reflection groups with imprimitive complexifications (without
any groups repeated).

\begin{theorem}
\label{TOIfamilyclass}
(Classification)
The primitive quaternionic reflection groups which have imprimitive complexifications
are of rank one, given by the dicyclic (binary dihedral) groups
\begin{equation}
\label{rankoneimpcomp}
\cD_n=\cH_n(C_{2})=\inpro{\zeta_{2n},j}, \qquad n\ge2,
\end{equation}
or of rank two, and have types $\cT,\cO,\cI$ as follows 
\begin{enumerate}[\rm(i)]
\item $\cT_a(G_5)=\cH_a(G_5)=\inpro{Z,Z^R,\zeta_{6a},s}$, $a\ge 1$.
\item $\cO_a(G_{12})=\cH_a(G_{12})=\inpro{F,F^{Z},F^{Z^2},\zeta_{2a},s}$, $a\ge3$, \\ 
	$\cO_b(G_8)=\cH_b(G_8)=\inpro{R,R^F,\zeta_{4b},s}$, $b\ge1$, $b$ is odd.
\item $\cI_a(\mI)=\cH_a(\mI)=\inpro{iR^2,\zeta_3 Z,\zeta_5^3 M,\zeta_{2a},s}$,
	$a\ge 3$,
	where $2\divides a$, $3\divides a$ or $5\divides a$.
\end{enumerate}		
The structure of these groups is given by Lemma \ref{internaldirectsum},
and their inclusions by (\ref{HaGinclusions}).
\end{theorem}

\begin{proof} 
The rank one complex reflection groups are  
$G=C_n=\inpro{\zeta_n}\subset\CC^*$, $n\ne 1$.
With $s=j$, we have $-1=j^2=s^2\in\cH_a(G)$, so that (\ref{rankoneimpcomp})
gives all the possible groups (not withstanding the fact that
$C_p$, $p$ a prime, are the minimal reflection groups).

We have already given all the possible rank two groups, and so it remains 
to give generators (from Table \ref{GjTablerank2}), 
and to restrict the parameters
in such a way that no group is repeated (up to conjugacy),
and to then exclude the imprimitive groups of (\ref{HaGnotprimitive}).

For $A_4$ there is a single family (\ref{Ta-family}), and we need no
restrictions on $a\ge1$ (the groups all have different orders).
In \cite{C80}, the groups are given as $\inpro{\zeta_{6a}}\mT=\cT_a(G_5)=\cH_a(\mT)$.

For $S_4$ there are two families (\ref{Oa-family}). We will take
	the $\cO_a(G_{12})$ family (of groups with different orders) without any restriction on $a\ge1$ (it could
be done the other way round, but this is simplest and closest to the
presentation in \cite{C80}).
As soon as $b\ge1$ 
becomes even, $\cO_b(G_8)$ becomes equal to an $\cO_a(G_{12})$ group
(otherwise not), with (\ref{OaOb-crossover}) being the first instance.
This is because 
$\cO_b(G_8)$ contains $\inpro{\zeta_8}$, and so 
	(see Figure \ref{G11-Zlat}) we have that
$$ \cO_b(G_{8})
= \inpro{G_8,\zeta_{4b},s}
= \inpro{G_9,\zeta_{4b},s}
	= \inpro{G_{12},\zeta_{2(2b)},s}
	=\cO_{2b}(G_{12}), \qquad
\hbox{$b$ even}. $$
	For $b$ odd, 
	$\cO_b(G_8)$ does not contain $\inpro{\zeta_8}$,
	and hence
	is not an $\cO_a(G_{12})$ group.

	For $A_5$ there are three families (\ref{Iabc-families}). We could impose
conditions on $a,b,c\ge1$ so that a given group appears only once. However, it
is much neater to start with $\mI$ (as does \cite{C80}) which sits below the three
minimal reflection groups $G_{22},G_{20},G_{16}$
	(see (\ref{Igrpsnotminimal}))
and then impose a condition 
on the parameter $a$ of $\cH_a(\mI)$ which insures that it
contains one of these minimal groups, 
	i.e., $2\divides a$, $3\divides a$ or $5\divides a$ (respectively).
Since there is single parameter $a\ge1$, this ensures that no groups are repeated.

Finally, we further restrict the parameters for $\cO_a(G_{12})$
	and $\cI_a(\mI)$ to $a\ge3$ to exclude the imprimitive groups of
	(\ref{HaGnotprimitive}).
\end{proof}

This gives the groups in Theorem 3.6 of \cite{C76} 
(which are seen to be distinct), 
in addition to their structure and inclusions.
In the Theorem 3.6 of \cite{C76} 
the imprimitive groups (\ref{HaGnotprimitive}) are listed as primitive groups.

Our presentation of the $\mO$ groups
differs the most from those developed in (\ref{OTgrouplabels}), which are not
of the $\cH_a(G)$ form. For Taylor's presentation in (\ref{OTgrouplabels}), 
we observe that
$$ \hbox{$d=4a$ but $4a$ not divisible by $16$} \Iff
\hbox{$d=4a$ with $a$ not divisible by $4$}, $$
so that
$$ E\bigl((\mu_{4d}|\mu_{2d};\mO|\mT)\bigr)
 = E(\cC_{4d}\boxdot_2\cO)
       = E(\mO\mT_{2d}) \qquad \hbox{($d$ an odd multiple of $1$ or $2$)}, $$
with these groups appearing in Theorem \ref{TOIfamilyclass} as
(symplectic groups)
\begin{equation}
\label{OTHaform}
E(\mO\mT_{2a})=\cO_a(G_{12}), \quad\hbox{$a$ odd}, \qquad
E(\mO\mT_{4b})=\cO_{b}(G_{8}), \quad\hbox{$b$ odd}.
\end{equation}

\begin{example}
In Section 4.1.1 of \cite{S23}, various inclusions between the ``primitive'' 
quaternionic reflection groups $\cH_a(G)$ of Theorem \ref{TOIfamilyclass} were given. These are instances 
of (\ref{HaGinclusions}), as illustrated in Figure \ref{HaGInclusions}.

The Table 4.2.2 of \cite{S23} gives a list of groups 
	$M_a(G)$ for 
$G=G_0$ a primitive rank two complex reflection group with $\hat G=G$ 
	(i.e., $G\ne G_4,G_6$)
which contain no reflections other than those in $G$.
As an example, the $\mu_4\mO$ and $\mathsf{OT}_6$ entries are
\begin{align*}
G_0 =\mu_4\mO=G_{13}=M_2(G_{12}),
    & \qquad \mu_d\mO=M_{d\over 2}(G_{12}),\quad d\in\{4,20,28\}, \cr
G_0 =\mathsf{OT}_6=G_{14}=M_3(G_{12}),
    & \qquad \mathsf{OT}_d=M_{d\over 2}(G_{12}),\quad d\in\{6,18,30\}. 
\end{align*}
For $G_0=G_{13}$, the groups are 
$\mu_d\mO=M_{d\over4}(G_{13})=M_{d\over2}(G_{12})$, 
where $d=4a$ and $a$ is not divisible 
by $2$ or $3$ (see Figure \ref{G11-Zlat}), i.e.,
$$ \mu_d\mO=M_{d\over 2}(G_{12}), \quad
d\in \{4a:\gcd(a,2)=\gcd(a,3)=1\} 
= \{ 4, 20, 28, 44, 52, 68, 76, 
\ldots\}  $$
and, similarly, for $G_0=G_{14}$ the groups are
$$ \mathsf{OT}_d=M_{d\over 2}(G_{12}), \quad
d\in \{6a:\gcd(a,2)=1\}
=\{ 6, 18, 30, 42, 54, 66, 78, 90, 102, \ldots\}. $$
\end{example}

We finish our investigation for the rank two primitive 
complex reflection groups $G$ by giving the lattice showing their inclusions
(Figure \ref{GjInclusions}),
and the corresponding lattice of quaternionic reflection groups
$\cH_1(G)$ that evolve from them (Figure \ref{HaGInclusions}).

We observe that it is somewhat by chance that the maximal reflection groups
that are obtained from the hidden reflections of 
$\mG=\mT,\mO,\mI$ do contain $\mG$. 
Indeed, $\mT$ and $\mO$ are minimal collineation preserving 
normal subgroups of $G_7$ and $G_{11}$, but they are not below all 
of the collineation preserving normal reflection subgroups.

\vfil\eject

\begin{figure}[!h]
\begin{center}
\begin{tikzpicture}
    \matrix (A) [matrix of nodes, row sep=0.5cm, column sep = 0.2 cm]
        { &&& \greybox{$G_{11}$} &&&& & && \greybox{$G_{19}$} &\\
        & $G_{15}$ && $G_{9}$ && $G_{10}$ &&&& $G_{21}$ & $G_{17}$ & $G_{18}$   \\
        $G_{14}$ && $G_{13}$ && $G_{8}$ && \greybox{$G_7$} &&& $G_{22}$ & $G_{20}$ & $G_{16}$ \\
	& $G_{12}$ && $\mO$ && $G_5$&& $G_6$ && & $\mI$ &  \\
	& && & $\mT$ && $G_4$ && &&& \\
    };
	\draw (A-3-3)--(A-4-4); \draw [blue] (A-4-4)--(A-5-5);  \draw (A-4-6)--(A-5-5);   
    \draw (A-3-10)--(A-4-11); \draw (A-3-11)--(A-4-11); \draw (A-3-12)--(A-4-11); 
    \draw (A-1-4)--(A-2-2); \draw (A-1-4)--(A-2-4); \draw (A-1-4)--(A-2-6);
    \draw (A-2-2)--(A-3-1); \draw (A-2-2)--(A-3-3); \draw (A-2-4)--(A-3-3);
	\draw (A-2-4)--(A-3-5); \draw (A-2-6)--(A-3-5); \draw [blue] (A-2-6)--(A-3-7);
	\draw [blue] (A-2-2)--(A-3-7); \draw (A-3-1)--(A-4-2); \draw (A-3-3)--(A-4-2);
	\draw [blue] (A-3-1)--(A-4-6); \draw (A-3-7)--(A-4-6); \draw (A-3-7)--(A-4-8);
    \draw (A-4-6)--(A-5-7); \draw (A-4-8)--(A-5-7);
        \draw (A-1-11)--(A-2-10);
    \draw (A-1-11)--(A-2-11); \draw (A-1-11)--(A-2-12); \draw (A-2-10)--(A-3-10);
    \draw (A-2-10)--(A-3-11); \draw (A-2-11)--(A-3-10); \draw (A-2-11)--(A-3-12);
        \draw [red,thick,dashed] (A-4-6)--(A-3-11);
        \draw [red,thick,dashed] (A-2-10)--(A-3-7);
        \draw (A-2-12)--(A-3-12);
        \draw (A-2-12)--(A-3-11);
\end{tikzpicture}
\end{center}
\vskip-0.6truecm
\caption{
\label{GjInclusions}
The inclusions between the  rank two primitive complex reflection groups,
	with $\mT,\mO,\mI$ included.
The lines indicate normal reflection subgroups of the 
maximal reflection groups $G_{11},G_7,G_{19}$, with those in blue not being
collineation preserving.
The red dashed lines are inclusions which are neither normal or collineation 
preserving.
}
\end{figure}
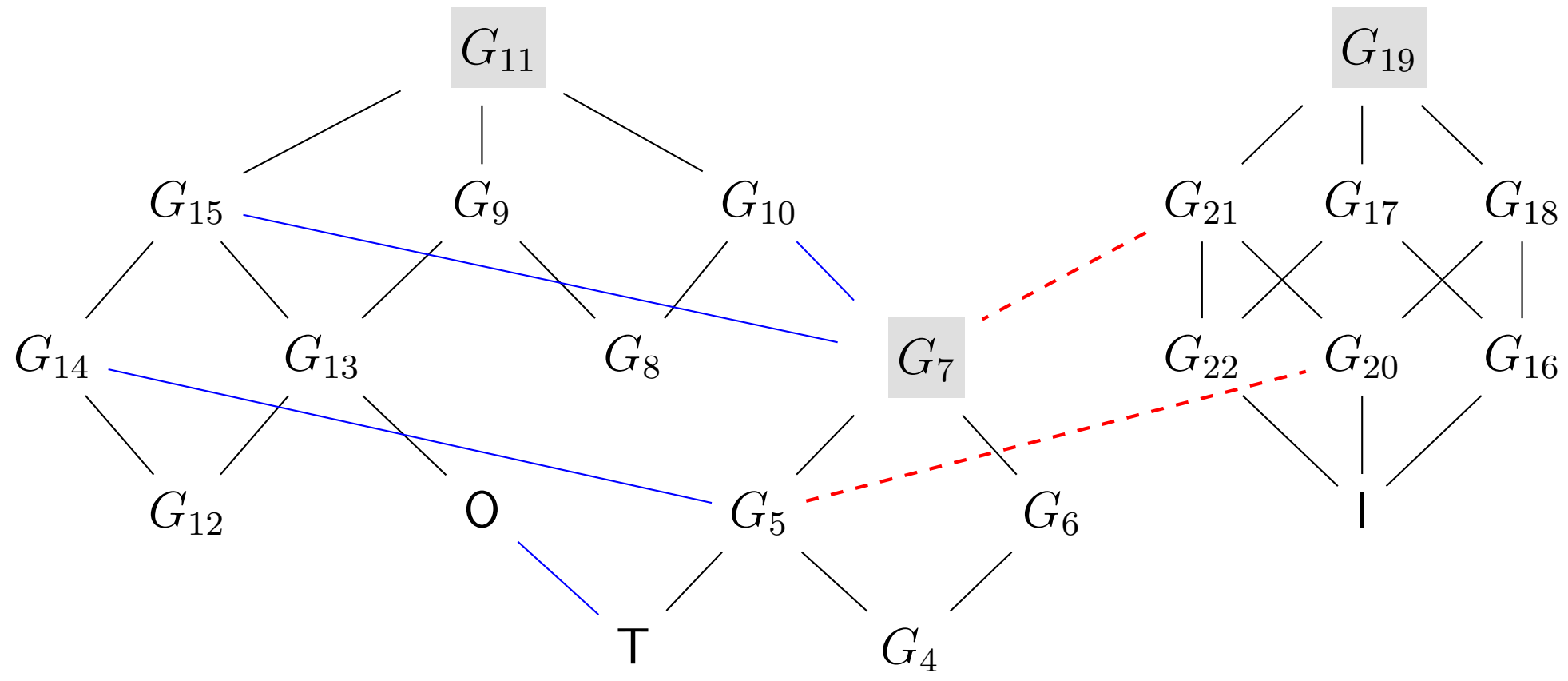

\begin{figure}[!h]
\begin{center}
\begin{tikzpicture}
    \matrix (A) [matrix of nodes, row sep=0.5cm, column sep = -0.6 cm]
	{ &&& $\cO_{12}(G_{12})$ &&&& & && $\cI_{30}(\mI)$ &\\
        & $\cO_{6}(G_{12})$ && $\cO_{4}(G_{12})$ && $\cO_{3}(G_{8})$ &&&& $\cI_{6}(\mI)$ & $\cI_{10}(\mI)$ & $\cI_{15}(\mI)$   \\
	$\cO_{3}(G_{12})$ && \greybox{$\cO_{2}(G_{12})$} && $\cO_{1}(G_{8})$ && $\cT_2(G_5)$ 
	& &\qquad\qquad\quad & 
	\quad\, \greybox{$\cI_{2}(\mI)$} \quad & \, \quad $\cI_{3}(\mI)$ \, \quad & \quad $\cI_{5}(\mI)$\, \quad \\
	& \greybox{$\cO_{1}(G_{12})$} && && $\cT_1(G_5)$&& $\cT_2(G_5)$ && &&  \\
        & && && & $\cT_1(G_5)$ && &&& \\
    };
    \draw (A-1-4)--(A-2-2); \draw (A-1-4)--(A-2-4); \draw (A-1-4)--(A-2-6);
    \draw (A-2-2)--(A-3-1); \draw (A-2-2)--(A-3-3); \draw (A-2-4)--(A-3-3);
    \draw (A-2-4)--(A-3-5); \draw (A-2-6)--(A-3-5); \draw (A-2-6)--(A-3-7);
    \draw  (A-2-2)--(A-3-7); \draw (A-3-1)--(A-4-2); \draw (A-3-3)--(A-4-2);
    \draw  (A-3-1)--(A-4-6); \draw (A-3-7)--(A-4-6); 
	\draw[line width=0.6 mm] (A-3-7)--(A-4-8);
    \draw[line width=0.6 mm] (A-4-6)--(A-5-7); \draw (A-4-8)--(A-5-7);
        \draw (A-1-11)--(A-2-10);
    \draw (A-1-11)--(A-2-11); \draw (A-1-11)--(A-2-12); \draw (A-2-10)--(A-3-10);
    \draw (A-2-10)--(A-3-11); \draw (A-2-11)--(A-3-10); \draw (A-2-11)--(A-3-12);
        \draw [red,thick,dashed] (A-4-6)--(A-3-12);
        \draw [red,thick,dashed] (A-2-10)--(A-3-7);
        \draw (A-2-12)--(A-3-12);
        \draw (A-2-12)--(A-3-11);
\end{tikzpicture}
\end{center}
\vskip-0.6truecm
\caption{
\label{HaGInclusions}
The inclusions between the quaternionic reflection groups
of the form $\cH_1(G)$ corresponding to those of Figure \ref{GjInclusions}
for the primitive complex reflection groups $G$
(using the notation of Theorem \ref{TOIfamilyclass}).
Note that the groups in grey are imprimitive, and 
	the thick black lines indicate groups that are equal. }
\end{figure}
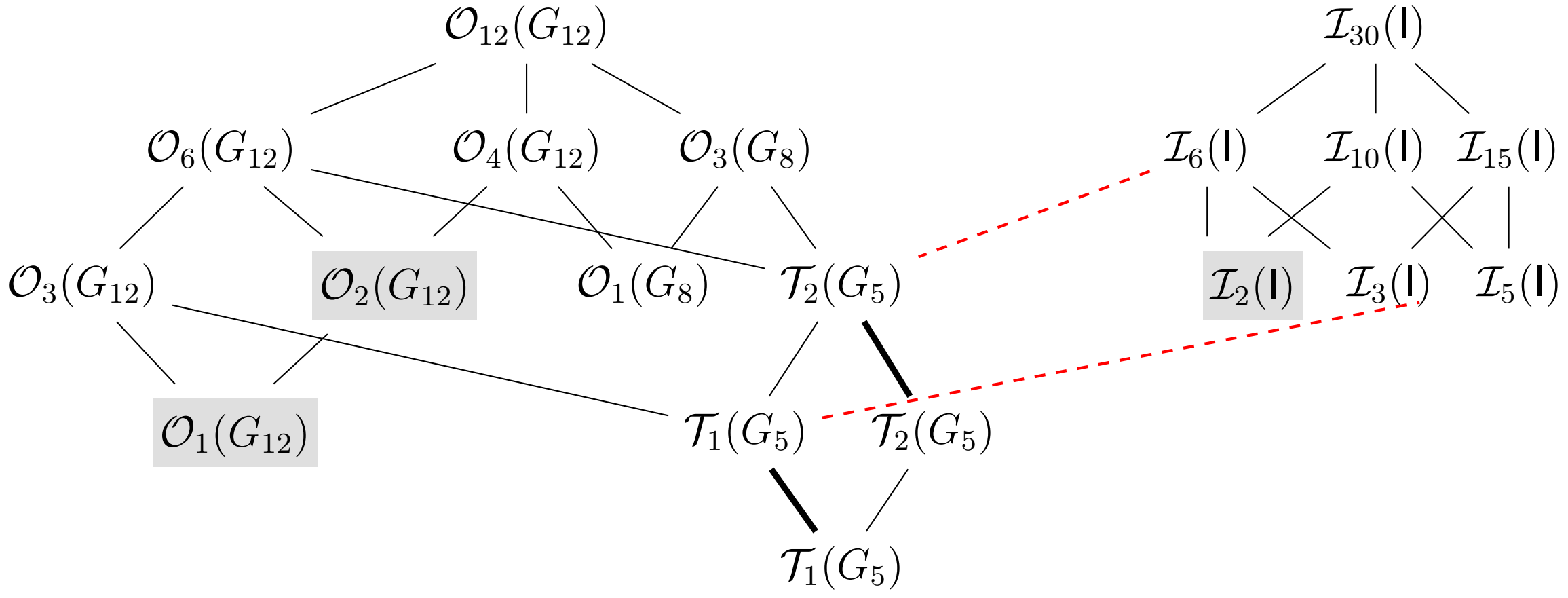


We now consider the quaternionic reflection groups 
obtained from $D_n$ (equivalently $\mD_n$). The corresponding
complex reflection groups are imprimitive, given by
$$ G(m,p,2)=\cG(\{1,\zeta_m\},\{\zeta_m^p\})
=\inpro{\pmat{0&1\cr1&0}, \pmat{0&\zeta_m\cr\zeta_m^{-1}&0}, 
\pmat{\zeta_m^p&0\cr0&1}}, \qquad p\divides m, $$
and so the corresponding quaternionic reflection groups are also imprimitive.
In view of Lemma \ref{hatGformulas}, 
we may restrict ourselves to $m$ even. 

The collineation group of $G(m,p,2)$ has order
\begin{equation}
\label{Gmp2collgrp}
{|G(m,p,2)|\over|Z(G(m,p,2))|} 
={ 2m^2/p\over (m/p)\gcd(p,2)}= {2m\over\gcd(p,2)}
\Implies n={m\over\gcd(p,2)}.
\end{equation}
The groups of interest are therefore:
\begin{equation}
\label{Gampdefn}
G_{a,m,p}:=\cH_a(G(m,p,2))
=\inpro{G(m,p,2),\zeta_{{m\over p}\gcd(p,2)a},s}, \qquad
\hbox{($m$ even)}.
\end{equation}
These have
\begin{equation}
\label{Gamporder}
|G_{a,m,p}|
= 2| G(m,p,2)|a= {4m^2 a\over p}, \qquad 
	|Z(G_{a,m,p})|={m\over p}\gcd(p,2)a,
\end{equation}
and, since the number of reflections in a complex reflection 
group is the sum of its degrees minus its rank, adding these reflections 
to those induced by the complex scalars gives
\begin{equation}
\label{Gampreflections}
\hbox{$\#$ of reflections in $G_{a,m,p}$} = \Bigl( m+{2m\over p}-2\Bigr) + {m\over p}\gcd(p,2)a. 
\end{equation}

To understand the groups $G_{a,m,p}$, 
we now consider the groups $G(n,a,b,r)$ of \cite{W25}, which
consist of monomial matrices (thereby being imprimitive) with entries
in the binary dihedral group $\cD_n:=\inpro{\go,j}$, $\go:=\zeta_{2n}$, $n\ge2$,
and indices $[n,a,b,r]\in\gL_n$, where
\begin{equation}
\label{gLndefn}
\gL_n= \bigcup_{(a,b)\in\gO_n}\{ \hbox{$[n,a,b,{n\over ab}]$}\} \cup
\bigcup_{(a,b)\in\gO_n\atop ab\, {\rm is}\, {\rm odd}} \{\hbox{$[n,a,b,{2n\over ab}]$}\},
\end{equation}
\begin{equation}
\label{Andefn}
\gO_n  := \{(a,b):1\le a\le b\le n,\ a\divides n,\  b\divides n,\ \gcd(a,b)=1\}.
\end{equation}
With the definition of $\cG$ being self evident, we calculate
\begin{align}
\label{G(n,a,b,r)formula}
G(n,a,b,r) & = \cG(\{1,\go^a,j,\go^bj\},\{\go^{2n\over r}\}) \cr
	& = \inpro{ \pmat{0&1\cr1&0},\pmat{0&\go^a\cr\go^{-a}&0},\pmat{0&j\cr-j&0},
	\pmat{0&\go^b j\cr -\go^b j&0},\pmat{\go^{2n\over r}&0\cr0&1} } \cr
& = \inpro{ \pmat{0&1\cr1&0},
	\pmat{0&\go^a\cr\go^{-a}&0},
	\pmat{\go^{2n\over r}&0\cr0&1},
	\pmat{\go^b&0\cr0&\go^b},
	\pmat{0&j\cr-j&0} } \cr
& = \inpro{ \pmat{0&1\cr1&0},
	\pmat{0&\zeta_{2n\over a}\cr\zeta_{2n\over a}^{-1}&0},
	\pmat{(\zeta_{2n\over a})^{2n\over ar}&0\cr0&1},
	\zeta_{2n\over b},s 
	} \cr
& = \inpro{G({2n\over a},{2n\over ar},2),\zeta_{2n\over b},s}, \qquad
	[n,a,b,r]\in\gL_n.
\end{align}
The centre of the complex reflection group above has order
\begin{align*}
|Z(G({2n\over a},{2n\over ar},2))|
&={{2n\over a}\over {2n\over ar}}\gcd(2,{2n\over ar})
= r\gcd(2,{2n\over ar})
=r\begin{cases}
\gcd(2,2b), & r={n\over ab}; \cr
\gcd(2,b), & r={2n\over ab},\, \hbox{$ab$ is odd}; \cr
\end{cases} \cr
&={2n\over ab},
\end{align*}
and since $({2n\over b})/({2n\over ab})=a$, we have
\begin{equation}
\label{G(n,a,b,r)isGamp}
G(n,a,b,r) = \cH_a(G({2n\over a},{2n\over ar},2))
=G_{a,{2n\over a},{2n\over ar}},
\end{equation}
i.e., each $G(n,a,b,r)$ is a $G_{a,m,p}$ group.
As we would expect from previous examples, there may be several
groups $G_{a,m,p}$ that are conjugate to each other.

By (\ref{Gmp2collgrp}), the corresponding ``$n$'' of the collineation group 
$D_n$ for $G({2n\over a},{2n\over ar},2)$ is
\begin{equation}
\label{gcdformula}
{{2n\over a}\over\gcd({2n\over ar},2)}
=\begin{cases}
{n\over a}, & r={n\over ab}; \cr
{2n\over a}, & r={2n\over ab}
\end{cases}
={r b}.
\end{equation}
Now $G(n,a,b,r)$ and $G(n,b,a,r)$ are conjugate to each other 
(the choice $a\le b$ in the indexing $\gL_n$ of (\ref{gLndefn}) was 
made to avoid this repetition). Indeed, with 
$$ U:=\pmat{1&0\cr0&j}, $$
we calculate
\begin{align*}
U G(n,a,b,r) U^{-1}
&= U \cG(\{1,\go^a,j,\go^bj\},\{\go^{2n\over r}\}) U^{-1} \cr
&= \cG(\{-j,-\go^aj,1,\go^b\},\{\go^{2n\over r}\})
= G(n,b,a,r),
\end{align*}
with the last equality holding because
$$ \Bigl\{ \pmat{0&j\cr-j&0}\pmat{0&1\cr1&0}\Bigr\}^2
= \pmat{j&0\cr0&-j}^2=-1 \Implies -1\in G(n,a,b,r), $$
and $U$ is monomial, since $G(n,a,b,r)$ 
has no additional systems of imprimitivity \cite{W26}.
Thus, we have
\begin{equation}
\label{G(n,a,b,r)equivG(n,b,a,r)}
U \cH_a(G({2n\over a},{2n\over ar},2)) U^{-1}
=\cH_b(G({2n\over b},{2n\over br},2)), 
\end{equation}
i.e., for $a\ne b$, $G(n,a,b,r)$ comes from a second complex 
reflection group $G({2n\over b},{2n\over br},2)$ (with a different 
collineation group). We observe that, by (\ref{gcdformula}), 
the numbers of complex reflections and reflections
of the type $\xi s$ in the conjugate groups $G(n,a,b,r)=G_{a,{2n\over a},{2n\over ar}}$
and $G(n,b,a,r)=G_{b,{2n\over b},{2n\over br}}$ are
$$ {2n\over a}+2r-2, \quad {2n\over b}, \qquad\quad
{2n\over b}+2r-2, \quad {2n\over a}, $$
respectively, and the total number of reflections in each group is the same, 
as it must be.
However, conjugation does not simply interchange the reflections of the
two types (as one might expect).
It therefore seems, that deriving the groups $G(n,a,b,r)$ by identifying 
when the groups $G_{a,m,p}$ are conjugate, 
using (\ref{Gamporder}) and (\ref{Gampreflections}), 
may be quite involved.


The classification of \cite{W25} also gave two 
``exceptional'' groups of monomial matrices with entries in $\cD_n$, namely
\begin{align*}
G_{\cD_n}(\cD_n,\cD_n)&=\cG(\{1,\go,j,\go j\},\{\go,j\}), \cr
G_{\cD_n}(\cD_n,\cD_{n/2}) &=\cG(\{1,\go,j,\go j\},\{\go^2,j\}), \qquad
\hbox{($n$ is even)}.
\end{align*}
These have common generators
$$ \pmat{0&1\cr1&0},\pmat{0&\go\cr\overline{\go}&0},\pmat{0&j\cr-j&0},
\pmat{0&\go j\cr-\go j&0}, \quad \pmat{j&0\cr0&1}, $$
equivalently
$$ \pmat{0&1\cr1&0},\pmat{0&\go\cr\overline{\go}&0},
\pmat{\go&0\cr0&\go}, \quad
\pmat{0&j\cr-j&0},
\quad \pmat{j&0\cr0&1}, $$
together with (respectively)
$$ \pmat{\go&0\cr0&1}, \qquad \pmat{\go^2&0\cr0&1} \quad \hbox{($n$ is even)}. $$
The first two generators together with the above reflections (respectively) give
$$ G(2n,1,2), \qquad G(2n,2,2), \quad \hbox{($n$ is even)}. $$
Since the centre of these imprimitive complex reflection groups is $\inpro{\go}$, 
we conclude that
the exceptional groups are obtained by adding two quaternionic 
reflections (the last not having a block monomial symplectic form), i.e.,
\begin{align*}
G_{\cD_n}(\cD_n,\cD_n) &= \inpro{G(2n,1,2),s,\pmat{j&0\cr0&1}}
=\inpro{G(n,1,1,2n),\pmat{j&0\cr0&1}}, \cr
G_{\cD_n}(\cD_n,\cD_{n/2}) &=\inpro{G(2n,2,2),s,\pmat{j&0\cr0&1}}
=\inpro{G(n,1,1,n),\pmat{j&0\cr0&1}},
	\quad \hbox{($n$ is even)}.
\end{align*}

In summary, the groups $G(n,a,b,r)$ have an imprimitive complexification
with the block monomial form (i.e., a system of imprimitivity with two subspaces), whereas the imprimitive quaternionic 
reflection groups 
$G_{\cD_n}(\cD_n,\cD_n)$
and $G_{\cD_n}(\cD_n,\cD_{n/2})$ ($n$ even) have an imprimitive complexification 
(a system of imprimitivity with four subspaces) which is not block monomial.


In Example 5.2 of \cite{W26}, it was shown that the only groups $G(n,a,b,r)$ 
which are conjugate to complex reflection groups are 
$G(n,1,n,1)\cong_\HH G(2n,2,2)$, and accordingly, we 
modify the indexing to $\gL_n^*:=\gL_n\setminus\{[n,1,n,1]\}$.
By combining a number of results and observations, 
we obtain the following classification.

\begin{theorem}
\label{G(n,a,b,r)familyclass}
(Classification)
The rank two imprimitive quaternionic reflection groups which have 
block monomial complexifications, i.e., their complexification
	has a system of imprimitivity of two subspaces
	(in addition to one of four subspaces), are given by 
\begin{equation}
\label{G(n,a,b,r)familyclassI}
G(n,a,b,r) 
= \cH_a(G({2n\over a},{2n\over ar},2))
\cong_\HH \cH_b(G({2n\over b},{2n\over br},2)), \qquad [n,a,b,r]\in\gL_n^*,
\end{equation}
and 
\begin{equation}
\label{G(n,a,b,r)familyclassII}
\cO_1(G_{12})=\inpro{G_{12},s}, \qquad
\cO_2(G_{12})=\inpro{G_{13},s}, \qquad
\cI_2(\mI)=\inpro{G_{22},s}.
\end{equation}
\end{theorem}

\begin{proof}
The classification of the imprimitive quaternionic reflection groups
obtained from $D_n$
in terms of the $G(n,a,b,r)$ was given in \cite{W25} 
(which corrects that of \cite{C80}), and we have discussed 
the exclusion of those which are conjugate to complex reflection groups.
The alternative formulas were derived in (\ref{G(n,a,b,r)formula})
	and (\ref{G(n,a,b,r)equivG(n,b,a,r)}).
The three groups $\cH_a(G)$, with $G$ primitive, which are
imprimitive quaternionic reflection groups
	(with block monomial complexifications) were identified in Lemma \ref{HaGquatimprimitiveid}.
\end{proof}

\section{The finite subgroups of $GL(\CC^2)$}

In this section, we use our techniques to find all the 
minimal finite irreducible subgroups $H\subset GL(\CC^2)$ with a 
given collineation group, from which 
all of the finite irreducible subgroups with this collineation group 
can be obtained explicitly in form of (\ref{MaGdefn}), i.e., 
$$ M_a(H)=\inpro{H,\zeta_{|Z(H)|a}}, \qquad a\ge1, \qquad
M_1(H)=H, \quad |M_a(H)|=|H|\cdot a. $$
We then restrict the parameters so that no group is repeated (up to conjugation),
and identify all of those which are reflection groups. 
We will essentially follow
\cite{NVT08}.

The elements of a general finite irreducible group $H$ with collineation group 
$$ {H\over Z(H)}={\mG\over Z(\mG)}={\mG\over \inpro{-1}}, \qquad
\mG=\mD_n,\mT,\mO,\mI \quad ({\mG\over\inpro{-1}}=D_n,A_4,S_4,A_5), $$
can be written 
as scalar multiples $\gl g$ of the elements $g\in\mG$. 
These scalars form a finite group
\vskip-0.8truecm
\begin{align}
\label{multHdefn}
	\mult(H) &:=\{\gl\in\CC^*:\gl g\in H,g\in\mG\} \cr
&\,\,=\{\pm\sqrt{\det(h)}:h \in H\}, \quad\hbox{when $-1\in Z(H)$}.
\end{align}
It obvious that $\mult(H)$ is a finite group, 
and it contains the centre $Z(H)$ of $H$, since if $\xi\in Z(H)$, then 
$\xi=\xi I$, $I\in\mG$.
The second formula follows from the fact that each $h\in H$
can be written as a scalar multiple of an element of $\mG$ in exactly two possible ways,
i.e.,
$$ h = \pm \sqrt{\det{h}} \Bigl( {h\over\pm\sqrt{\det{h}}}\Bigr), \qquad
{h\over\pm\sqrt{\det{h}}}\in \mG, $$
and it appears in both ways if and only if $-1\in Z(H)$.
Further, $\mult(H)$ is generated by the scalars $\pm\sqrt{\det{h}}$ that
it contains for $h$ in any generating set for $H$.

Understanding $\mult(H)\subset\CC^*$ will help us understand the possible
choices for $H$.
There is a group homomorphism onto $\mult(H)/Z(H)$ 
(essentially the $h$ of \cite{NVT08}) 
given by
\begin{equation}
\label{gThHhomo}
\gTh_H:\mG\to {\mult(H)\over Z(H)}: g\mapsto \gl Z(H), \qquad \gl g\in H.
\end{equation}
Elementary calculations show that $\gTh_H$ is well defined and is a homomorphism. 
Indeed
$$ \gl_1 g,\gl_2 g\in H \Implies \gl_1\gl_2^{-1}\in H\Implies \gl_1\gl_2^{-1}\in Z(H)
\Implies \gl_1 Z(H) =\gl_2 Z(H), $$
$$ \gl_1 g_1, \gl_2 g_2\in H \Implies \gl_1\gl_2 g_1g_2\in H
\Implies \gTh_H(g_1)\gTh_H(g_2)  
=\gl_1\gl_2 Z(H)=\gTh_H(g_1g_2). $$
The kernel of $\gTh_H$ is the set of elements in $H$ with determinant $1$, i.e.,
\begin{equation}
\label{specHeqkergThH}
\spec(H):=\ker(\gTh_H)=\{g\in\mG:\gl g\in H,\forall \gl\in Z(H)\}=H\cap SL(\CC^2),
\end{equation}
and
\begin{equation}
\label{multHquotbyZH}
{\mult(H)\over Z(H)} \cong {\mG\over\spec(H)}.
\end{equation}
Since the image of $\gTh_H$ of (\ref{gThHhomo}) is abelian, 
$\gTh_H$ must factor through $\mG/[\mG,\mG]$, which is given by
\begin{equation}
\label{factorquots}
{\mD_n\over[\mD_n,\mD_n]}={\mD_n\over\mC_n}=C_4, \quad
{\mT\over[\mT,\mT]}={\mG\over\mD_2}=C_3, \quad 
{\mO\over[\mO,\mO]}={\mO\over\mT}=C_2, \quad 
{\mI\over[\mI,\mI]}=1.
\end{equation}
Thus $\mult(H)/Z(H)$ is (isomorphic to) a subgroup of $\mG/[\mG,\mG]$,
and we do have the implied inclusion
\begin{equation}
\label{specHinclusion}
[\mG,\mG]\subset\spec(H)\subset \mG.
\end{equation}
If  $\mult(H)/Z(H)=1$, i.e., $\mult(H)=Z(H)$, then 
\begin{equation}
\label{trivialspecH}
H=M_a(\mG)=\inpro{\zeta_{2a}}\mG, \qquad a\ge1,
\end{equation}
with $\mG =M_1(\mG) =\inpro{\zeta_2}\mG$ being minimal. 
This accounts for all groups $H$ for $\mG=\mI$. 

Before treating the cases where the $\mult(H)/Z(H)$ is not trivial,
we observe that generators for any $\mG$ can be obtained by taking 
any set of reflections which generate a reflection group with the same collineation
group (see Table \ref{GjTablerank2}), 
and then multiplying each of them by a scalar so that it lies in $\mG$, i.e.,
has determinant $1$. In this way, we will choose a generating set for $\mG$ whose
image under $\gTh_H$ has the fewest number of nonidentity elements as is possible.
We observe that the reflections
$F,R,Z,M$ of (\ref{FRZdefn}) and (\ref{Mdefn}) can be 
``hidden'' in $SU(\CC^2)$ (which has no reflections) as $iF$ and
\begin{equation}
\label{RZMhide}
\sqrt{i}R=\pmat{\sqrt{i}&0\cr0&\sqrt{i}^{-1}}, \quad \zeta_3 Z={1\over2}\pmat{-1+i&-1+i\cr1+i&-1-i}, 
\quad \zeta_5 M={1\over2} \pmat{\tau^{-1}-\tau i&-1\cr 1& \tau^{-1}+\tau i},
\end{equation}
and that we have
$$ \mT=\inpro{iR^2,\zeta_3 Z} \subset \mO=\inpro{\sqrt{i}R,\zeta_3 Z}
=\inpro{iF,\zeta_3 Z}. $$


For $G=\mT=\inpro{iR^2,\zeta_3 Z}$, there could be (and are) groups with
$|\mult(H)|=3|Z(H)|$. To see what these are, we suppose, 
without loss of generality, that
$$ |Z(H)|=3^r\cdot 2k, \quad \gcd(3,k)=1, \qquad r\ge0. $$
To have $|\mult(H)|=3\cdot 3^r\cdot 2k$, the group 
$H$ must have an element of
$\mT$ of order divisible by $3$, multiplied by $\zeta_{3^{r+1}}$
(multiplying an element of $\mG$ of order not divisible by $3$
would introduce the scalar $\zeta_{3^{r+1}}$ to the centre of $H$).
The only such elements are two conjugacy classes of elements of order three
given by $\zeta_3 Z$ and $(\zeta_3 Z)^2$, and two of elements of order six
given by $-\zeta_3 Z$ and $-(\zeta_3 Z)^2$. 
Since $-1\in\mT$ and 
$\mT$ is a normal subgroup of $\mO$ with $\zeta_3 Z$ and $(\zeta_3 Z)^2$ 
conjugate in $\mO$, 
it suffices to choose any of the four elements above, and we take $\zeta_3 Z$.
For $H$ to have collineation group $A_4$, it must contain
a scalar multiple of $iR^2\in\mT$, and the scalar multiplier 
cannot have an order divisible by $3^{r+1}$ (which would introduce 
$\zeta_{3^{r+1}}$ to the centre of $H$), and so this scalar belongs to $Z(H)$.
Therefore, such an $H$ contains the group
$$ \inpro{\zeta_{3^r\cdot 2k},iR^2,\zeta_{3^{r+1}}(\zeta_3 Z)}
= \inpro{\zeta_{2k}} \inpro{iR^2,\zeta_{3^{r+1}}(\zeta_3 Z)}, 
\quad\hbox{(since $(\zeta_{3^{r+1}}(\zeta_3 Z))^3=\zeta_{3^{r}}$)}.  $$
Thus the minimal groups with $\mult(H)/Z(H)=C_3$ are
\begin{equation}
\label{HrA4defn}
H_r^{A_4}:=\inpro{iR^2,\zeta_{3^{r+1}}(\zeta_3 Z)}, 
\qquad |H_r^{A_4}|=3^r|\mT|=3^r\cdot 24, 
\qquad r\ge0.
\end{equation}
Of these, only one 
could be a reflection group based on there
being a minimal collineation preserving normal reflection subgroup of $G_7$
with the same order and centre, which it is, i.e.,
\begin{equation}
\label{A4refgroup}
H_0^{A_4} =\inpro{iR^2,\zeta_3^2 Z} =\inpro{\zeta_3^2 Z,(iR^2)\zeta_3^2 Z (iR^2)^{-1}} 
=G_4\ne \mT.
\end{equation}

For $\mG=\mO=\inpro{iF,\zeta_3 Z}=\inpro{\sqrt{i}R,\zeta_3 Z}$, there are groups with $|\mult(H)|=2|Z(H)|$. 
These have 
$$ |Z(H)|=2^r\cdot 2k = 2^{r+1}k, \quad \gcd(2,k)=1, \qquad r\ge0. $$
Since $|\mO|=2^4\cdot 3$, the previous reasoning implies that $H$ must contain
$\zeta_{2^r\cdot 2k}$ and $\zeta_3 Z$.
To have an element of order $2^{r+2}$ in $\mult(H)$, 
$H$ must also contain an element $\zeta_{2^{r+2}}g$, $g\in\mO$,
where order of $g$ is a power of $2$ and
$\inpro{g,\zeta_3 Z}=\mO$. The only possibilities for $g$ up to conjugation in $\mO$
are $iF$ (order $4$) and $\pm\sqrt{i}R$ (order $8$).
We first consider $iF$, i.e., that $H$ contains
$$ H_{r,k}:=\inpro{\zeta_{2^r\cdot 2k},\zeta_{2^{r+2}}(iF),\zeta_3 Z}. $$
For $r\ge1$, we have
$$ (\zeta_{2^{r+2}}(iF))^{2+2^{r+1}} 
=\zeta_{2^{r+1}} \Implies H_{r,k}=\inpro{\zeta_{k}} \inpro{\zeta_{2^{r+2}}(iF),\zeta_3 Z}, \quad r\ge1. $$
The above formula for $H_{r,k}$ also holds for $r=0$, since in this case
$$ \zeta_{2^{r+1}}=-1\in\inpro{\zeta_{2^{r+2}}(iF),\zeta_3 Z}
= \inpro{-F,\zeta_3 Z}. $$
The relation
$$ 
\zeta_3 Z = (\sqrt{i}R)(iF), 
$$
implies that 
$$ \inpro{\zeta_{2^{r+2}}(iF),\zeta_3 Z}=\inpro{\zeta_{2^{r+2}}(\sqrt{i}R),\zeta_3 Z}. $$
Thus the minimal groups with $\mult(H)/Z(H)=C_2$ are
\begin{align}
\label{HrS4grps}
H_r^{S_4}
&:=\inpro{\zeta_{2^{r+2}}(iF),\zeta_3 Z} \cr 
&\,\,=\inpro{\zeta_{2^{r+2}}(\sqrt{i}R),\zeta_3 Z}, \qquad
|H_r^{S_4}|=2^r\cdot|\mO|=2^r\cdot 48,  \qquad r\ge0.
\end{align}
By examination of the minimal reflection groups for $S_4$, we have
\begin{equation}
\label{S4refgroups}
H_0^{S_4} =\inpro{{i}(iF),\zeta_3 Z}
=G_{12}\ne\mO,  \quad |H_0^{S_4}|=|\mO|, \qquad
H_1^{S_4}=\inpro{i\sqrt{i}F,\zeta_3 Z} =G_8.
\end{equation}

\begin{remark}
Our group $\mult(H)$ of (\ref{multHdefn}) is in the same spirit 
as the Schur multiplier, but it is not same concept.
Indeed, the Schur multiplier group of $A_4,S_4,A_5$ and $D_n$ ($n$ even)
is $C_2$, and it is trivial for $D_n$ ($n$ odd).
We can recognise the Schur covering groups as minimal subgroups, e.g., 
$\mT$ and $G_4$ are the unique Schur cover of $A_4$,
$\mO$ and $G_{12}$ nonisomorphic Schur covers of $S_4$,
and $\mI$ is the unique Schur cover of $A_5$.
\end{remark}

We can now prove the following.

\begin{theorem} 
\label{GL2Csubgroups-primitive} (Classification)
The minimal finite primitive subgroups of $GL(\CC^2)$ for the
possible collineation groups $A_4,S_4,A_5$ are (up to conjugation)
\begin{enumerate}[\rm (i)]
\item $A_4$:  \quad $\mT=\inpro{iR^2,\zeta_3 Z}$ and 
$H_r^{A_4}=\inpro{iR^2,\zeta_{3^{r+1}}(\zeta_3 Z)}$, $r\ge0$.
\item $S_4$:  \quad $\mO=\inpro{\sqrt{i}R,\zeta_3 Z}$ and
$H_r^{S_4} 
=\inpro{\zeta_{2^{r+2}}(\sqrt{i}R),\zeta_3 Z}$, $r\ge0$.
\item $A_5$:  \quad $\mI=\inpro{iR^2,\zeta_5 M}$.
\end{enumerate}
with the generators given by (\ref{RZMhide}).
Exactly three of these are reflection groups, namely
\begin{equation}
\label{A4S4refgroups}
H_0^{A_4}={G_4}, \qquad H_0^{S_4}=G_{12}, \qquad H_1^{S_4}=G_{8}.
\end{equation}
Basic properties of these minimal groups are collected in Table 
\ref{GL(2,C)-minimalgroups-table}.
\end{theorem}

\begin{proof}
We have established the possible minimal subgroups 
(\ref{HrA4defn}), (\ref{HrS4grps}), in addition to the minimal groups $\mT,\mO,\mI$.
It remains to show that these are indeed minimal.  


A collineation preserving subgroup $H$ of $H_r^{A_4}$ must contain
a scalar multiple of $iR^2$ and one of $\zeta_3 Z$. 
Since $iR^2$ has order $2$, and $-1$ is commutator in $\mT$, and hence 
is a commutator in $H$, we have (after taking appropriate powers) that
$iR^2\in H$.
The only scalar multiples of $\zeta_3 Z$ 
in $H_r^{A_4}$ are those by the scalars in $\zeta_{3^{r+1}}\inpro{\zeta_{3^r}}$,
and so $\zeta_{3^{r+1}}\gl\in H$, for $\gl\in \inpro{\zeta_{3^r}}$.
	Thus $H$ contains $\inpro{iR^2,\zeta_{3^{r+1}}\gl(\zeta_3)Z}=H_r^{A_4}$, 
	and $H_r^{A_4}$ is minimal.

Now suppose that $H$ is a collineation preserving subgroup of $H_r^{S_4}$,
and so contains scalar multiples of $\sqrt{i}R$ and $\zeta_3 Z$.
The scalar multiples of $\sqrt{i}R$ in $H_r^{S_4}$ are those 
by scalars $\zeta_{2^{r+2}}\gl$, $\gl\in Z(H_r^{S_4})=\inpro{\zeta_{2^{r+1}}}$,
and $\zeta_3 Z$ appears as $\gl(\zeta_3 Z)$ (a power of which equals $\zeta_3 Z$). 
Thus $H$ contains $\inpro{\zeta_{2^{r+2}}\gl\sqrt{i}R,\zeta_3 Z}
	=\inpro{\zeta_{2^{r+2}}\sqrt{i}R,\zeta_3 Z}=H_r^{S_4}$, 
and so $H_r^{S_4}$ is minimal.

The minimal groups which are reflection groups are given in (\ref{A4refgroup})
	and (\ref{HrS4grps}).
\end{proof}

In \cite{NVT08}, it is claimed that the groups $\mO$ and $H_0^{S_4}=G_{12}$ of 
order $48$ are conjugate, and so there is a missing group in (4) (c) of their Theorem 4.

We now consider the dihedral groups $D_n$, $n\ge 2$, for which 
the binary dihedral group is monomial (imprimitive), and given by
\begin{equation}
\label{BinDihedralgens}
\mD_n=\inpro{A,B}, \qquad A 
=\pmat{\zeta_{2n}&0\cr0&\zeta_{2n}^{-1}}, \quad 
B=\pmat{0&1\cr-1&0}, \quad \Bigl\{\hbox{or }\ B=\pmat{0&i\cr i&0}\Bigr\}.
\end{equation}
For $H$ a minimal group with collineation group $D_n$, $\mult(H)/Z(H)$ is a
subgroup of $C_4$, with the trivial subgroup giving $H=\mD_n$.
We now suppose that 
$\mult(H)/Z(H)$ is not trivial.  
Since the collineation
group of $H$ is $D_n$, it must contain a subgroup
$$ H_{\gl,\mu}:=\inpro{\gl A,\mu B}, \qquad \gl,\mu\in\mult(H). $$
We observe that $\mD_n$ is also generated by $B$ and $AB$, which have
the same algebraic structure within the group.
In particular, $B^2=(AB)^2=-1$, so that
$$ (\mu B)^2=-\mu^2, \quad (\gl A\mu B)^2 =-\gl^2\mu^2
\Implies \inpro{\gl^2,-\mu^2}\subset Z(H)\subset\inpro{\gl,\mu}=\mult(H_{\gl,\mu}).$$

By the symmetry in the generators $B$ and $AB$, we 
may assume, without loss of generality, that $\mu$ is a 
$2^{r+2}\cdot k$ primitive root of unity, 
where $\gcd(2,k)=1$, $r\ge-1$.
Now
$$ (\mu B)^{2^{r+2}} =\hat\zeta_k, \quad r\ge0, \qquad
(\mu B)^{2(k+1)}=(\zeta_{2k}B)^{2(k+1)}=\hat\zeta_k, \quad r=-1, $$
with $\hat\zeta_k$ a primitive $k$ root of unity,
so that 
$$ \inpro{\zeta_k}\inpro{\gl A,\mu^k B}\subset\inpro{\gl A,\mu B} \subset H. $$
By the minimality of $H$, we must have $k=1$,
and $\mu$ is a primitive $2^{r+2}$ root of unity.
By the same reasoning for $AB$, $\gl\mu$ is is a primitive $2^{r+2}$ root of unity,
and hence $\gl=(\gl\mu)\mu^{-1}$ is a $2^{s}$ root of unity,
where $0\le s \le r+2$. 
For $r=-1$, we have $\mu=-1$, $\gl=\pm1$, so that
$H_{\gl,\mu}
=\inpro{\pm A,-B}=\mD_n$, 
since $Z(\mD_n)=\inpro{-1}$,
and so we must have $r\ge 0$.

The $2^{r+2}$ primitive roots of unity are obtained by multiplying
any given one $\mu$ by the $2^{r+1}$ roots of unity. Now $-\mu^2\in Z(H)$,
and so $\zeta_{2^{r+1}}$ must be in the centre of $H$,
except for possibly when $r=0$ (i.e., $\mu=\pm i$).
Hence for $r\ne0$, the choices for $\gl$ and $\mu$ can be 
reduced to $\mu=\zeta_{2^{r+2}}$ (i.e., any primitive root can be taken)
and $\gl=1$ or $\gl=\zeta_{2^{r+2}}$. 
Accordingly, we define groups
\begin{equation}
\label{HrKrDn-defn}
H_r^{D_n}:=\inpro{A,\zeta_{2^{r+2}}B}, \qquad
K_r^{D_n}:=\inpro{\zeta_{2^{r+2}}A,\zeta_{2^{r+2}}B}, \qquad r\ge 0,
\end{equation}
where for $r=0$ there may be other groups.
We note that from our discussion
\begin{equation}
\label{HrKrDn-possmin}
Z(H)=\inpro{\zeta_{2^{r+1}}}, \quad \mult(H)=\inpro{\zeta_{2^{r+2}}}, 
\qquad H=H_r^{D_n},K_r^{D_n}, \qquad r\ge1, 
\end{equation}
and that these groups will be seen to be minimal. 

We now consider all the possible groups for $r=0$, i.e.,
$\inpro{\gl A,\mu B}$, $\gl,\mu\in\{1,-1,i,-i\}$.
We have the conjugation
$$ P=\pmat{1&0\cr0&-1}: \qquad
PAP^{-1}=A, \quad PBP^{-1} = -B \quad\hbox{(for either choice of $B$)}, $$
and so the following pairs of groups are conjugate
$$ H_0^{D_n}=\inpro{A,iB}, \inpro{A,-iB}, \qquad
K_0^{D_n}=\inpro{iA,iB}, \inpro{iA,-iB}, $$
$$
\inpro{-A,iB}, \inpro{-A,-iB}, \qquad
\inpro{-iA,iB}, \inpro{-iA,-iB}.  $$
They are in fact all equal, except for $\inpro{-A,iB}$ and $\inpro{-A,-iB}$
when $n$ is odd.

We now suppose that $n$ is odd. 
Since $(iB)^2=1$, the $2n$ elements of $H=\inpro{-A,iB}$ are
$$ \{1,-A,A^2,-A^3,\ldots,-A^{2n-1}\}\{1,iB\}, $$
so that $Z(H)=1$, $\mult(H)=\inpro{i}$, i.e., $\mult(H)/Z(H)=C_4$.
Since $Z(H)=1$, this group is minimal and equal to $D_n$. Moreover, (by considering orders) we have
$$ \inpro{-1}D_n=\inpro{-1}\inpro{-A,iB}=\inpro{A,iB}=H_0^{D_n}, 
\qquad \inpro{i}D_n=\inpro{i}\inpro{-A,iB}=K_0^{D_n}, 
\inpro{-iA,iB},
$$
and so there is just one minimal group $D_n=\inpro{-A,iB}$,
which is a real reflection group.

We now suppose that $n$ is even. Then
\begin{align*}
(\pm A)^n=-1 
&\Implies H_0^{D_n} =\inpro{A,iB} =\inpro{-A,iB} =\inpro{A,-iB} =\inpro{-A,-iB},\cr
(\pm i B)^2=-1
&\Implies K_0^{D_n} =\inpro{iA,iB} =\inpro{-iA,iB} =\inpro{iA,-iB} =\inpro{-iA,-iB},
\end{align*}
and we will see that $H_0^{D_n}$ and $K_0^{D_n}$ are minimal groups.

\begin{theorem}
\label{GL2Csubgroups-imprimitive} (Classification)
The minimal finite imprimitive irreducible subgroups of $GL(\CC^2)$ for the
possible collineation groups $D_n$, $n\ge2$, are (up to conjugation)
\begin{enumerate}[\rm (i)]
	\item $\mD_n=\inpro{A,B}$ (binary dihedral group).
	\item $D_n=\inpro{-A,iB}$ (dihedral group) and $H_r^{D_n}=\inpro{A,\zeta_{2^{r+2}}B}$, $r\ge1$, 
	for $n$ odd.
\item $H_r^{D_n}=\inpro{A,\zeta_{2^{r+2}}B}$ 
	and $K_r^{D_n}=\inpro{\zeta_{2^{r+2}}A,\zeta_{2^{r+2}}B}$, $r\ge0$, for $n$ even.
\end{enumerate}
These groups are distinct, except for when $n=2$,
$H_r^{D_2}$ and $K_r^{D_2}$ are conjugate.
Of these groups, the following are reflection groups
\begin{align}
\label{MinimalDn-refgroups}
D_n &=G(n,n,2),  &&\hbox{$n$ odd}, \cr
K_r^{D_n} &=G(n,b,2),  && \hbox{$n$ even}, \ n=2^{r+1}b, \ \gcd(2,b)=1, \cr
H_0^{D_n} &=G(2n,2n,2), && \hbox{$n$ even}.
\end{align}
Basic properties of these minimal groups are collected in Table
\ref{GL(2,C)-minimalgroups-table}.
\end{theorem}

\begin{proof}
We have shown that $\mD_n$ and $D_n$ (for $n$ odd) are minimal groups.
It therefore remains to determine whether the additional 
	possible minimal groups of (\ref{HrKrDn-possmin}) and $H_0^{D_n}$ and $K_0^{D_n}$
	(for $n$ even) are minimal, and 
	whether there are any conjugacies between them.

For $r=0$ and $n$ even, we have $A^n=(\zeta_{2^{r+2}}A)^n=-1$, so that the
centre of $H_r^{D_n}$ and $K_r^{D_n}$ is $\inpro{\zeta_{2^{r+1}}}$,
in this case, as well as for $r\ge1$. Now suppose that
$r\ge1$ or $r=0$ and $n$ is even, and let
	$H$ be a collineation preserving subgroup of $H_r^{D_n}$
or $K_r^{D_n}$. Then $H$ must contain a scalar multiple of the generators, 
i.e.,
$$ \gl A,\quad \gl\zeta_{2^{r+2}}A\quad \hbox{(respectively)},\qquad
\mu \zeta_{2^{r+2}}B, \qquad \gl,\mu\in\inpro{\zeta_{2^{r+1}}}. $$
Since $(\mu\zeta_{2^{r+2}}B)^2=-\mu\zeta_{2^{r+1}}$ generates 
$\inpro{\zeta_{2^{r+1}}}$, for $r\ge1$, and
$(\gl A)^n=(\gl\zeta_{2^{r+2}}A)^n=-1$ generates $\inpro{\zeta_{2^{r+1}}}$, 
for $r=0$ and $n$ even,
these elements generate the corresponding groups $H_r^{D_n}$
and $K_r^{D_n}$ (respectively), which are therefore minimal.

For $n$ odd and $r\ge1$, we have
$$ (\zeta_{2^{r+2}}A)^n = \zeta_{2^{r+2}}^n(-1), $$
which is a primitive $2^{r+2}$ root of unity,
so that $H_r^{D_n}$ is contained within $K_r^{D_n}$,
and, since they are both minimal groups, we have
$H_r^{D_n}=K_r^{D_n}$, which gives (ii).

For $n$ even and $r\ge0$, we have (by examination 
and the fact $|H|=2|\spec(H)|$) that
$$ \spec(H_{r}^{D_n})=\inpro{A}=\mC_{2n}, \qquad
	\spec(K_{r}^{D_n})=\inpro{A^2,AB}\cong\mD_{n/2}. $$
Since the determinant of a matrix is preserved under conjugation and
	$\mD_n$ is not cyclic for $n\ne2$, 
we conclude that $H_{r}^{D_n}$ and $K_{r}^{D_n}$ are not conjugate
for $n\ne 2$. For $n=2$, $A^2=-1$, and we have
	$\inpro{A}\cong C_4\cong\inpro{A^2,AB}=\inpro{AB}$.
It is easy to verify that $H_{r}^{D_2}$ and $K_{r}^{D_2}$ are conjugate.
If we take
$$ P=\pmat{1&1\cr1&-1}, \quad \Bigl\{\hbox{or }\ P=\pmat{1&i\cr-i&1}\Bigr\}, $$
for the choices of $A$ and $B$ in (\ref{BinDihedralgens}), then
$$ PAP^{-1}=AB, \qquad PBP^{-1}=-B. $$
Since $-1$ in the centre of both $H_{r}^{D_2}$ and $K_{r}^{D_2}$, we therefore have
	$$ P(K_{r}^{D_2})P^{-1}=P\inpro{\zeta_{2^{r+2}}A,-\zeta_{2^{r+2}}B}P^{-1}
	= \inpro{\zeta_{2^{r+2}}AB,\zeta_{2^{r+2}}B}
= \inpro{A,\zeta_{2^{r+2}}B}=H_{r}^{D_2}.$$
The above two cases for $n$ give (iii).

We now consider which of these 
groups are reflection groups.
For an imprimitive reflection group $G(m,p,2)$ to have collineation group $D_n$ (of
order $2n$), we must have
$$ {|G(m,p,2)|\over |Z(G(m,p,2))|}= {2m\over\gcd(p,2)}=2n, $$
which gives two cases: $p$ is odd and $m=n$, or $p$ is even and $m=2n$.
Given the collineation preserving inclusion
$$ G(m,p,2) \lhd G(m,\gcd(p,2),2), $$
we have, for the cases $p$ odd and $p$ even, the possible minimal
reflection groups (based on their orders being the minimal possible)
\begin{align*}
        & G(n,n,2), && \hbox{$n$ odd}, \cr
        & G(n,b,2), && \hbox{$n$ even}, \ n=2^{r+1}b, \ \gcd(2,b)=1, \cr
        & G(2n,2n,2), && \hbox{$n\ge 2$}.
\end{align*}
These are easily seen to give the minimal groups in (\ref{MinimalDn-refgroups}),
and we will see later 
that there are no more cases
(Corollary \ref{G(m,p,2)minimalgroups}), i.e.,
no other $G(m,p,2)$ are minimal groups.
\end{proof}

For $n\ne2$ even, $\mD_n:=\inpro{A,B}$ 
has two normal subgroups $\inpro{A^2,B},\inpro{A^2,AB}\cong\mD_n$,
and we observe from the proof that it is the second of these 
that appears as $\spec(K_r^{D_n})$. 
We also observe that for $n=2$, we have 
$$ H_0^{D_2}=G(4,4,2) \cong G(2,1,2)=K_0^{D_2}. $$

\begin{remark}
For $n$ odd, \cite{NVT08} define a single family of groups
$$ H_k:=\inpro{\pmat{\zeta_n&0\cr0&\zeta_n^{-1}},\zeta_{2^k}\pmat{0&1\cr1&0}},
\qquad k\ge 1. $$
These relate to our presentation as follows. 
Since $(\zeta_{2^{r+2}}B)^2=-\zeta_{2^{r+1}}$, for $r\ge1$, we have
$$
H_{r}^{D_n} 
=\inpro{\pmat{\zeta_{2n}&0\cr0&\zeta_{2n}^{-1}},\zeta_{2^{r+2}}\pmat{0&i\cr i&0}}
=\inpro{\pmat{\zeta_{n}&0\cr0&\zeta_{n}^{-1}},\zeta_{2^{r+2}}\pmat{0&1\cr 1&0}}
=H_{r+2}, \quad r\ge1 \quad\hbox{($n$ odd)},
$$
and
\begin{align*}
& H_2=\inpro{\pmat{\zeta_n&0\cr0&\zeta_n^{-1}},i\pmat{0&1\cr1&0}}
=\inpro{\pmat{\zeta_{2n}&0\cr0&\zeta_{2n}^{-1}},\pmat{0&i\cr i&0}}
=\inpro{A,B}=\mD_n, \cr
& H_1=\inpro{\pmat{\zeta_n&0\cr0&\zeta_n^{-1}},-\pmat{0&1\cr1&0}}
=\inpro{A^2,iB}\cong D_n. 
\end{align*}
\end{remark}

\begin{table}[!h]
\caption{The minimal finite subgroups $H$ of $GL(\CC^2)$ corresponding
to the collineation groups $A_4,S_4,A_5$ (primitive) and
$D_n$, $n\ge2$ (imprimitive).
All appear just once, except for $H_r^{D_2}\cong K_r^{D_2}$
	(here $\mC_4\cong\mD_1$).
	}
        \vskip0.3truecm
\label{GL(2,C)-minimalgroups-table}
\begin{center}
\begin{tabular}{ |  >{$}l<{$} | >{$}l<{$} | >{$}l<{$} | >{$}l<{$} | >{$}l<{$} | >{$}l<{$} | >{$}l<{$} | }
\hline 
&&&&&&\\[-0.4cm]
	H/Z(H) & H & |H| & |Z(H)| & {\mult(H)\over Z(H)} & \mG &
	\spec(H) \\[0.1cm]
\hline
&&&&&&\\[-0.3cm]
A_4 & \mT & 24 & 2 & 1 & \mT & \mT \\
A_4 & H_r^{A_4},\, r\ge0 & 3^r\cdot 24 & 3^r\cdot 2 & C_3 & \mT & \mD_2 \\
&&&&&&\\[-0.3cm]
S_4 & \mO & 48 & 2 & 1 & \mO & \mO \\
S_4 & H_r^{S_4},\, r\ge0 & 2^r\cdot 48 & 2^r\cdot 2 & C_2 & \mO & \mT \\
&&&&&&\\[-0.3cm]
A_5 & \mI & 120 & 2 & 1 & \mI & \mI \\
&&&&&&\\[-0.3cm]
D_n,\, \hbox{$n$ odd}  & \mD_n & 4n & 2 & 1 & \mD_n & \mD_n \\
& D_n & 2n & 1 & C_4 & \mD_n & \mC_{n} \\
& H_r^{D_n}, \, r\ge 1 & 2^r\cdot 4n & 2^{r}\cdot 2 & C_2 & \mD_n & \mC_{2n} \\
&&&&&&\\[-0.3cm]
D_n, \, \hbox{$n$ even} & \mD_n & 4n & 2 & 1 & \mD_n & \mD_n \\
& H_r^{D_n}, \, r\ge 0 & 2^r\cdot 4n & 2^{r}\cdot 2 & C_2 & \mD_n & \mC_{2n} \\
	& K_r^{D_n}, \, r\ge 0 & 2^r\cdot 4n & 2^{r}\cdot 2 & C_2 & \mD_n & \mD_{n/2} \\ [0.1cm]
\hline
\end{tabular}
\end{center}
\end{table}

We now seek to describe all the finite (collineation preserving) subgroups
of $GL(\CC^2)$ uniquely in terms of the minimal groups.
For a given minimal group $H$, all the subgroups obtained from it by 
adding additional scalars, i.e.,
$$ M_a(H)=\inpro{H,\zeta_{|Z(H)|a},}, \qquad a\ge1, $$
have different orders, and hence are different groups.
For $A_5$, there is just one minimal group $\mI$, and this is the end of the
story. In the other cases, we need to determine when $M_a(H)=M_b(K)$
for different minimal groups $H$ and $K$, and then restrict the parameters
so that a given group is counted exactly once.
We begin with an example.

\begin{example}
The reflection group $G_{13}$ can be obtained from minimal groups in two ways, i.e.,
$$ G_{13}=\inpro{\zeta_4}G_{12}=\inpro{\zeta_4}\mO \Iff
M_2(H_0^{S_4})=M_2(\mO), \quad G_{12}=H_0^{S_4}. $$
Here, 
we observe that ${\mult(G)/Z(G)}$ can reduce as we add scalars
to a minimal group
$${\mult(M_a(\mO))\over Z(M_a(\mO))}={\inpro{\zeta_{2a}}\over\inpro{\zeta_{2a}}}=1,
\qquad
	{\mult(G_{12})\over Z(G_{12})} ={\inpro{i}\over\inpro{-1}}=C_2,\quad
	{\mult(M_2(G_{12}))\over Z(M_2(G_{12}))} ={\inpro{i}\over\inpro{i}}=1,
$$
and, equivalently, by (\ref{multHquotbyZH}), $\spec(G)$ can increase
$$ \spec(\mO)=\spec(M_a(\mO))=\mO, \qquad
	\spec(G_{12})=\mT, \quad \spec(M_2(G_{12}))=\mO. $$
\end{example}

In the our classification of the finite irreducible subgroups of
$GL(\CC^2)$ below, we will choose (when a choice must be made) the
minimal group $H$ giving a subgroup $G$, so that
$$ {\mult(G)\over Z(G)} = {\mult(H)\over Z(H)} \cong {\mG\over\spec(H)}, $$
i.e., $G_{13}$ would appear as $M_2(\mO)$, rather than as $M_2(H_0^{S_4})$.

It is convenient to formalise the observation of (\ref{trivialspecH}), 
and give a direct proof which does not utilise the isomorphism  (\ref{specHeqkergThH}).


\begin{lemma}
\label{mult(G)/Z(G)-lemma}
Let $G$ be a finite irreducible subgroup of $GL(\CC^2)$. 
Then
$$ {\mult(G)\over Z(G)}= 1 \Iff G=M_a(\mG),\quad a\ge1,\ \mG=\mD_n,\mT,\mO,\mI.
$$
\end{lemma}

\begin{proof} Suppose that ${\mult(G)\over Z(G)}= 1$, i.e., $\mult(G)=Z(G)$,
and that $G=M_a(H)$ for some minimal group $H$. If $H\ne D_n$, for $n$ odd,
then $-1\in Z(H)\subset Z(G)=\mult(G)$, so that $\inpro{-1}\mG=\mG\subset G$,
and $G=M_a(\mG)$, for $a={1\over2}|Z(G)|={1\over2}|\mult(G)|$. 
For the only other possible case, since $\mult(D_n)=\inpro{i}$, for $n$ odd, 
we have
\begin{equation}
\label{Ma(Dn)cdn}
{\mult(M_a(D_n))\over Z(M_a(D_n))}
={\inpro{\zeta_a}\inpro{i}\over\inpro{\zeta_a}} 
=1 \Iff a=4m \Iff M_a(D_n)=M_{2m}(\mD_n).
\end{equation}
The converse follows 
	from  
	$\mult(M_a(\mG))=\inpro{\zeta_{2a}}\inpro{-1}
	=\inpro{\zeta_{2a}}=Z(M_a(\mG))$.
\end{proof}

\begin{theorem}
\label{GL2Csubgroups-all} (Classification)
The finite irreducible subgroups of $GL(\CC^2)$ are precisely
	the primitive groups (for $A_4$, $S_4$, $A_5$)
\begin{enumerate}[\rm (i)]
\item $M_a(\mT)$, \quad $a\ge1$, \qquad $M_c(H_r^{A_4})$, \quad 
		$c\ge1$, \ $\gcd(3,c)=1$, \ $r\ge0$.
\item $M_a(\mO)$, \quad $a\ge1$, \qquad $M_b(H_r^{S_4})$, \quad 
		$b\ge1$, \ $\gcd(2,b)=1$, \ $r\ge0$.
\item $M_a(\mI)$, \quad $a\ge1$.
\end{enumerate}
	and the imprimitive groups (for $D_n$, $n\ge2$)
\begin{enumerate}[\rm (i)]
\item $M_a(\mD_n)$, \quad $a\ge1$\quad (binary dihedral group).
\item For $n$ odd: \ $M_b(D_n)$ and $M_b(H_r^{D_n})$,
        \quad $b\ge1$, \ $\gcd(2,b)=1$, \ $r\ge2$.
\item For $n$ even: \ $M_b(H_r^{D_n})$ and $M_b(K_r^{D_n})$, \quad $b\ge1$, \ $\gcd(2,b)=1$, \  $r\ge1$.
\end{enumerate}
These groups are distinct, except for when $n=2$,
$M_b(H_r^{D_2})$ and $M_b(K_r^{D_2})$ are conjugate.
Basic properties of these groups are given in Table \ref{Ma(H)groups-table}.
\end{theorem}

\begin{proof}
We consider the cases in order, observing first that 
$G=M_a(\mG)$, $a\ge1$, 
are groups of different orders, with $\mult(G)/Z(G)=\mult(\mG)/Z(\mG)=1$.
We start by taking each of these groups.

We also note that if $G$ contains (generators for) a collineation 
preserving subgroup $H$ which is a minimal group, then $G=M_a(H)$, 
for $a=|Z(G)|/|Z(H)|$.

We recall that $\zeta_{3^r}=(\zeta_{3^{r+1}}(\zeta_3 Z))^3\in H_r^{A_4}$
	and $\zeta_{3^{r+1}}\not\in H_r^{A_4}$, 
so that
\begin{align*}
3\divides c & \Iff \zeta_{3^{r+1}}\in M_c(H_r^{A_4}) \cr
& \Iff \zeta_3 Z=\zeta_{3^{r+1}}^{-1}(\zeta_{3^{r+1}}(\zeta_3 Z))\in M_c(H_r^{A_4}) \cr
	& \Iff \mT\subset M_c(H_r^{A_4}) \cr
	& \Iff  M_c(H_r^{A_4})=M_{3^r c}(\mT).
\end{align*}
Thus, we require the condition $\gcd(c,3)=1$ to obtain additional groups
$G=M_c(H_r^{A_4})$, which satisfy
$$ {\mult(G)\over Z(G)} 
={Z(M_c(H_r^{A_4}))\mult(H_r^{A_4})\over Z(M_c(H_r^{A_4}))}
= {\inpro{\zeta_{3^r\cdot 2c}}\inpro{\zeta_{3^{r+1}\cdot 2}}\over\inpro{\zeta_{3^r\cdot 2c}}}
= {\inpro{\zeta_{3^{r+1}\cdot 2c}}\over\zeta_{3^r\cdot 2c}}
=C_3. $$
This condition also ensures that the groups obtained have different orders, and hence
	are not conjugate to each other, i.e., if $\gcd(c,3)=\gcd(c',3)=1$, then 
	$$ |M_c(H_r^{A_4})|=|M_{c'}(H_s^{A_4})| 
\Iff 3^r \cdot 24c = 3^s \cdot 24c' 
\Iff r=s, \quad c=c'. $$
Thus, we obtain (i). Similarly,
$\zeta_{2^{r+1}}\in H_r^{S_4}$, $\zeta_{2^{r+2}}\not\in H_r^{S_4}$,
which gives
\begin{align*}
2\divides b 
&\Iff \zeta_{2^{r+1}}\in M_b(H_r^{S_4}) 
\Iff \sqrt{i}R\in M_b(H_r^{S_4}) 
\Iff \mO\in M_b(H_r^{S_4}) \cr
&\Iff M_b(H_r^{S_4}) = M_{2^r b}(\mO), 
\end{align*}
giving additional groups $G=M_b(H_r^{S_4})$ for $\gcd(b,2)=1$,
	with $\mult(G)/Z(G)=C_2$, and
different orders, since
$$ |M_b(H_r^{S_4})|=|M_{b'}(H_s^{S_4})| 
\Iff 2^r \cdot 48b = 2^s \cdot 48b'
\Iff r=s, \quad b=b'. $$

For the groups $M_a(H)$ for $D_n$, we will use Lemma \ref{mult(G)/Z(G)-lemma}
to obtain the restriction on $a$ giving additional groups, and then
show they are all different. For $H=H_r^{D_n},K_r^{D_n}$, we have
$$ {\mult(M_a(H))\over Z(M_a(H))}
= {Z(M_a(H))\mult(H)\over Z(M_a(H))}
= {\inpro{\zeta_{2^{r+1}a}}\inpro{\zeta_{2^{r+2}}}\over\inpro{\zeta_{2^{r+1}a}}}
= \begin{cases}
C_2,  &\gcd(2,a)=1;\cr
	1, & 2\divides a,
\end{cases} $$
so that we require that $a$ be odd. The corresponding additional groups
$M_b(H_r^{D_n})$, $\gcd(2,b)=1$, have different orders, and so are different groups.

For $n$ even, there are also groups $M_b(K_r^{D_n})$, $\gcd(2,b)=1$, each
different from each other, which we claim are different from 
the corresponding group $M_b(H_r^{D_n})$ of the same order. To distinguish 
these groups, we observe that
$$ A\in \mC_{2n} = \spec(H_r^{D_n})\subset\spec(M_b(H_r^{D_n})), \quad
A\not\in \mD_{n/2}=\spec(K_r^{D_n}), $$
and, since $K_r^{D_n}=\inpro{\zeta_{2^{r+2}}A,\zeta_{2^{r+2}}B}$, we have
$$ A\in M_b(K_r^{D_n})
\Iff \zeta_{2^{r+2}} \in Z(M_b(K_r^{D_n}))
\Iff 2^{r+2} \divides 2^r\cdot 2b
\Iff \gcd(b,2)=2, $$
i.e., the groups $M_b(H_r^{D_n})$ and $M_b(K_r^{D_n})$ are not conjugate
when $\gcd(b,2)=1$, with their $\spec(G)$ subgroups given by $\mC_{2n}$ and 
$\mD_{n/2}$ (which differ for $n\ne2$).

Finally, for $n$ odd, we have the groups $M_a(D_n)$ to consider.
By (\ref{Ma(Dn)cdn}), the condition $\gcd(4,a)=4$ gives $M_a(D_n)=M_{a/2}(\mD_n)$,
so we may suppose that $\gcd(4,a)=1,2$. In these cases, we have
$$ {\mult(M_a(D_n))\over Z(M_a(D_n))}
={\inpro{\zeta_a}\inpro{i}\over\inpro{\zeta_a}} 
=\begin{cases}
C_4, & \gcd(4,a)=1; \cr
C_2, & \gcd(4,a)=2.
\end{cases}
$$
For $\gcd(4,a)=1$, i.e., $a$ odd, we get additional groups $G$ with $\mult(G)/Z(G)=C_4$.
For $\gcd(4,a)=2$, i.e., $a$ even with $a/2$ odd, using
$A^n=-1$ and $Z(H_0^{D_n})=\inpro{-1}$, we have
$$ M_a(D_n)=\inpro{-A,iB,\zeta_a}=\inpro{A,iB,\zeta_{a/2}}=M_{a/2}(H_0^{D_n}), $$
so that there are no additional groups in this case.
\end{proof}

As we have shown,
the classification of the finite subgroups of 
$GL(\CC^2)$, equivalently of $U(\CC^2)$, naturally follows 
from that of $PGL(\CC^2)$, and it is described in \cite{C76} 
(Theorem on page 393) as being ``well-known''.
In addition to the presentation of \cite{NVT08}, we have found
Theorem 2.2 of \cite{FP04} (which considers minimal numbers of
generators -- these must contain generators for the collineation group) and  Theorem 11.4 of \cite{BFO21} (which considers
the imprimitive/monomial subgroups of $GL(\CC^p)$ for $p$ prime).
Other more general references include \cite{D64} (see the final section) 
and \cite{C91}.

The following table summarises the 
properties of the irreducible groups $G=M_a(H)$, with 
$H$ minimal and $\mult(G)/Z(G)=\mult(H)/H$.

\begin{table}[H]   
\caption{The finite irreducible subgroups $G$ 
of $GL(\CC^2)$, 
in terms of the minimal groups.
Here $a,b,c\ge1$, with $\gcd(2,b)=1$, $\gcd(3,c)=1$, 
($a,b,c=1$ give the minimal groups),
	and $\hat R=\sqrt{i}R$, $\hat Z=\zeta_3 Z$, $\hat M=\zeta_5 M$.
The only repeated groups are $M_b(K_r^{D_2})\cong M_b(K_r^{D_2})$. }
\vskip-0.3truecm
\label{Ma(H)groups-table}
\begin{center}
\begin{tabular}{ | >{$}l<{$} | >{$}l<{$} | >{$}l<{$} | >{$}l<{$} | >{$}l<{$} | >{$}l<{$} |>{$}l<{$} |}
\hline
&&&&&&\\[-0.3cm]
G & |G| & G/Z(G) & \spec(G) & |Z(G)| & {\mult(G)\over Z(G)} & \hbox{generators of $G$} \\[0.1cm]
\hline
&&&&&&\\[-0.3cm]
M_a(\mT) & 24a & A_4 & \mT & 2a & 1 & \hat R^2,\hat Z,\zeta_{2a} \\
M_c(H_r^{A_4}),\, r\ge0 & 3^r\cdot 24 c & A_4 & \mD_2 & 3^r\cdot 2c & C_3 
	& \hat R^2,\zeta_{3^{r+1}}\hat Z, \zeta_{3^r\cdot 2c} \\
&&&&&&\\[-0.3cm]
M_a(\mO) & 48a & S_4 & \mO & 2a & 1 & \hat R,\hat Z,\zeta_{2a} \\
M_b(H_r^{S_4}),\, r\ge0 & 2^r\cdot 48b & S_4 & \mT & 2^r\cdot 2b & C_2 & \zeta_{2^{r+2}}\hat R,\hat Z,\zeta_{2^r\cdot 2b} \\
&&&&&&\\[-0.3cm]
M_a(\mI) & 120a & A_5 & \mI & 2a & 1 & \hat R^2,\hat M,\zeta_{2a} \\
&&&&&&\\[-0.3cm]
M_a(\mD_n) & 4n  a & D_n & \mD_n & 2a & 1 & A,B,\zeta_{2a} \\
&&&&&&\\[-0.4cm]
M_b(H_r^{D_n}),\, r\ge 1 & 2^r\cdot 4nb & D_n,\, \hbox{$n$ odd}	& \mC_{2n} & 2^{r} \cdot 2b & C_2 & A,\zeta_{2^{r+2}}B, \zeta_{2^{r} \cdot 2b} \\
	M_b(D_n) & 2nb & D_n,\, \hbox{$n$ odd}	& \mC_{n} & b & C_4 & -A,iB,\zeta_{b}  \\
&&&&&&\\[-0.4cm]
	M_b(H_r^{D_n}), \, r\ge 0 & 2^r\cdot 4nb & D_n, \, \hbox{$n$ even} & \mC_{2n} & 2^{r} \cdot 2b & C_2 & A,\zeta_{2^{r+2}}B,\zeta_{2^{r}\cdot 2b} \\
	M_b(K_r^{D_n}), \, r\ge 0 & 2^r\cdot 4nb & D_n, \, \hbox{$n$ even} & \mD_{n/2} & 2^{r} \cdot 2b & C_2 & \zeta_{2^{r+2}}A,\zeta_{2^{r+2}}B,\zeta_{2^r\cdot 2b}  \\ [0.1cm]
\hline
\end{tabular}
\end{center}
\end{table}


The following simple example illustrates 
that for each collineation group ``there are not many 
irreducible groups of a given order'' -- indeed there are either $1$, $2$, or $3$ groups.

\begin{example}
If $G$ is subgroup of $GL(\CC^2)$ with collineation group $A_4$, then 
$|G|=3^r\cdot 24b$, $\gcd(3,b)=1$, and there are two groups with this order, i.e.,
$$ M_{3^rb}(\mT),\qquad 
M_{b}(H_r^{A_4}). $$
Since all reflection groups with collineation group $A_4$ (and also $A_5$) 
are obtained from just 	one minimal group, it follows that there are 
no reflection groups of the same order in these cases. Similarly, it is possible 
that there are two reflection groups with collineation group $S_4$ of the same order,
and this does occur in two cases: 
\begin{align*}
G_8&=M_1(H_1^{S_4}), \quad G_{13}=M_2(\mO), \quad\hbox{(order $96$)}, \cr
 G_{10}&=M_3(H_1^{S_4}), \quad G_{15}=M_6(\mO), \quad\hbox{(order $288$)}.
\end{align*}
\end{example}

Motivated by this example, we now give an alternative indexing of the
irreducible subgroups of $GL(\CC^2)$ by their order and collineation
group (this is essentially what was done in the theorem of page 393 in \cite{C76}).
Let
$$ G=G(a,\mG,\mH), \qquad a\ge1, \quad \mG\in\{\mT,\mO,\mI,\mD_n\}, $$
be the unique group from Theorem \ref{GL2Csubgroups-all}, with 
\begin{itemize}
\item The collineation group of $G$ is given by $\mG$
\item $\mH=\spec(G) \ne \mC_n$.
\item $|G|=|\mG|a$, equivalently $|Z(G)|=2a$.
\end{itemize}
This accounts for all groups from Theorem \ref{GL2Csubgroups-all}, except
for $M_b(D_n)$ for $nb$ odd, which we denote similarly by
$$ G(b,\mD_n,\mC_n), \qquad \hbox{($n$ odd, $b$ odd)}. $$

\begin{corollary} 
\label{ClassificationbyOrder}
(Classification by order)
The finite irreducible subgroups $G$ of $GL(\CC^2)$ have two types
	\begin{enumerate}[\rm (i)]
\item The centre of $G$ has odd order, i.e., $-1\not\in G$, equivalently $\gcd(|G|,4)=2$, and
$$ G=G(b,\mD_n,\mC_n)=\inpro{\zeta_b}\times D_n, 
		\quad |G|=2nb, \quad |Z(G)|=b,  \qquad \hbox{($n$ odd, $b$ odd)}, $$
where $\mult(G)/Z(G)={\mD_n/\mC_n}=C_4$.
\item The centre of $G$ has even order, i.e., $-1\in G$, equivalently $\gcd(|G|,4)=4$, and
	$$ G=G(a,\mG,\mH), 
		\quad |G|=|\mG|a,
		\quad |Z(G)|=2a,
		\qquad \mG=\mD_n,\mT,\mO,\mI, \quad
		\mG/\mH=1,C_2,C_3. $$
\end{enumerate}
The above groups are distinct, and their properties and presentation in terms
of the minimal groups is given in Table \ref{G(a,G,H)-table}. 
We note, in particular, that
the minimal groups are
\begin{equation}
\label{G(a,G,H)-minimal}
G(1,\mD_n,\mC_n)=D_n,\quad \hbox{for $n$ odd}, \qquad
G(a,\mG,\mH), \quad \hbox{$a=|\mG/\mH|^r$, $r\ge0$ }.
\end{equation}
\end{corollary}

\begin{proof}
By examination of Table \ref{G(a,G,H)-table}, 
every group appears exactly once, either as type (i) or (ii).
Here the conjugate groups $M_b(H_r^{D_2})$ and $M_b(H_r^{D_2})$
appear as $G(a,\mD_2,\mC_4)$ and $G(a,\mD_2,\mD_1)$, which are same 
	(since $\mC_4\cong \mD_1$).

The minimal groups are
$$ G(1,\mT,\mT),\  G(3^r,\mT,\mD_2), \qquad 
 G(1,\mO,\mO),\ G(2^r,\mO,\mT), \qquad G(1,\mI,\mI) \qquad r\ge 0, $$
\begin{align*}
\hbox{$n$ odd:} & \qquad G(1,\mD_n,\mD_n),\quad G(2^r,\mD_n,\mC_{2n}),\ r\ge 1, 
	\quad G(1,\mD_n,\mC_n), \cr
\hbox{$n$ even:} & \qquad G(1,\mD_n,\mD_n),\quad 
	G(2^r,\mD_n,\mC_{2n}),\ r\ge 0, \quad
	G(2^r,\mD_n,\mD_{n/2}),\ r\ge 0.
\end{align*}
and these can be summarised by (\ref{G(a,G,H)-minimal}).
\end{proof}

From Corollary \ref{ClassificationbyOrder},
we observe that there are no irreducible subgroups of $GL(\CC^2)$ of odd order,
an instance of Frobenius's dimension divisibility theorem for representations.

\begin{table}[!h]
	\caption{The finite irreducible subgroups $G=G(a,\mG,\mH)$ of $GL(\CC^2)$, 
as they evolve from the minimal groups.
Here $a,b,c\ge1$, with $\gcd(2,b)=1$, $\gcd(3,c)=1$. The 
	notation of Cohen \cite{C76} is included in the last column, for 
	comparison.
We observe that every family with $\mH\ne\mG$ contains infinitely many 
minimal subgroups.
	}
     \vskip0.3truecm
	\label{G(a,G,H)-table}
\resizebox{\textwidth}{!}{
\begin{tabular}{ |  >{$}l<{$} | >{$}l<{$} | >{$}l<{$} | >{$}l<{$} | >{$}l<{$} | >{$}l<{$} | }
\hline
&&&&&\\[-0.3cm]
	G=G(a,\mG,\mH) & G & |G| & |Z(G)| & {\mG/\mH} & \hbox{Cohen's notation} \\[0.1cm]
\hline
&&&&&\\[-0.3cm]
G(a,\mT,\mT) & M_a(\mT) & 24a & 2a & 1 & \mu_{2a}\mT \\
G(a,\mT,\mD_2) & M_c(H_r^{A_4}),\ a=3^rc 
	& 24a & 2a & C_3 & (\mu_{6a}|\mu_{2a};\mT|\mD_2) \\
&&&&&\\[-0.3cm]
G(a,\mO,\mO) & M_a(\mO) & 48a & 2a & 1 & \mu_{2a}\mO  \\
	G(a,\mO,\mT) & M_b(H_r^{S_4}),\ a=2^rb 
	& 48a & 2a & C_2 & (\mu_{4a}|\mu_{2a};\mO|\mT) \\
&&&&&\\[-0.3cm]
G(a,\mI,\mI) & M_a(\mI) & 120a & 2a & 1 & \mu_{2a}\mI  \\
&&&&&\\[-0.3cm]
G(a,\mD_n,\mD_n) & M_a(\mD_n) & 4na & 2a & 1 & \mu_{2a}\mD_n \\
&&&&&\\[-0.4cm]
	G(a,\mD_n,\mC_{2n}) & M_b(H_r^{D_n}),\, a=2^rb,\, 
	\hbox{$a$ even, $n$ odd} & 
	4na & 2a & C_2 & (\mu_{4a}|\mu_{2a};\mD_n|\mC_{2n}) \\
	& M_{2a}(D_n),\, \hbox{$a$ odd, $n$ odd} & & & &  \\
	G(b,\mD_n,\mC_n) & M_b(D_n),\, \hbox{$b$ odd, $n$ odd}  & 2nb & b & C_4 & (\mu_{4b}|\mu_{b};\mD_n|\mC_{n}) \\
&&&&&\\[-0.4cm]
G(a,\mD_n,\mC_{2n}) & M_b(H_r^{D_n}),\, a=2^rb,\, 
\hbox{$n$ even} & 4na & 2a & C_2 & (\mu_{4a}|\mu_{2a};\mD_n|\mC_{2n}) \\
G(a,\mD_n,\mD_{n/2}) & M_b(K_r^{D_n}),\, a=2^rb,\,  
	\hbox{$n$ even} & 4na & 2a & C_2 & (\mu_{4a}|\mu_{2a};\mD_n|\mD_{n/2}) \\[0.1cm]
\hline
\end{tabular}
   }  
\end{table}

For the $G(a,\mD_n,\mC_{2n})$ entry of Table \ref{G(a,G,H)-table},
we note that $M_{2a}(D_n)=M_a(H_0^{D_n})$ for $a$ odd and $n$ odd,
and so $G(a,\mD_n,\mC_{2n})=M_b(H_r^{D_n})$, where $a=2^rb$ and $b$ odd,
in both cases.

\begin{remark}
This presentation is essentially that of the (unproved) theorem stated on page 393 
of Cohen \cite{C76}, where the notation there can be seen to correspond to ours 
as follows
$$ G=G(a,\mG,\mH) \Iff (\mult(G)|Z(G);\mG|\spec(G)). $$
There, $\mult(G)$, $Z(G)$, and $\spec(G)$ are first defined abstractly as
$$ \mu_{wd}=\inpro{\zeta_{wd}}=\mult(G), \qquad
\mu_{d}=\inpro{\zeta_{d}}=Z(G), \qquad
\mH=\spec(G)\lhd\mG, $$
an isomorphism $\varphi:{\mult(G)/Z(G)}\to{\mG/\mH}$ is supposed to exist,
and the groups are defined as
$$ \{\gl g:\gl\in\mult(G),g\in\mG,\varphi(\gl\inpro{\zeta_{|Z(G)|}})= g\mH\}. $$
The conditions on $wd$, $d$ and $\mH$ for such a group to exist,
e.g., $\mG/\mH$ is cyclic, are given, and it is concluded the 
particular choice of the isomorphism $\varphi$ does not affect the group,
and so it is omitted in the notation. It is observed that every finite
	irreducible subgroup $G$ of $GL(\CC^2)$ has this form, where, for example,
$$ \inpro{\zeta_{wd}}=G\cdot SL(\CC^2) \cap Z(GL(\CC^2))
       =\{\gl\in\CC^*:\gl\in G\cdot SL(\CC^2)\}.  $$
To illustrate the asserted correspondence, we show that
$\mult(G)=\inpro{\zeta_{wd}}$. We have
$$ \gl \in \mult(G)
\Iff \gl g\in G, \ g\in\mG
\Iff \gl \in G g^{-1}, \ g\in\mG
\Implies \gl \in G\cdot SL(\CC^2). $$
For the converse of the last implication
\begin{align*}
\gl \in G\cdot SL(\CC^2) & \Implies \gl = hU, \quad h\in G,\ U\in SL(\CC^2)\cr
	&\Implies \gl= G g^{-1}, \quad g=U^{-1}=\gl^{-1}h\in\mG.
\end{align*}
There are two small oversights in the list of groups given:
\begin{itemize}
\item Groups $G=G(a,\mD_n,\mC_n)$ for $n$ odd and $a$ even, with 
$\mult(G)/Z(G)=C_4$, are listed.
However, no such groups exist, since as soon as $a$ becomes even (for $n$ odd),
the scalar $-1$ is in the group, which effectively then becomes
		$G(a,\mD_n,\mC_{2n})$, which has $\mult(G)/Z(G)=C_2$.
\item It is not observed that the groups $G(a,\mD_2,\mC_4)$ and $G(a,\mD_2,\mD_1)$ 
are conjugate.
\end{itemize}
\end{remark}

It is relatively easy to identify ``groups in the wild'' using the
$G(a,\mG,\mH)$ indexing:
\begin{itemize}
\item For $G$ a finite subgroup of $GL(\CC^2)$ it is easily 
checked whether it is irreducible, e.g., the condition
		$$ \sum_{g\in G}\trace(g)\trace(g^{-1})=|G|, $$
which can be calculated at the level of conjugacy classes if need be.
\item The group $\mG$ is determined by the collineation
	group $G/Z(G)$, and, in all by a few cases, just by the order of the collineation group.
\item The group $\mH=\spec(G)$, in the cases where there could be more than one 
	choice, is determined by $\mH=G\cap \mG$ (with $\mG$ given by
		scalar multiples of elements of $G$),
		or as $\mH=\{g\in G:\det(g)=1\}$.
\item The scalar $a$ is obtained from $|Z(G)|$ or $|G|$
		using
	the formulas in Table \ref{G(a,G,H)-table}.
\end{itemize}
In other words, each finite irreducible subgroup $G$ of $GL(\CC^2)$ is 
uniquely determined by 
$$ |G|, \quad {G\over Z(G)}, \quad \spec(G). $$

\begin{example}
The (primitive) reflection groups with collineation groups $A_4$ and $S_4$ are
\begin{align*}
G_5&=G(3,\mT,\mT)=M_3(\mT),\quad && G_7=G(6,\mT,\mT)=M_6(\mT), \\[0.1cm]
\greyeq{G_4}&=G(1,\mT,\mD_2)=M_1(H_0^{A_4}),\quad && G_6=G(2,\mT,\mD_2)=M_2(H_0^{A_4}), \\[0.1cm]
G_9&=G(4,\mO,\mO)=M_4(\mO),\quad && G_{11}=G(12,\mO,\mO)=M_{12}(\mO),\quad \\
G_{13}&=G(2,\mO,\mO)=M_2(\mO),\quad &&	G_{15}=G(6,\mO,\mO)=M_6(\mO),\quad \\[0.1cm]
\greyeq{G_{12}}&=G(1,\mO,\mT)=M_1(H_0^{S_4}),\quad 
&& G_{14}=G(3,\mO,\mT)=M_3(H_0^{S_4}), \\
\greyeq{G_8}&=G(2,\mO,\mT)=M_1(H_1^{S_4}),\quad 
&& G_{10}=G(6,\mO,\mT)=M_3(H_1^{S_4}),
\end{align*}
where the reflection groups in grey are minimal subgroups of $GL(\CC^2)$.
\end{example}

We now seek to identify the imprimitive reflection groups $G(m,p,2)$.
We observed in (\ref{Gmp2collgrp}) that these have order $2m^2/p$ and
collineation group $D_n$, where
\begin{equation}
\label{G(m,p,2)-Dn}
n={m\over\gcd(p,2)} \Implies 
\mG=
\begin{cases}
\mD_m, & \hbox{$p$ odd}; \cr
\mD_{m/2}, & \hbox{$p$ even}.
\end{cases}
\end{equation}
If $m$ (and hence $p$) is odd, then $G(m,p,2)$ has order $2nb=2mb$, with
$b=m/p$ odd, so (see Table \ref{G(a,G,H)-table})
$$ G(m,p,2)=G({m\over p},\mD_m,\mC_m), \qquad 
\hbox{$p$ odd,\quad ${m\over p}$ odd}, $$
otherwise, $G(m,p,2)$ has order
$$ |G(m,p,2)|= 4na= 4{m\over\gcd(p,2)}a, \quad
a={m\over p}{\gcd(p,2)\over2} \Implies
G(m,p,2)= G(a,\mG,\mH), $$
where $\mH$ is to be determined. In this regard, by (\ref{G(m,p,2)form}), we have
$$ G(m,p,2) 
= \inpro{ \pmat{\zeta_m&0\cr0&\overline{\zeta_m}},\pmat{0&1\cr1&0},
\pmat{\zeta_{m\over p}&0\cr0&1}}, \qquad 
\mC_m=\mC_{\gcd(p,2)n}\subset G(m,p,2). 
$$
We recall $D_n=\inpro{A,B}$ as given in (\ref{BinDihedralgens}).
For ${m\over p}$ even, we have
$$ B=\pmat{0&1\cr-1&0}
= \pmat{0&1\cr1&0}\pmat{\zeta_{m\over p}&0\cr0&1}^{m\over 2p} 
\in G(m,p,2), $$
so that $\mH$ is not cyclic in this case, and,
since $\mC_{2m}\not\subset G(m,p,2)$,
we have
$$ G(m,p,2)=G({m\over 2p},\mD_m,\mD_{m/2}), \qquad
\hbox{$p$ odd,\quad ${m\over p}$ even}, $$
$$ G(m,p,2)=G({m\over p},\mD_{m/2},\mD_{m/2}), \qquad
\hbox{$p$ even,\quad ${m\over p}$ even}. $$
Finally, for ${m\over p}$ odd (and $p$ even), we cannot 
have $B\in G(m,p,2)$, 
and so we have
$$ G(m,p,2)=G({m\over p},\mD_{m/2},\mC_{m}), \qquad
\hbox{$p$ even,\quad ${m\over p}$ odd}. $$

\begin{example}
\label{Gmp2-G(a,G,H)}
The (irreducible) imprimitive reflection subgroups of $GL(\CC^2)$ are
$$ G(m,p,2) =
\begin{cases}
G({m\over p},\mD_m,\mC_m), & \hbox{$p$ odd, \ ${m\over p}$ odd}; \cr
G({m\over 2p},\mD_m,\mD_{m/2}), &\hbox{$p$ odd, \ ${m\over p}$ even}; \cr
G({m\over p},\mD_{m/2},\mD_{m/2}), & \hbox{$p$ even, \ ${m\over p}$ even}; \cr
G({m\over p},\mD_{m/2},\mC_{m}), & \hbox{$p$ even, \ ${m\over p}$ odd}. 
\end{cases}
$$
This agrees with the theorem on page 394 of \cite{C76}, except for the
group for $p$ even and ${m\over p}$ even, which is given as
	$G({m\over p},\mD_{m},\mD_{m})$ (which has the wrong order).
	The case $(m,p)=(2,2)$ must be excluded, since $G(2,2,2)$ is
	not irreducible.
\end{example}

We now give a generalisation of Lemma \ref{hatGformulas},
which will allow for the classification of all rank two quaternionic 
reflection groups with imprimitive complexifications.

\begin{lemma} 
\label{G(a,G,H)sclosed}
All the finite irreducible subgroups of $U(\CC^2)$ are
closed under conjugation by the quaternionic matrix $s$ of (\ref{generalform}), 
except for all groups $G$ in the two families
\begin{enumerate}[\rm (i)]
\item $G(a,\mT,\mD_2)=M_c(H_r^{A_4})$, $a\ge1$.
\item $G(b,\mD_n,\mC_n)=M_b(D_n)$, $n$ odd, $b$ odd.
\end{enumerate}
For these groups, $\hat G=\inpro{G,G^s}=\inpro{G,\det(G)}$ 
	(the closure under conjugation by $s$) is
$$ \hat G(a,\mT,\mD_2)=G(3a,\mT,\mT), \qquad 
	\hat G(b,\mD_n,\mC_n)=G(b,\mD_n,\mC_{2n}). $$
\end{lemma}

\begin{proof} As in Lemma \ref{hatGformulas}, we 
verify whether or not $\det(G)\subset G$, to determine whether $G$ is closed under 
conjugation by $s$ or not (respectively).
This can easily be done by checking whether $\det(g_1),\ldots,\det(g_m)\in Z(G)$ 
	for any generators $g_1,\ldots,g_m$ of $G$. 

The condition $\det(g_j)\in Z(G)$ is trivially satisfied for a generator
with $\det(g_j)=1$ or a generator which is a scalar matrix.  
The relevant calculations for the groups listed in Table \ref{Ma(H)groups-table}
(with their centres and generators) are
\begin{align*}
M_c(H_r^{A_4}): & \qquad \det(\zeta_{3^{r+1}}\hat Z)=\zeta_{3^{r+1}}^2 \not\in \inpro{\zeta_{3^r\cdot 2c}}, \cr
M_b(H_r^{S_4}): & \qquad \det(\zeta_{2^{r+2}}\hat R)=\zeta_{2^{r+2}}^2=\zeta_{2^{r+1}} \in \inpro{\zeta_{2^r\cdot 2b}}, \cr
M_b(H_r^{D_n}): & \qquad \det(\zeta_{2^{r+2}} B)=\zeta_{2^{r+2}}^2
	=\zeta_{2^{r+1}} \in \inpro{\zeta_{2^r\cdot 2b}}, \cr
M_b(D_n): & \qquad \det(-A)=1\in \inpro{\zeta_b}, \ 
	 \det(iB)=-1\not\in \inpro{\zeta_b}, \cr
M_b(H_r^{D_n}),M_b(K_r^{D_n}): & \qquad \det(\zeta_{2^{r+2}} A),\det(\zeta_{2^{r+2}} B)
	=\zeta_{2^{r+1}} \in \inpro{\zeta_{2^r\cdot 2b}},
\end{align*}
For the two families not closed under conjugation by $s$, we calculate
	$\hat G=\inpro{G,\det(G)}$
\begin{align*} \hat G(a,\mT,\mD_2) &=\hat M_c(H_r^{A_4})
=\inpro{\hat R^2,\zeta_{3^{r+1}}\hat Z,\zeta_{3^r\cdot 2c},\zeta_{3^{r+1}}}
=\inpro{\hat R^2,\hat Z,\zeta_{3^{r+1}\cdot 2c}} \cr
	&= M_{3^{r+1}\cdot c}(\mT)=G(3a,\mT,\mT), \\[0.2cm]
\hat G(b,\mD_n,\mC_n)
	&=\hat M_b(D_n)
        =\inpro{-A,iB,\zeta_b,-1}
        =\inpro{-A,iB,\zeta_{2b}} \cr
	&= M_{2b}(D_n)
        =G(b,\mD_n,\mC_{2n}).
\end{align*}
\end{proof}

\begin{remark}
During our calculations, we observed that the groups given
in Table \ref{Ma(H)groups-table} are all closed under the transpose,
equivalently (since they are unitary) closed under entry-wise conjugation,
except those in the family $G(a,\mT,\mD_2)$ (which is not 
closed under conjugation by $s$). The mechanics for this appears 
to be related to the fact
	$$ J\det(g)g^{-1} J^{-1}=g^T, \qquad J:=\pmat{0&1\cr-1&0} \quad
\hbox{($s=jJ$)}, $$
so that if $G$ is closed under conjugation by $s$ and $J\in G$ then
$G$ is closed under the transpose, and that our groups are chosen so
that they satisfy natural inclusions. 

As an example, the maximal 
complex reflection group $G_7$ appears $G_7=G(6,\mT,\mT)$, which is invariant
under conjugation by $s$ and under the transpose. 
This contains the normal reflection subgroups
	$G_4=G(1,\mT,\mD_2)$ and $G_6=G(2,\mT,\mD_2)$, which are not
	invariant under $s$ or the transpose. The groups by applying 
	conjugation by $s$ or the transpose to these is the 
	second (normal) copy of $G_4$ and $G_6$ in $G_7$ (as depicted in Figure \ref{G7-Zlat}).
\end{remark}

By our previous discussion (see Lemma \ref{internaldirectsum}), 
all of the rank two finite irreducible subgroups
of $GL(\HH^2)$ with imprimitive complexifications have the 
form $\inpro{G,s}=\cH_1(G)$, with $G$ a finite subgroup of $U(\CC^2)$ which is
closed under conjugation by $s$ (i.e., $\hat G=G$). 
In light of Lemma \ref{G(a,G,H)sclosed},
we define quaternionic groups
\begin{equation}
\label{GH(a,G,H)-defn}
G_\HH(a,\mG,\mH):=\inpro{G(a,\mG,\mH),s}, 
\qquad a\ge1, \quad \mG=\mD_n,\mT,\mO,\mI, \quad {\mG/\mH}=1,C_2,
\end{equation}
where the condition ${\mG/\mH}=1,C_2$ excludes the families (i),(ii), which have
$\mG/\mH=C_3,C_4$.
This has order
$$ |G_\HH(a,\mG,\mH)|=2|G(a,\mG,\mH)|=2a|\mG|. $$
and quotients
$$ {G_\HH(a,\mG,\mH)\over\inpro{\zeta_{2a},s}}=\mG, \qquad
{G_\HH(a,\mG,\mH)\over\inpro{\zeta_{2a}}}=C_2\times {\mG\over\inpro{-1}}, \qquad
{G_\HH(a,\mG,\mH)\over\inpro{s}}=G(a,\mG,\mH).
$$

We can now generalise the classification of Theorems
\ref{TOIfamilyclass} and \ref{G(n,a,b,r)familyclass}.

\begin{lemma}
\label{HaGquatimprimitiveidII}
Let $G\subset U(\CC^2)$ be a primitive group.
Then $\cH_a(G)$ is an imprimitive
quaternionic group 
in precisely the seven cases
\begin{equation}
\label{HaGnotprimitiveII}
	\cH_a(\mG)=G_\HH(a,\mG,\mG), \quad \mG=\mT,\mO,\mI,\ a=1,2, \qquad
	\cH_1(H_0^{S_4})=G_\HH(1,\mO,\mT).
\end{equation}
\end{lemma}

\begin{proof} Since the quaternionic systems of imprimitivity for a group
are shared by its subgroups (which could have more), we first consider
the quaternionic systems of imprimitivity for the minimal 
primitive complex groups. In this regard, since
$$ G_{12}=G(1,\mO,\mT)=M_1(H_0^{S_4}), \quad G_{13}=G(2,\mO,\mO)=M_2(\mO), \quad
	G_{22}=G(2,\mI,\mI)=M_2(\mI), $$
have quaternionic systems of imprimitivity, so do 
the minimal groups $H_0^{S_4}, \mO,\mI$.
Each of the primitive minimal groups $H$ 
(see Table \ref{GL(2,C)-minimalgroups-table}) contains
the monomial subgroup
$$   \inpro{\hat R^2,\hat Z^{-1}\hat R^2\hat Z}
=\inpro{\pmat{i&0\cr0&-i},\pmat{0&i\cr i&0}} \subset \mT. $$
By the proof of Lemma 4.1 of \cite{W26} (a simple calculation), this
implies that if $H$ has a quaternionic system of imprimitivity, 
then it has the form
(\ref{Uimpsysform}), and by the calculation 
in the proof of Lemma \ref{HaGquatimprimitiveid}, it is therefore
given by 
	$$ U={1\over\sqrt{2}} \pmat{1&-j\cr j&-1}. $$
All that remains is to determine when this $U$ conjugates
generators for $G=M_a(H)$ 
	to monomial matrices (it conjugates $s$ to a monomial matrix).
We have already seen that $\hat R,\hat Z,\hat M$ conjugate to
monomial matrices, and that a scalar matrix $\xi$ does if and only 
if $\xi$ is a $4$-th root of unity. Hence, by considering the generators 
given in Table \ref{Ma(H)groups-table}, the conditions for
$G=M_a(H)$ to have a quaternionic system of imprimitivity are
\begin{align*}
M_a(\mT): &\quad a=1,2, \cr
M_c(H_r^{A_4}): &\quad \hbox{not possible}, \cr
M_a(\mO): &\quad a=1,2, \cr
M_b(H_r^{S_4}): &\quad r=0, \ b=1, \cr
M_a(\mI): &\quad a=1,2.
\end{align*}
Thus we obtain the seven imprimitive quaternionic groups of 
	(\ref{HaGnotprimitiveII}).
\end{proof}

\begin{theorem}
\label{GL(H^2)classification}
(Classification) The finite irreducible subgroups 
of $GL(\HH^2)$ which have imprimitive complexifications are precisely the
groups
$$ G_\HH(a,\mG,\mH)=\inpro{G(a,\mG,\mH),s}, \qquad
a\ge1, \quad \mG=\mD_n,\mT,\mO,\mI, \quad {\mG\over\mH}=1,C_2. $$
which are primitive, except for when $\mG=\mD_n$ and the seven 
	exceptional cases of (\ref{HaGnotprimitiveII}).
Of these, the following are primitive (quaternionic) reflection groups
\begin{align}
\label{GH(a,H,G)primrefgroups}
G_\HH(3a,\mT,\mT), & \qquad a\ge1, \cr
G_\HH(2a,\mO,\mO), & \qquad a\ge2, \cr
G_\HH(a,\mO,\mT), & \qquad a\ge2, \quad\gcd(a,4)\ne 4, \cr
G_\HH(a,\mI,\mI), & \qquad a\ge3,  \quad 2\divides a, \ 3\divides a, \ \hbox{or } \ 5\divides a,
\end{align}
	and imprimitive reflection groups
\begin{equation}
\label{GH(a,H,G)imprimrefgroups-exceptions}
G_\HH(1,\mO,\mT), \qquad G_\HH(2,\mO,\mO), \qquad
G_\HH(2,\mI,\mI), 
\end{equation}
\begin{align}
\label{GH(a,H,G)imprimrefgroups}
G_\HH(2a,\mD_n,\mD_n), & \qquad a\ge1, \quad \hbox{$n$ even}, \cr
G_\HH(b,\mD_n,\mC_{2n}), & \qquad \hbox{$b$ odd}, \cr
G_\HH(2^r b,\mD_n,\mD_{n/2}), & \qquad \hbox{$b$ odd}, \quad \hbox{$n$ even}, 
	\quad \hbox{$r:=\nu_2({n\over2})$}.
\end{align}
\end{theorem}

\begin{proof}
It remains to show that the groups of (\ref{GH(a,G,H)-defn}) are different
(not conjugate) and to identify those which are the quaternionic 
reflection groups of Theorem \ref{TOIfamilyclass}.

By Lemma \ref{internaldirectsum}, each $G_\HH(a,\mG,\mH)=\cH_a(G(a,\mG,\mH))$ has a unique complex subgroup
$G(a,\mG,\mH)$ of index $2$ (the type $A$ elements), and so by
Corollary \ref{ClassificationbyOrder} these (normal) subgroups are distinct,
and hence none of groups themselves are conjugate to each other.

We first consider the primitive quaternionic reflection groups
of Theorem \ref{TOIfamilyclass}.
For the $A_4$ groups,
$$ \cT_a(G_5)=\inpro{M_3(\mT),\zeta_{6a},s}
	= \inpro{ M_{3a}(\mT),s} = G_\HH(3a,\mT,\mT). $$
For the $S_4$ groups, it is convenient to split them into those that 
	evolve from $\mO$, and the others (see Figure \ref{G11-Zlat}), i.e..
$\cH_{2a}(\mO)$ and $\cH_b(G_{12}),\cH_b(G_8)$, 
for which 
\begin{align*}
\cH_{2a}(\mO) &= \inpro{\mO,\zeta_{4a},s} = \inpro{M_{2a}(\mO),s} 
	= G_\HH(2a,\mO,\mO),\cr
\cH_b(G_{12}) &=\inpro{H_0^{S_4},\zeta_{2b},s} =\inpro{M_b(H_0^{S_4}),s}
	=G_\HH(b,\mO,\mT), \cr
\cH_b(G_{8}) &=\inpro{H_1^{S_4},\zeta_{4b},s} =\inpro{M_b(H_1^{S_4}),s}
	=G_\HH(2b,\mO,\mT),
\end{align*}
and so we obtain the groups $G_\HH(a,\mO,\mT)$, where $a$ is odd or $a$ is an 
	even multiple of an odd number, i.e., $\gcd(a,4)\ne 4$.
For the next case, we have
$$ \cI_a(\mI)
=\inpro{\mI,\zeta_{2a},s}
	=\inpro{M_a(\mI),s}= G_\HH(a,\mI,\mI), $$
with conditions for $M_a(\mI)$ to be a reflection group.
Lastly, we use Lemma \ref{HaGquatimprimitiveidII} to
restrict the parameters so as to split these groups
into the primitive and imprimitive groups.

We now consider the (imprimitive) $D_n$ groups. It suffices to find those 
$G=G(m,p,2)$, 
with collineation group $D_n$ and $\hat G=G$,
	which are minimal, and to then obtain the 
quaternionic reflection groups from them. 
With $n$ given by (\ref{G(m,p,2)-Dn}),
the relevant cases of Example \ref{Gmp2-G(a,G,H)} are
$$
	G(m,p,2)=\begin{cases}
G({n\over 2p},\mD_n,\mD_{n/2}), &\hbox{$p$ odd, \ ${n\over p}$ even}; \quad
                \hbox{($n$ even)} \cr
G({2n\over p},\mD_{n},\mD_{n}), & \hbox{$p$ even, \ ${2n\over p}$ even}; \quad
		\hbox{($n$ even)}\cr
G({2n\over p},\mD_{n},\mC_{2n}), & \hbox{$p$ even, \ ${2n\over p}$ odd}.
\end{cases}
$$
(Here the case $p$ and ${m\over p}$ both odd is excluded,
	as this gives a group $G=G(b,\mD_n,\mC_n)$ with $\hat G\ne G$), 
For a given $n$, and fixed parity of $p$ and ${n\over p}$, we may minimise
${n\over p}$ (equivalently maximise $p$) to the obtain minimal
reflection groups from which all others are obtained by adding scalars.
This minimisation gives
$$
\begin{cases}
	G(2^r,\mD_n,\mD_{n/2}), &\hbox{$p=b$,  \ ${n\over p}=2^{r+1}$}; \quad
  \hbox{(where $n=2^{r+1}b$, $b$ odd)} \cr
	G({2},\mD_{n},\mD_{n}), & \hbox{$p=n$, \ ${2n\over p}=2$}; \quad
                \hbox{($n$ even)}\cr
G(1,\mD_{n},\mC_{2n}), & \hbox{$p=2n$, \ ${2n\over p}=1$}.
\end{cases}
$$
Above $2^r$ is the $2$-part of ${n\over2}$, i.e., $r=\nu_2({n\over2})$ is the
$2$-adic valuation of ${n\over2}$.
From Table \ref{G(a,G,H)-table}, we can identify these 
minimal complex reflection groups, and then write down the
families of quaternionic reflection groups that they give rise to
\begin{align*}
G(2^r,\mD_n,\mD_{n/2})=M_1(K_r^{D_n}) &\Implies
G_\HH(2^rb,\mD_n,\mD_{n/2})=\inpro{M_b(K_r^{D_n}),s}, \cr
G(2,\mD_n,\mD_n)=M_2(\mD_n) &\Implies
G_\HH(2a,\mD_n,\mD_n)=\inpro{M_{2a}(\mD_n),s}, \cr
G(1,\mD_n,\mC_{2n})=M_1(H_0^{D_n}) &\Implies
	G_\HH(b,\mD_n,\mC_{2n})=\inpro{M_{b}(H_0^{D_n}),s},
\end{align*}
where the restriction on the parameters, i.e., $b$ is odd, 
is so that no group is repeated (see Theorem \ref{GL2Csubgroups-all}).
\end{proof}

The above minimisation of reflection groups $G$ with collineation
group $D_n$ and $\hat G=G$ (in the proof) gives two which 
are minimal groups (see Theorem \ref{GL2Csubgroups-all}), i.e.,
\begin{equation}
\label{twominrefgps}
K_r^{D_n}=G(n,{n\over 2^{r+1}},2) \quad \hbox{($n$ even)}, \qquad
H_0^{D_n}=G(2n,2n,2) \quad \hbox{($n$ even)}.
\end{equation}
We now complete the minimisation (for reflection groups with $\hat G\ne G$)
to obtain the minimal reflection groups for $D_n$.
We recall from \cite{W26c}, that all reflection groups with collineation group
$D_n$ are (collineation preserving) normal reflection subgroups of the 
maximal reflection group $G(2n,2,2)$, which can be 
indexed by the divisors of $(2,2,n)$.

\begin{corollary}
\label{G(m,p,2)minimalgroups}
(Minimal reflection groups)
The minimal imprimitive complex reflection groups with collineation group $D_n$ are
the following 
normal reflection subgroups 
of the maximal reflection group $G(2n,2,2)$
\begin{align*}
\hbox{$n$ odd}: &\qquad G(n,n,2) \quad \hbox{(twice)}, \cr
\hbox{$n$ even}: &\qquad G(2n,2n,2) \quad G(n,{n\over 2^r},2), \quad r=\nu_2(n), 
\quad \hbox{(twice)}.
\end{align*}
	They 
give all reflection groups with collineation 
	group $D_n$ (see Figures \ref{G(2n,2,2)-coll-odd}, \ref{G(2n,2,2)-coll-even}), as follows
\begin{align*}
\hbox{$n$ odd:} & \qquad
\inpro{\zeta_q} G(n,n,2) =
\begin{cases}
G(n,{n\over q},2), & \hbox{$q$ odd}; \cr
G(2n,2{2n\over q},2), & \hbox{$q$ even}, \cr
\end{cases}
	\qquad q\divides 2n \\[0.4cm]
\hbox{$n$ even:} & \qquad
\inpro{\zeta_{2q}}G(2n,2n,2)=G(2n,2{n\over q},2), \qquad q\divides n, \cr
& \qquad
\inpro{\zeta_{2^rp}}G(n,{n\over 2^r},2)=G(n,{n\over 2^r p},2),
\qquad p\divides {n\over 2^r}.
\end{align*}
\end{corollary}

\begin{proof} We firstly observe that the middle index in the first
	group in (\ref{twominrefgps}) is
	$$ {n\over 2^{r+1}}, \quad r=\nu_2({n\over 2})
	\Implies r+1=\nu_2(n),$$
and so we redefine the index $r$ for a more natural formula.

We now consider $G(m,p,2)$ for $p$ and ${m\over p}$ odd ($m=n$), which 
was not treated 
in the proof of Theorem \ref{GL(H^2)classification}, i.e.,
$$ G= G(m,p,2) = G({n\over p},\mD_{n},\mC_{n}),\qquad 
\hbox{$p$ odd, \ ${n\over p}$ odd}; \quad \hbox{($n$ odd)}, $$
	for which minimising ${n\over p}$ gives
\begin{equation}
\label{minrefgps-odd}
G(1,\mD_{n},\mC_{n}) = D_n=G(n,n,2),\qquad 
\hbox{$p=n$, \ ${n\over p}=1$}; \quad \hbox{($n$ odd)}. 
\end{equation}
Thus (\ref{twominrefgps}) and (\ref{minrefgps-odd}) 
give all of the minimal reflection groups for $D_n$. 
The lattice of reflection groups between these and the maximal group $G(2n,2,2)$,
together with their indices and centres, 
is given in Figure \ref{G(2n,2,2)-coll-odd} ($n$ odd)
	and Figure \ref{G(2n,2,2)-coll-even} ($n$ even).
	We note, in particular, that  $G(n,n,2)$ and $G(n,{n\over 2^r},2)$ 
	appear twice, and those groups that evolve from them appear with 
	multiplicity one or two, as shown in the Figures.
\end{proof}

\begin{figure}[!htb]   
\begin{center}
\begin{tikzpicture}
    \matrix (A) [matrix of nodes, row sep=0.8cm, column sep = 1.0 cm]
    {
	    &$\mat{G^{(2,2,n)},\inpro{\zeta_{2n}}\cr G(2n,2,2)}$ & \\
	$\mat{G^{(1,2,n)},\inpro{\zeta_{n}}\cr G(n,1,2)}$ & $\mat{G^{(2,2,q)},\inpro{\zeta_{2q}}\cr G(2n,2{n\over q},2)}$ & $\mat{G^{(2,1,n)},\inpro{\zeta_{n}}\cr G(n,1,2)}$ \\
	$\mat{G^{(1,2,q)},\inpro{\zeta_{q}}\cr G(n,{n\over q},2)}$ & $\mat{G^{(2,2,1)},\inpro{\zeta_{2}}\cr G(2n,2n,2)}$ & $\mat{G^{(2,1,q)},\inpro{\zeta_{q}}\cr G(n,{n\over q},2)}$ \\
	\greybox{$\mat{G^{(1,2,1)},\inpro{1}\cr G(n,n,2)}$} && \greybox{$\mat{G^{(2,1,1)},\inpro{1}\cr G(n,n,2)}$} \\
    };
    \draw (A-1-2)--(A-2-1);
    \draw (A-1-2)--(A-2-2);
    \draw (A-1-2)--(A-2-3);
    \draw (A-2-2)--(A-3-1);
    \draw (A-2-2)--(A-3-2);
    \draw (A-2-2)--(A-3-3);
    \draw (A-2-1)--(A-3-1);
    \draw (A-2-3)--(A-3-3);
    \draw (A-3-2)--(A-4-1);
    \draw (A-3-2)--(A-4-3);
    \draw (A-3-1)--(A-4-1);
    \draw (A-3-3)--(A-4-3);
\end{tikzpicture}
\end{center}
\vskip -0.7 truecm
\caption{\label{G(2n,2,2)-coll-odd}
The $\tau(2n)$ reflection groups with collineation group
$D_n$ for $n$ odd, as they evolve from the minimal reflection
	group $G(n,n,2)=D_n$ (in grey), where $q\divides n$.
They are given as they appear as collineation preserving normal 
	reflection subgroups of the maximal reflection group $G(2n,2,2)$, 
	with the reflection orbit index 
	and centre given.
	}
\end{figure}
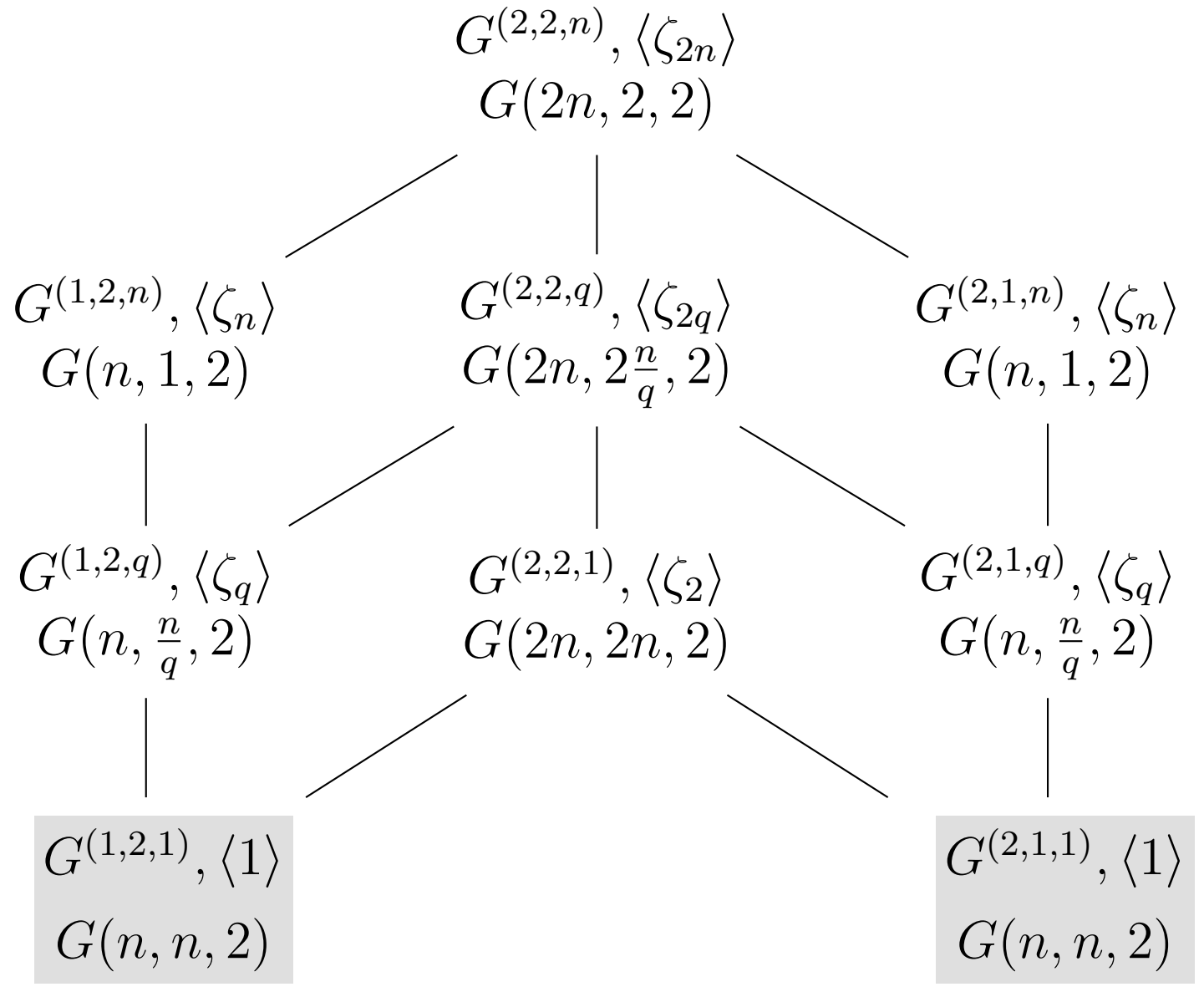
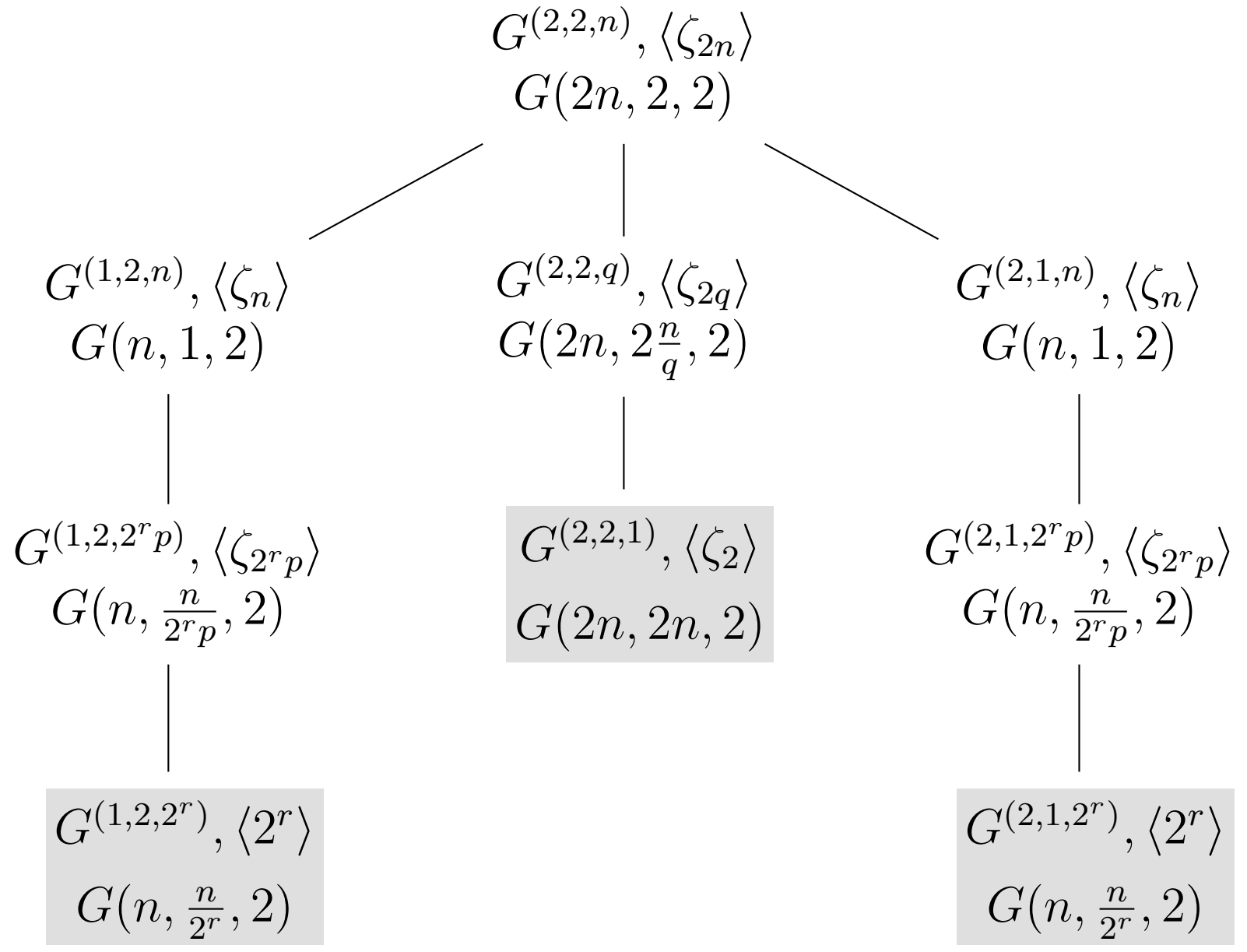
\begin{figure}[!h]
\begin{center}
\begin{tikzpicture}
    \matrix (A) [matrix of nodes, row sep=0.8cm, column sep = 1.0 cm]
    {
    &$\mat{G^{(2,2,n)},\inpro{\zeta_{2n}}\cr G(2n,2,2)}$ & \\
	$\mat{G^{(1,2,n)},\inpro{\zeta_{n}}\cr G(n,1,2)}$ & $\mat{G^{(2,2,q)},\inpro{\zeta_{2q}}\cr G(2n,2{n\over q},2)}$ & $\mat{G^{(2,1,n)},\inpro{\zeta_{n}}\cr G(n,1,2)}$ \\
	$\mat{G^{(1,2,2^r p)},\inpro{\zeta_{2^r p}}\cr G(n,{n\over 2^r p},2)}$ & \greybox{$\mat{G^{(2,2,1)},\inpro{\zeta_{2}}\cr G(2n,2n,2)}$} & $\mat{G^{(2,1,2^r p)},\inpro{\zeta_{2^r p}}\cr G(n,{n\over 2^r p},2)}$ \\
	\greybox{$\mat{G^{(1,2,2^r)},\inpro{2^r}\cr G(n,{n\over 2^r},2)}$} && \greybox{$\mat{G^{(2,1,2^r)},\inpro{2^r}\cr G(n,{n\over 2^r},2)}$} \\
    };
    \draw (A-1-2)--(A-2-1);
    \draw (A-1-2)--(A-2-2);
    \draw (A-1-2)--(A-2-3);
    \draw (A-2-2)--(A-3-2);
    \draw (A-2-1)--(A-3-1);
    \draw (A-2-3)--(A-3-3);
    \draw (A-3-1)--(A-4-1);
    \draw (A-3-3)--(A-4-3);
\end{tikzpicture}
\end{center}
\vskip -0.7 truecm
\caption{\label{G(2n,2,2)-coll-even}
The $\tau(n)+\tau({n\over 2^r})$ reflection groups with collineation group
$D_n$ for $n$ even, as they evolve from the minimal reflection
groups $G(n{n\over2^r},2)$ and $G(2n,2n)$ (in grey), where
	${n\over 2^r}$ is odd ($\nu_2(n)=r$),
	$p\divides {n\over 2^r}$ and $q\divides n$.
They are given as they appear as collineation preserving normal 
	reflection subgroups of the maximal reflection group $G(2n,2,2)$, 
	with the reflection orbit index 
	and centre given.
	For $n=2$, the minimal groups are equal: $G(2,1,2)\cong G(4,4,2)$.
	}
\end{figure}

\begin{example} (Twins)
For $n$ even, there are two minimal groups, and so there can be two 
	collineation preserving normal reflection subgroups of $G(2n,2,2)$ of
	the same order, namely when 
$$ |G(n,{n\over 2^rp},2)|=|G(2n,{2n\over q},2)|
\Iff {2 n^2\over{n\over 2^rp}} = {2 (2n)^2\over {2n\over q}}
\Iff q = 2^{r-1}p.  $$
Indeed, corresponding to each group $G(n,{n\over 2^rp},2)$, there is 
group $G(2n,{2n\over 2^{r-1}p},2)$  of the same order, and these groups are different
	(since $\spec(G)$ differs), except 
	in the case $n=2$ where $G(2,1,2)\cong G(4,4,2)$.
\end{example}

We now seek to identify the imprimitive quaternionic 
reflection groups $G(n,a,b,r)$ within the classification 
of Theorem \ref{GL(H^2)classification}. 
By (\ref{G(n,a,b,r)formula}) and (\ref{Gmp2-G(a,G,H)}), we have
$$ G(n,a,b,r) 
 = \inpro{G({2n\over a},{2n\over ar},2),\zeta_{2n\over b},s}, \qquad
        [n,a,b,r]\in\gL_n,
$$
where 
\begin{equation}
\label{G(2n/a,2n/ar,2)formula}
G({2n\over a},{2n\over ar},2)=
\begin{cases}
G(r,\mD_{2n/a},\mC_{2n/a}), & \hbox{${2n\over ar}$ odd, \ $r$ odd}; \cr
G({r\over 2},\mD_{2n/a},\mD_{n/a}), &\hbox{${2n\over ar}$ odd, \ $r$ even}; \cr
G(r,\mD_{n/a},\mD_{n/a}), & \hbox{${2n\over ar}$ even, \ $r$ even}; \cr
G(r,\mD_{n/a},\mC_{2n/a}), & \hbox{${2n\over ar}$ even, \ $r$ odd}.
\end{cases}
\end{equation}
The groups with $r={n\over ab}$ are called the {\em base groups},
and those with $r={2n\over ab}$ ($ab$ odd), which can be obtained 
from them, by adding a scalar matrix are the {\em higher order groups}.
Now
\begin{align*}
\hbox{base group:} &\qquad && r={n\over ab} 
	\Implies {2n\over ar}=2b \ \hbox{is even}, \cr
\hbox{higher order group:} &\qquad && r={2n\over ab} \Implies {2n\over ar}=b \ 
	\hbox{is odd}, \quad r=2{n\over ab} \ \hbox{is even},
\end{align*}
and so there are only three cases that appear in (\ref{G(2n/a,2n/ar,2)formula}), i.e.,
$$ G({2n\over a},{2n\over ar},2)=
\begin{cases}
	G({r\over 2},\mD_{2n/a},\mD_{n/a}), & 2n=abr, \ \hbox{$r$ even}; \cr
	G(r,\mD_{n/a},\mD_{n/a}), & n=abr, \ \hbox{$r$ even}; \cr
	G(r,\mD_{n/a},\mC_{2n/a}), & n=abr, \ \hbox{$r$ odd}.
\end{cases} $$
which we can write as
\begin{equation}
\label{G(2n/a,2n/ar,2)formulaII}
G({2n\over a},{2n\over ar},2)=
\begin{cases}
G({n\over ab},\mD_{br},\mD_{br/2}), & 2n=abr; \cr
G({n\over ab},\mD_{br},\mD_{br}), & n=abr, \ \hbox{$r$ even}; \cr
G({n\over ab},\mD_{br},\mC_{2br}), & n=abr, \ \hbox{$r$ odd}.
\end{cases}
\end{equation}
Therefore
$$ \inpro{G({2n\over a},{2n\over ar},2),\zeta_{2n\over b}} =
G({n\over b},\mD_{br},\mH),
$$
with $\mH$ as in (\ref{G(2n/a,2n/ar,2)formulaII}),  
except when the scalar $-1\in SL(\CC^2)$ is
introduced, which can happen in the one case
$$ 
\hbox{$r={n\over ab}$ is odd}, \qquad 
\hbox{$ra={n\over b}$ is even}, $$
for which $\mC_{2br}$ enlarges to $\inpro{-1}\mC_{2br}=\mD_{br}$.
Thus, we obtain
\begin{equation}
\label{G(2n/a,2n/ar,2)formulaIII}
\inpro{G({2n\over a},{2n\over ar},2),\zeta_{2n\over b}} =
\begin{cases}
G({n\over b},\mD_{br},\mD_{br/2}), & 2n=abr; \cr
G({n\over b},\mD_{br},\mD_{br}), & n=abr, \ \hbox{$ar$ even}; \cr
G({n\over b},\mD_{br},\mC_{2br}), & n=abr, \ \hbox{$ar$ odd}.
\end{cases}
\end{equation}

We finish this section with a technical calculation of
the somewhat complicated correspondence between the $G(n,a,b,r)$ and $G_\HH(a,\mD_n,\mH)$
quaternionic reflection groups.

\begin{proposition} (Correspondence)
The imprimitive quaternionic reflection groups $G(n,a,b,r)$ can be written in
terms of those in (\ref{GH(a,H,G)imprimrefgroups}) as follows
\begin{equation}
\label{G(2n/a,2n/ar,2)formulaIV}
G(n,a,b,r) =
\begin{cases}
G_\HH({n\over b},\mD_{br},\mD_{br}), & n=abr, \ \hbox{${n\over b}$ even}; \cr
G_\HH({n\over b},\mD_{br},\mC_{2br}), & n=abr, \ \hbox{${n\over b}$ odd}; \cr
G_\HH({n\over b},\mD_{br},\mD_{br/2}), & 2n=abr,
\end{cases}
\end{equation}
and conversely those of (\ref{GH(a,H,G)imprimrefgroups}) as
\begin{align*}
G_\HH(2a,\mD_n,\mD_n)
	&=G({2an\over r},{2a\over r},{n\over r},r), && \quad r:=\gcd(2a,n), \
	\hbox{($n$ even)}, \cr
G_\HH(b,\mD_n,\mC_{2n})
	&=G({bn\over r},{b\over r},{n\over r},r), && \quad r:=\gcd(b,n), 
	\cr
G_\HH(2^r b,\mD_n,\mD_{n/2})
	&=G({bn\over 2q},{b\over q},{n\over 2^{r+1}q},2^{r+1}q), && 
	\quad q:=\gcd(b,n), \ \hbox{$r:=\nu_2({n\over 2})$}, \ \hbox{($n$ even)}, 
\end{align*}
where ${b\over q}$ and ${n\over 2^{r+1}q}$ should be swapped 
	when $2^{r+1}b>n$, so as
to obtain an index in in 
	$\gL_{bn\over 2q}$.
\end{proposition}

\begin{proof}
Since
$$ G(n,a,b,r) = \inpro{G({2n\over a},{2n\over ar},2),\zeta_{2n\over b},s}, 
\qquad [n,a,b,r]\in\gL_n,$$
we obtain (\ref{G(2n/a,2n/ar,2)formulaIV}) from (\ref{G(2n/a,2n/ar,2)formulaIII}), 
where we write $ar$ as ${n\over b}$ and change the ordering of cases for aesthetic purposes.

For the converse formulas, we solve for the parameters in the
three separate cases, using the fact that 
	the fraction
	${a\over b}$ is reduced, since $\gcd(a,b)=1$.
For the $(\mD_n,\mD_n)$ case, we write
$$ G(2A,\cD_N,\cD_N)=G_\HH({n\over b},\mD_{br},\mD_{br}), \qquad n=abr, \quad
	\hbox{${n\over b}$ even}, \quad \hbox{($N$ even)}, $$
and equating coefficients gives
\begin{align*} 
{n\over b}=2A, \quad br={n\over a}=N  &\Implies
{a\over b}={2A\over N} 
	={{2A\over\gcd(2A,N)}\over{N\over\gcd(2A,N)}} 
	\cr
&\Implies a={2A\over\gcd(2A,N)}, \quad b={N\over\gcd(2A,N)}, \cr
&\Implies r={N\over b}= \gcd(2A,N), \quad n=abr = {2AN\over\gcd(2A,N)},
\end{align*}
which gives the first equation, i.e., 
$$ G_\HH(2a,\mD_n,\mD_n)
=G({2an\over r},{2a\over r},{n\over r},r), \qquad r:=\gcd(2a,n), \quad
	\hbox{($n$ even)},  $$
where ${2a\over r}\le {n\over r}$.
Similarly, for the $(\mD_n,\mC_{2n})$ case
$$ G_\HH(B,\cD_N,\mC_{2N})=G_\HH({n\over b},\mD_{br},\mC_{2br}), \qquad n=abr, 
\quad \hbox{${n\over b}$ odd}, \quad \hbox{($B$ odd)}, $$
	gives (replace $2A$ by $B$ in the previous argument)
\begin{align*}
{n\over b}=B, \quad br={n\over a}=N 
	&\Implies n={BN\over r}, \quad a={B\over r}, \quad b={N\over r}, 
	\quad r = \gcd(B,N),
\end{align*}
and we obtain the second equation.

For the higher order groups, the calculations are slightly more involved.
Let 
$$ G_\HH(2^R B,\mD_N,\mD_{N/2})=G_\HH({n\over b},\mD_{br},\mD_{br/2}),
\qquad \hbox{$B$ odd}, \quad \hbox{$N$ even}, \quad \hbox{$R=\nu_2(N/2)$}. $$
	Since $(N/2)/2^R=N/2^{R+1}$ is odd, we have
\begin{align*}
{n\over b} =2^R B, \quad br={2n\over a}=N 
&\Implies {a\over b} 
= {2^R B\over N/2} 
	= {B/q\over (N/2^{R+1})/q}, \qquad q=\gcd(B,N), \cr
&\Implies a={B\over q}, \quad b={N\over 2^{R+1}q}, \quad
n 
={BN\over 2q}, \quad
r 
	= 2^{R+1}q,
\end{align*}
and we have
	$$ G_\HH(2^r b,\mD_n,\mD_{n/2})
=G({bn\over 2q},{b\over q},{n\over 2^{r+1}q},2^{r+1}q), \qquad q=\gcd(b,n), \quad
\hbox{$r=\nu_2({n\over 2})$}. $$
To ensure that the index is $\gL_{bn\over 2q}$, 
we must exchange ${b\over q}$ and ${n\over 2^{r+1}q}$
when ${b\over q}>{n\over 2^{r+1}q}$, to obtain
the index 
$[{bn\over 2q},{\min(2^{r+1}b,n)\over 2^{r+1} q},{\max(2^{r+1}b,n)\over 2^{r+1}q},2^{r+1}q]\in\gL_{bn\over 2q}$.
\end{proof}



\vfil\eject

\section{Reflection orbits and normal reflection subgroups}

We now consider the reflection orbits of the rank two quaternionic 
reflection groups 
$$ \cH_a(G) = \inpro{G,\zeta_m,s}, \qquad m:=|Z(G)|a, \quad \hat G=G, $$
where $G$ is a finite irreducible unitary complex reflection group with $\hat G=G$.
For any reflection group, the reflections with a given root
form a (parabolic) subgroup, and the group acts on these subgroups by conjugation.
We call the orbits the ``reflection orbits'' of the reflection group,
and use the symbolic notation (see \cite{W26c})
$$ n_1 H_{1},\ldots, n_\ell H_\ell \qquad\hbox{(reflection type)} $$
to indicate that there are $\ell$ orbits of sizes $n_1,\ldots,n_\ell$ with 
the reflection subgroup for a given root being isomorphic 
to the abstract groups $H_1,\ldots,H_\ell$. In the case of complex 
reflection groups, and hence the quaternionic reflection groups that we
consider (which have additional reflections of order $2$), the 
groups $H_1,\ldots,H_\ell$ are all cyclic, and the notation
$$ (k_1,\ldots,k_\ell) \qquad\hbox{(reflection type)}  $$ 
is adopted, where $k_j$ is the order of the cyclic 
group $H_j$. The normal reflection subgroups of a given reflection 
group are obtained by taking the reflections corresponding to taking subgroups 
$C_{\ga_j}$ of the $H_j=C_{k_j}$, i.e., divisors $(\ga_1,\ldots,\ga_\ell)$ 
of $(k_1,\ldots,k_\ell)$ (see Lemma 2.1 of \cite{W26c}).
This labelling of the normal reflection subgroups is $1$--$1$ 
for complex reflection groups (Theorem 2.1 of \cite{W26c}),
but may not be so for quaternionic reflection groups, as 
we will see in this section.
We recall that the reflections of $\cH_a(G)$ have two types
\begin{itemize}
	\item The complex reflections in $\inpro{G,\zeta_m}$. (type $A$)
	\item The $m$ quaternionic reflections $\xi s$, where $\xi\in\inpro{\zeta_m}$. (type $Bj$)
\end{itemize}

We now consider the orbit of $s$ in $\cH_a(G)$. 
For any $g\in\inpro{G,\zeta_m}$, by (\ref{gsconjrelation}), we have 
\begin{align*}
s^g &= gsg^{-1}=g(sg^{-1}s)s=g\det(g)g^{-1}s=\det(g)s, \cr
s^{gs}&=gs(s)(gs)^{-1}=gsg^{-1}=s^g=\det(g)s,
\end{align*}
so that the orbit of $s$ is
\begin{equation}
\label{s-orbits}
s^{\cH_a(G)} = \det(\inpro{G,\zeta_m}) s.
\end{equation}
This orbit contains at least the half of the type $Bj$ reflections given
by $\inpro{\zeta_m^2}s$, and the orbit of $\zeta_m s$ contains (at least)
the other half. Thus, there is either one
$m C_2$ reflection orbit of type $Bj$ reflections or 
two ${m\over 2}C_2, {m\over 2}C_2$ orbits, as
determined by calculating (\ref{s-orbits}).

\begin{lemma} 
\label{typeBorbits}
The number of reflection orbits of quaternionic reflections
(type $Bj$) in the groups of Theorems \ref{TOIfamilyclass} and \ref{G(n,a,b,r)familyclass} are as follows
\begin{align*}
\hbox{One orbit:} &\qquad \hbox{$\cO_b(G_8)$ and $\cO_1(G_{12})$}, \quad
	\cr
&\qquad \hbox{$G(n,a,b,r)$, ${n\over b}$ odd or $\nu_2(r)=\nu_2({2n\over b})\ge2$}, \cr
\hbox{Two orbits:} &\qquad \hbox{$\cT_a(G_5)$, $\cO_a(G_{12})$ with $a\ne1$, and $\cI_a(\mI)$}, \cr
&\qquad    \hbox{$G(n,a,b,r)$, ${n\over b}$ even and $\nu_2(r)\ne\nu_2({2n\over b})$}.
\end{align*}
\end{lemma}

\begin{proof}
For each of the groups, we calculate (\ref{s-orbits}), with there
being one orbit if $\det(\inpro{G,\zeta_m})=\inpro{\zeta_m}$ and 
otherwise two. We observe that the determinant of a reflection of
order $k$ is a primitive $k$--th root of unity. 
Since $G$ is generated by reflections, and the determinant is
	multiplicative, we have $\det(G)=\inpro{\zeta_t}$, where $t$ is
the least common multiple of the orders of the reflections in $G$,
and hence
$$ \det(\inpro{G,\zeta_m})=\det(G)\inpro{\zeta_m^2}
=\inpro{\zeta_t}\inpro{\zeta_m^2}. $$
For the primitive groups, we calculate
\begin{align*}
\cT_a(G_5) &:\quad \inpro{\zeta_t}\inpro{\zeta_m^2}=\inpro{\zeta_3}\inpro{\zeta_{6a}^2} =\inpro{\zeta_{3a}}\ne \inpro{\zeta_{6a}}, \cr
\cO_a(G_{12}) &:\quad \inpro{\zeta_t}\inpro{\zeta_m^2}
=\inpro{\zeta_2}\inpro{\zeta_{2a}^2} =\inpro{\zeta_{a}}\ne \inpro{\zeta_{2a}} \iff a\ne 1, \cr
\cO_b(G_{8}) &:\quad \inpro{\zeta_t}\inpro{\zeta_m^2}
	=\inpro{\zeta_4}\inpro{\zeta_{4b}^2} =\inpro{\zeta_4}\inpro{\zeta_{2b}}=\inpro{\zeta_{4b}}, \quad \hbox{(since $b$ is odd)}, \cr
\cI_a(\mI) &:\quad \inpro{\zeta_t}\inpro{\zeta_m^2}
=\inpro{1}\inpro{\zeta_{2a}^2} = \inpro{\zeta_{a}} \ne \inpro{\zeta_{2a}}.
\end{align*}
For the imprimitive groups $G(n,a,b,r)$ generated by the reflections
given in (\ref{G(n,a,b,r)formula}), we have $r\divides m= {2n\over b}$,
$t=r$ for $r$ even, and $t=2r$ for $r$ odd. If $r$ is even and $r$
does not divide ${n\over b}$ then $\nu_2(r)=\nu_2({2n\over b})=k$. 
Accordingly, we calculate the cases
\begin{align*}
	\hbox{$r$ even, $r\divides {n\over b}$:} & \quad
\inpro{\zeta_r}\inpro{\zeta_m^2}
        =(\inpro{-1}\inpro{\zeta_r})\inpro{\zeta_{2n\over b}^2}
        = \inpro{-1}\inpro{\zeta_{n\over b}}
        = \inpro{\zeta_{2n\over b}} \iff \hbox{${n\over b}$ is odd},
	\cr
\hbox{$r$ even, 
	$r\ndivides {n\over b}$:} &\quad
\inpro{\zeta_r}\inpro{\zeta_m^2} 
=\inpro{\zeta_{2^k}} \inpro{\zeta_{r\over 2^k}} \inpro{\zeta_{m\over 2}} =\inpro{\zeta_{2^k}}\inpro{\zeta_m}
	= \inpro{\zeta_m}, \cr
	\hbox{$r$ odd:} &\quad
	\inpro{\zeta_{2r}}\inpro{\zeta_m^2}
	=\inpro{\zeta_{2r}}\inpro{\zeta_{n\over b}}
	=\inpro{\zeta_{2}}\inpro{\zeta_{n\over b}}
        = \inpro{\zeta_{2n\over b}} \iff \hbox{${n\over b}$ is odd}.
\end{align*}
The middle condition above for there to be a single orbit is that
$$ \hbox{$\nu_2(r)=\nu_2({2n\over b})\ge 1$} \Iff
\hbox{$r$ is even and ${n\over b}$ is odd, 
or $\nu_2(r)=\nu_2({2n\over b})\ge 2$}, $$
and we can combine the cases as in the theorem, which completes the proof.
\end{proof}

The result for the $G(n,a,b,r)$ is consistent with that of (ii) in Lemma 5.3 
of \cite{W25}, which was stated and proved at the level of reflection systems.

We now consider the orbits of the complex reflections $g\in G$ in $\cH_a(G)$.
Since conjugation by $\zeta_m$ fixes $g$, its $\cH_a(G)$-orbit does not depend on the
value of $a$. Moreover, the $\cH_a(G)$-orbit of $g$ contains its $G$-orbit,
and can only be larger when $s$ conjugates $g$ into another $G$-orbit.
We therefore consider when (not often) this can happen.

\begin{lemma}
\label{typeAorbits}
The reflection orbits of all the irreducible unitary complex 
reflection groups $G\subset U(\CC^2)$ are closed under conjugation by $s$, 
except for the groups
\begin{align*}
	& \hbox{$G_5$ and $G_7$}, &&\ \hbox{(the $4C_3,4C_3$ orbits combine to a $8C_3$ orbit)}, \cr
	& G(m,p,2), \ 
	\gcd(m,4)=\gcd(p,4)=2, 
	&&\ 
	\hbox{(the ${m\over 2}C_2,{m\over2}C_2$ orbits combine to a $mC_2$ orbit)}.
\end{align*}
\end{lemma}

\begin{proof}
We consider when $s$ conjugates elements of one $n_j C_{k_j}$ 
reflection orbit into
another. Clearly, this could only happen when there are two orbits
	with reflections for a given root of the same order, i.e., the
	indices ${k_j}$ are the same.

For $G_5=\inpro{Z,Z^R}$ and the larger group $G_7$, there are two 
$4C_3$ reflection orbits given by
$$ O_Z=\{(Z^a)^g\}_{0\le a\le 2,g\in G_5},
	\qquad O_{Z^R}=\{((Z^R)^a)^g\}_{0\le a\le 2,g\in G_5}. $$
Since 
$$ ((Z^a)^g)^s =\det( (gZ^a g^{-1})^*) gZ^a g^{-1}
=g(\det(Z^*)Z)^a g^{-1} =g (Z^s)^a g^{-1}, $$
it follows that if $s$ conjugates $Z$ to an element of $O_{Z^R}$
	then it conjugates $O_Z$ to $O_{Z^R}$, and this is indeed the case, since
$$ Z^s = \det(Z^*)Z = g (Z^R)^{-1} g^{-1}, \qquad g=(Z^2 Z^R)^2\in G_5. $$
The only other primitive groups where there are reflection orbits
with the same orders of reflections are the $12C_2,6C_2$ reflection 
orbits of $G_{13}$. However, these come from $12C_2,6C_4$ reflection orbits 
of the larger group $G_{11}$, which are fixed by $s$, and so 
can be disregarded.

For the imprimitive reflection groups $G(m,p,2)$, the reflection orbits were
calculated in Lemma 8.1 of \cite{W26c} (there can be one, two or three).
The relevant cases are for $p$ (and hence $m$) even, where there 
are ${m\over 2}C_2,{m\over 2}C_2$ reflection orbits, given by
$$ O_1 =
\bigl\{ \pmat{0&\go^{2j}\cr\overline{\go}^{2j}&0}\bigr\}_{1\le j\le{m\over2}} ,\qquad
O_2 =
\bigl\{ \pmat{0&\go^{2j-1}\cr\overline{\go}^{2j-1}&0}\bigr\}_{1\le j\le{m\over2}},
	\qquad \go=\zeta_m. $$
and when ${m\over p}=2$ there is a third $2C_2$  reflection orbit
	$$ O_3=\bigl\{ \pmat{-1&0\cr0&1},\pmat{1&0\cr0&1} \bigr\}. $$
Since all these reflections have determinant $-1$, $O_3$ is closed under
conjugation by $s$, and 
$$ (O_1)^s =
\bigl\{ -\pmat{0&\go^{2j}\cr\overline{\go}^{2j}&0}\bigr\}_{1\le j\le{m\over2}}, \qquad
O_2=\bigl\{ \pmat{0&\go^{2j-1}\cr\overline{\go}^{2j-1}&0}\bigr\}_{1\le j\le{m\over2}},
$$
so that $s$ conjugates $O_1$ to $O_2$ to give a single $mC_2$ orbit
	if and only if ${m\over 2}$ is odd.
\end{proof}

For clarity, we modify the labelling of the normal reflection subgroups of 
$\cH_a(G)$ to
$$ (\ga_1,\ldots,\ga_\ell:\gb_1,\ldots,\gb_{\ell'}), \quad \ell'\in\{1,2\} \qquad\hbox{(reflection type)}, $$
where the ``:'' separates the type $A$ (complex) and type $Bj$ (quaternionic) reflections.
For $a>1$, the addition of the scalars $\inpro{\zeta_m}$ to $\cH_a(G)=\inpro{G,\zeta_m,s}$ may introduce
additional complex reflections (previously hidden), and so it is 
convenient to denote the largest complex reflection group contained in 
$\inpro{G,\zeta_m}$ by
$$ G_a:=\hbox{the reflection group generated by the reflections in  
$\inpro{G,\zeta_m}\subset \cH_a(G)$}. $$
Clearly, $G_1=G$, $G\subset G_a$, and $G_a$ is contained in the maximal
reflection group for $G$.

If $\cH_a(G)$ has a single orbit of quaternionic reflections, 
then its normal subgroups are relatively easy to describe.

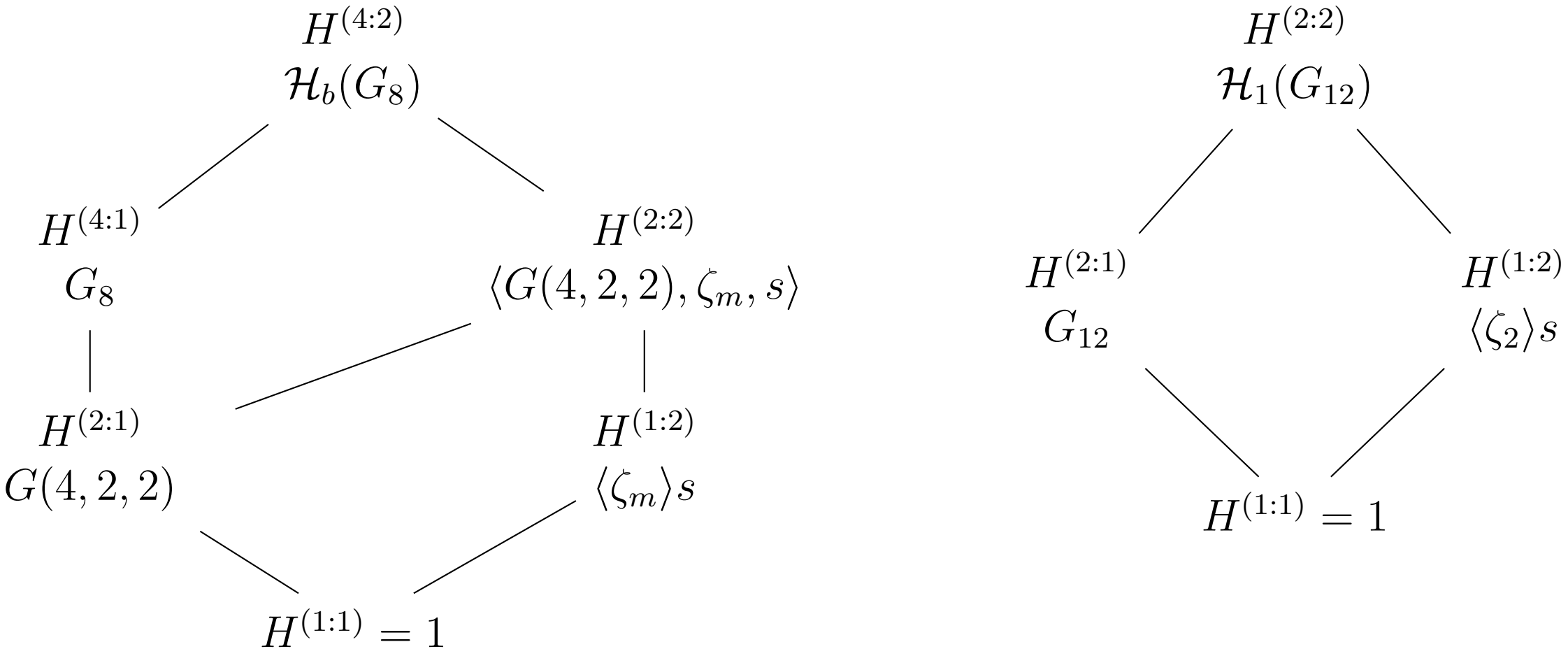
\begin{figure}[H] 
\centering
\begin{tikzpicture}[baseline=(current bounding box.north)]
    \matrix (A) [matrix of nodes, row sep=0.6cm, column sep = 0.1 cm]
    {
            & $\mat{H^{(4:2)}\cr\cH_b(G_8)}$ & \\
	$\mat{H^{(4:1)}\cr G_{8}}$ & & $\mat{H^{(2:2)}\cr \inpro{G(4,2,2),\zeta_m,s}}$ \\
        $\quad \mat{H^{(2:1)}\cr G(4,2,2)} \quad$ & & $\mat{H^{(1:2)}\cr \inpro{\zeta_m}s}$ \\
        & $H^{(1:1)}=1$  & \\
    };
    \draw (A-1-2)--(A-2-1);
    \draw (A-1-2)--(A-2-3);
    \draw (A-2-1)--(A-3-1);
    \draw (A-2-3)--(A-3-3);
    \draw (A-3-1)--(A-4-2);
    \draw (A-3-3)--(A-4-2);
    \draw (A-2-3)--(A-3-1);
\end{tikzpicture}%
\hfill
\begin{tikzpicture}[baseline=(current bounding box.north)]
    \matrix (A) [matrix of nodes, row sep=1.0cm, column sep = 0.4 cm]
    {
            & $\mat{H^{(2:2)}\cr\cH_1(G_{12})}$ & \\
        $\mat{H^{(2:1)}\cr G_{12}}$ & & $\mat{H^{(1:2)}\cr \inpro{\zeta_2}s}$ \\
        & $H^{(1:1)}=1$  & \\
    };
    \draw (A-1-2)--(A-2-1);
    \draw (A-1-2)--(A-2-3);
    \draw (A-2-1)--(A-3-2);
    \draw (A-2-3)--(A-3-2);
\end{tikzpicture}
\vskip-0.2truecm
\caption{\label{Ob(G_8andO1(G12)-lat} 
The lattice of the normal reflection subgroups of
$\cO_b(G_8)=\cH_b(G_8)=H^{(4:2)}$, $b$ odd 
	(note the lattice size is independent of $b$), and $\cO_1(G_{12})=\cH_1(G_{12})$.
        }
\end{figure}

\begin{figure}[!h]    
\centering
\begin{tikzpicture}[baseline=(current bounding box.north)]
    \matrix (A) [matrix of nodes, row sep=0.6cm, column sep = 0.1 cm]
    {
	    &&& $\mat{H^{(4,3:2)}\cr\cH_b(G_8)}$ && \\
	& $\mat{H^{(4,3:1)}\cr G_{10}}$ & $\mat{H^{(4,1:2) }\cr\cH_b(G_8)}$ && $\mat{H^{(2,3:2)}\cr\inpro{G_7,\zeta_m,s}}$ & \\
	$\mat{H^{(4,1:1)}\cr G_8}$ && $\mat{H^{(2,3:1)}\cr G_7}$ & $\mat{H^{(2,1:2)}\cr\inpro{G(4,2,2),\zeta_m,s}}$ && $\mat{H^{(1,3:2)}\cr\inpro{G_5,\zeta_m,s}}$  \\
	& $\mat{H^{(2,1:1)}\cr G(4,2,2)}$ && $\mat{H^{(1,3:1)}\cr G_5}$ & $\mat{H^{(1,1:2)}\cr\inpro{\zeta_m,s}}$ & \\
	&& $\mat{H^{(1,1:1)}\cr 1}$ &&& \\
    };
    \draw (A-1-4)--(A-2-2);
    \draw (A-1-4)--(A-2-3);
    \draw (A-1-4)--(A-2-5);
    \draw (A-2-2)--(A-3-1);
    \draw (A-2-2)--(A-3-3);
    \draw (A-2-3)--(A-3-4);
    \draw (A-2-5)--(A-3-4);
    \draw (A-2-5)--(A-3-6);
    \draw (A-3-1)--(A-4-2);
    \draw (A-3-3)--(A-4-2);
    \draw (A-3-3)--(A-4-4);
    \draw (A-3-4)--(A-4-5);
    \draw (A-3-6)--(A-4-5);
    \draw (A-4-2)--(A-5-3);
    \draw (A-4-4)--(A-5-3);
    \draw (A-4-5)--(A-5-3);

    \draw (A-2-3)--(A-3-1);
    \draw (A-2-5)--(A-3-3);
    \draw (A-3-4)--(A-4-2);
    \draw (A-3-6)--(A-4-4);
\end{tikzpicture}%
\vskip-0.2truecm
	\caption{\label{Ob(G_8)3dividesb-lat} 
The lattice of the normal reflection subgroups of
$\cO_b(G_8)=\cH_b(G_8)=H^{(4,2:2)}$, $3\divides b$, $b$ odd 
        (note the lattice size is independent of $b$).
        }
\end{figure}
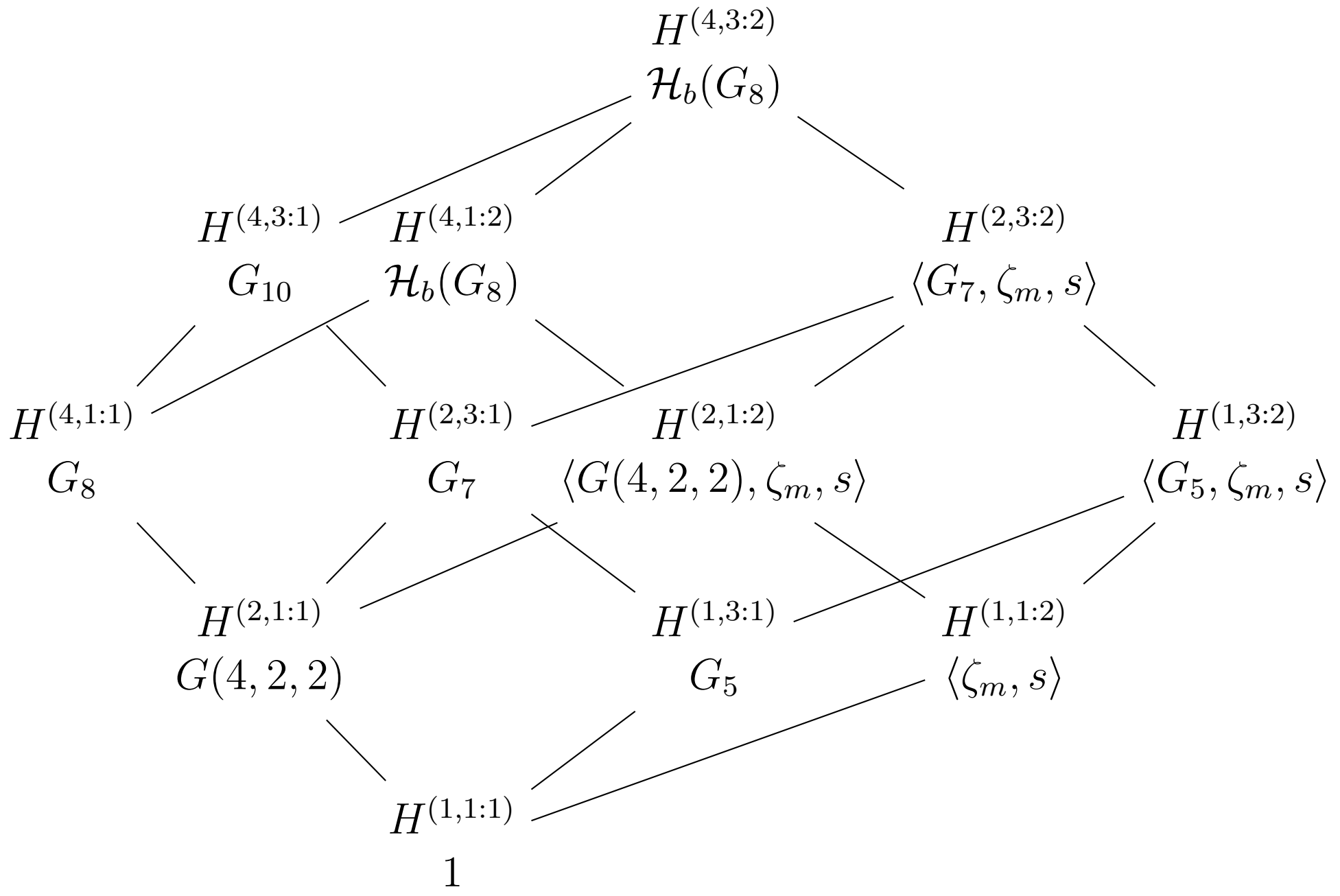

\begin{example}
\label{oneBjorbit}
(One $Bj$ orbit)
Suppose that $\cH_a(G)=H^{(k_1,\ldots,k_\ell:2)}$ has a single orbit
$\inpro{\zeta_m}s$ of type $Bj$ (quaternionic) reflections, i.e., $\ell'=1$.
Then its normal reflection subgroups are
$$ G^{(\ga_1,\ldots,\ga_\ell)}:=H^{(\ga_1,\ldots,\ga_\ell:1)}\lhd G_a, \qquad
H^{(\ga_1,\ldots,\ga_\ell:2)}  
	= \inpro{G^{(\ga_1,\ldots,\ga_\ell)},\zeta_m, s}, $$
where $(\ga_1,\ldots,\ga_\ell)$ (entrywise) divides $(k_1,\ldots,k_\ell)$.
Each $G^{(\ga_1,\ldots,\ga_\ell)}$ is a different normal reflection 
subgroup of $G_a=H^{(k_1,\ldots,k_\ell:1)}$ (uniquely determined by its reflection orbits), and it remains to determine whether any
groups $H^{(\ga_1,\ldots,\ga_\ell:2)}=\inpro{G^{(\ga_1,\ldots,\ga_\ell)},\zeta_m, s}$
	could have multiple labels.

For $G$ primitive, Lemma \ref{typeBorbits} provides just two 
possibilities:  
$$ \cO_b(G_8) 
=\cH_b(G_8), \qquad \cO_1(G_{12}) 
	=\cH_1(G_{12}), $$
where (in the first case)
$$ G_b=\begin{cases}
G_8, & \gcd(b,3)=1; \cr
G_{10}, & b\divides 3.
\end{cases}
$$
The reflection types of $G_8$, $G_{10}$ and $G_{12}$ are
$6C_4$,  $6C_4,8C_3$ and $12C_2$, giving the following 
$H^{(k_1,\ldots,k_\ell:2)}$
$$ H^{(4:2)}, \quad H^{(4,3:2)}, \quad H^{(2:2)}, $$
respectively. Since there are just a few divisors of $(k_1,\ldots,k_\ell:2)$,
we can check whether there are any multiple labels of the
type $(\ga_1,\ldots,\ga_\ell:1)$, and there is one case:
	$$ \cH_b(G_8)=H^{(4,3,:2)}=H^{(4,1,:2)}, \qquad
	\hbox{$3\divides b$}. $$
The lattices of the normal reflection subgroups of $\cO_b(G_8)$
and $\cO_1(G_{12})$ are given in Figures \ref{Ob(G_8andO1(G12)-lat} and \ref{Ob(G_8)3dividesb-lat}.
	From these, we can conclude that if $\cH_a(a)$ has one orbit of
	quaternionic reflections and is primitive, then its only 
	possible normal complex reflection subgroups are $G_8$, $G(4,2,2)$ and
	$G_{12}$.
\end{example}

\begin{example}
\label{TwoBjorbits}
(Two $Bj$ orbits)
Suppose that $\cH_a(G)=H^{(k_1,\ldots,k_\ell:2,2)}$ has two orbits
$\inpro{\zeta_m^2} s$ and $\inpro{\zeta_m^2} \zeta_m s$ of quaternionic reflections, 
i.e., $\ell'=2$.
Then its normal reflection subgroups are
$$ G^{(\ga_1,\ldots,\ga_\ell)}:=H^{(\ga_1,\ldots,\ga_\ell:1,1)}\lhd G, \qquad
H^{(\ga_1,\ldots,\ga_\ell:2,2)}  
        = \inpro{G^{(\ga_1,\ldots,\ga_\ell)},\zeta_m, s}, $$
and the pairs of isomorphic subgroups
$$ H^{(\ga_1,\ldots,\ga_\ell:2,1)}  
= \inpro{G^{(\ga_1,\ldots,\ga_\ell)},\zeta_{m/2}, s}, \qquad
H^{(\ga_1,\ldots,\ga_\ell:1,2)}  
= \inpro{G^{(\ga_1,\ldots,\ga_\ell)},\zeta_{m/2}, \zeta_m s}.  $$
The above structural form of the pairs follows from the calculations
$$ (\zeta_m^2 s)s^{-1}=\zeta_m^2=\zeta_{m/2}, \qquad
	(\zeta_m^2 \zeta_m s)(\zeta_m s)^{-1}=\zeta_m^2=\zeta_{m/2}. $$
However, this does not preclude the possibility that $\zeta_m$ belongs 
to one, and hence both of the pairs, which would imply the multiple labelling
$$ H^{(\ga_1,\ldots,\ga_\ell:2,2)}
= H^{(\ga_1,\ldots,\ga_\ell:2,1)}
= H^{(\ga_1,\ldots,\ga_\ell:1,2)}. $$
This can happen, e.g., for $\cH_a(G_5)=H^{(3:2,2)}$ with $a$ odd, the centre of $G_5$
has order $6$ and the complex orbit of $\cH_a(G_5)$ generates $G_5$, 
so that $m=6a$ and $\zeta_m=\zeta_2 \zeta_m^2\in H^{(3:2,1)}$, giving
$$ \cH_a(G_5)=H^{(3:2,2)}=H^{(3:2,1)}=H^{(3:1,2)}. $$

On the other hand, if $\xi_m$ does not belong to one (hence both) of the
pairs then they are different groups containing either $s$ or $\zeta_m s$,
but not both, with the particular case
$$ H^{(1,\ldots,1:2,2)} = \inpro{\zeta_{m}}s, \qquad
H^{(1,\ldots,1:1,2)} = \inpro{\zeta_{m\over2}}s, \qquad
H^{(1,\ldots,1:2,1)} = \inpro{\zeta_{m\over2}}\zeta_{m}s. $$
\end{example}

We now give the reflection types for a general quaternionic reflection group
from Theorem \ref{TOIfamilyclass}.
Here we use the notation $ n_j C_{k_j}=[g]$ to indicate that $g$ is
a reflection of order $k_j$ on the $ n_j C_{k_j}$ reflection orbit.


\begin{example}
\label{Ta(G5)-orbits}
($\cT_a(G_5)$ orbits)
The reflection orbits for $\cT_a(G_5)$ are
	$$ 8C_3=[Z] \ \hbox{(a combined orbit)}, \qquad
	6C_2=[R^2] \  \hbox{($a$ even)}, \qquad
	3aC_2=[s],3aC_2=[\zeta_{6a}s], $$
giving $8\cdot 2+6(\gcd(a,2)-1)+6a$ reflections on 
$3+(\gcd(a,2)-1)$ orbits. Here the reflection type of
$G_5$ is $4C_3=[Z],4C_3=[Z^R]$, and the orbits combine 
	by Lemma \ref{typeAorbits}.
\end{example}

\begin{example} 
\label{Oa(G12)-orbits}
($\cO_a(G_{12})$ orbits)
The reflection orbits for $\cO_a(G_{12})$ are
$$ 12C_2=[F],\qquad 6C_2=[R^2], \ 2\divides a,\ 4\ndivides a \quad\hbox{or}\quad
6C_4=[R], \ 4\divides a, \qquad 8C_3=[Z], \ 3\divides a, $$
and the quaternionic reflections
$$ aC_2=[s], aC_2=[\zeta_{2a}s] \ \hbox{(for $a$ even)}
	\quad\hbox{or}\quad
2aC_2=[s] \ \hbox{(for $a$ odd)}, $$
for a total number of reflections
$$ 12+6(\gcd(4,a)-1)+8(\gcd(3,a)-1)+2a. $$
The number of reflection orbits is
$$ 1+(\gcd(2,a)-1)+{1\over2}(\gcd(3,a)-1)+\gcd(2,a), $$
and there can be $2,3,4,5$ orbits, e.g.,
$$ \cO_1(G_{12})=H^{(2:2)}, \quad 
	\cO_2(G_{12})=H^{(2,2:2,2)}, \quad 
	\cO_3(G_{12})=H^{(2,3:2)}, \quad 
	\cO_{12}(G_{12})=H^{(2,4,3:2,2)}.
$$  
\end{example}

\begin{example}
\label{Oa(G8)-orbits}
($\cO_b(G_8)$ orbits)
The reflection orbit for $\cO_b(G_8)$ ($b$ odd or even) are
$$ 6C_4=[R], \qquad 8C_3=[Z],\ 3\divides b, \qquad 12C_2=[F],\ 2\divides b, $$
$$ 4bC_2=[s]\ \hbox{($b$ odd)},  \quad 2bC_2=[s],2bC_2=[\zeta_{4b}s]\ \hbox{($b$ even)}, $$
for a total number of reflections 
	$$ 6\cdot 3+8(\gcd(3,b)-1)+12(\gcd(2,b)-1)+4b, $$
on $1 +{1\over2}(\gcd(3,b)-1) +(\gcd(2,b)-1) +\gcd(2,b)$ orbits.
\end{example}

\begin{example}
\label{Ia(I)-orbits}
	($\cI_a(\mI)$ orbits) The reflection orbits
	of $\cH_a(\cI)$ for
	$a=2,3,4,5,6,8,9,\ldots$ 
	($a$ is divisible by $2$, $3$ or $5$) are
$$ 
30C_2=[R^2],\ 2\divides a, \quad
20C_3=[F],\ 3\divides a, \quad
12C_5=[M],\ 5\divides a, \quad
	aC_2=[s],aC_2=[\zeta_{2a}s]
	$$
so the number of reflections is
	$$ 30(\gcd(2,a)-1)+20(\gcd(3,a)-1)+12(\gcd(5,a)-1)+2a, $$
on $3,4,5$ orbits as given by
$$ (\gcd(2,a)-1) +{1\over2}(\gcd(3,a)-1) +{1\over4}(\gcd(5,a)-1) +\gcd(2,a). $$
\end{example}

Curiously, many of these groups can have the same number of reflections, e.g.,
\begin{align*}
\hbox{$46$ reflections:} &\quad \cH_a(G_{12}), \quad a= 6, 8, 9, 14, 17, \cr
\hbox{$70$ reflections:} &\quad \cH_a(G_{12}), \quad a= 12, 18, 20, 21, 26, 29, \cr
\hbox{$94$ reflections:} &\quad \cH_a(G_{12}), \quad a= 24, 30, 32, 33, 38, 41, 
\end{align*}
which leads to the following:

\begin{example} (Sixes) 
The reflection group $\cH_a(G_{12})$ has $70+24k$ reflections for the six values
$$ a\in \{12,18,20,21,26,29\}+12k, \qquad k=0,1,2,\ldots.  $$
\end{example}

%
%
%
%

\begin{table}[H]
\caption{
\label{Gj-orbitsTable}
	The reflection types of groups $\cH_a(G)$ giving the primitive reflection groups
	of Theorem \ref{TOIfamilyclass},
as detailed in Examples \ref{Ta(G5)-orbits}, 
	\ref{Oa(G12)-orbits},
	\ref{Oa(G8)-orbits},
	\ref{Ia(I)-orbits}.
	}
\begin{center}
\begin{tabular}{ |  >{$}l<{$} | >{$}l<{$} | >{$}l<{$} | }
\hline
&&\\[-0.3cm]
\hbox{group} & \hbox{reflection orbits} & \hbox{reflection type} \\[0.1cm]
\hline
 &&\\[-0.3cm]
\cT_a(G_5),\, \hbox{$a$ odd} &  8C_3:3aC_2,3aC_2  & {(3:2,2)} \\
	\cT_a(G_5),\, \hbox{$a$ even} &  8C_3,6C_2:3aC_2,3aC_2  & {(3,2:2,2)}  \\[0.1cm]
\cO_a(G_{12}),\, \gcd(a,12)=1 &  12C_2:2aC_2  & {(2:2)}  \\
\cO_a(G_{12}),\, \gcd(a,12)=2 &  12C_2,6C_2:aC_2,aC_2  & {(2,2:2,2)}   \\
\cO_a(G_{12}),\, \gcd(a,12)=3 &  12C_2,8C_3:2aC_2  & {(2,3:2)}  \\
\cO_a(G_{12}),\, \gcd(a,12)=4 &  12C_2,6C_4:aC_2,aC_2  & {(2,4:2,2)}   \\
\cO_a(G_{12}),\, \gcd(a,12)=6 &  12C_2,6C_2,8C_3:aC_2,aC_2  & {(2,2,3:2,2)}   \\
\cO_a(G_{12}),\, \gcd(a,12)=12 &  12C_2,6C_4,8C_3:aC_2,aC_2  & {(2,4,3:2,2)}   \\[0.1cm]
\cO_b(G_{8}),\, \gcd(b,6)=1 &  6C_4:4bC_2  & {(4:2)}   \\
\cO_b(G_{8}),\, \gcd(b,6)=2 &  6C_4,12C_2:2bC_2,2bC_2  & {(4,2:2,2)}   \\
\cO_b(G_{8}),\, \gcd(b,6)=3 &  6C_4,8C_3:4bC_2  & {(4,3:2)}   \\
\cO_b(G_{8}),\, \gcd(b,6)=6 &  6C_4,12C_2,8C_3:2bC_2,2bC_2  & {(4,2,3:2,2)}   \\[0.1cm]
\cI_a(\mI),\, \gcd(a,30)=2 &  30C_2:aC_2,aC_2  & {(2:2,2)}   \\
\cI_a(\mI),\, \gcd(a,30)=3 &  20C_3:aC_2,aC_2  & {(3:2,2)}   \\
\cI_a(\mI),\, \gcd(a,30)=5 &  12C_5:aC_2,aC_2  & {(5:2,2)}   \\
\cI_a(\mI),\, \gcd(a,30)=6 &  30C_2,20C_3:aC_2,aC_2  & {(2,3:2,2)}   \\
\cI_a(\mI),\, \gcd(a,30)=10 &  30C_2,12C_5:aC_2,aC_2  & {(2,5:2,2)}   \\
\cI_a(\mI),\, \gcd(a,30)=15 &  20C_3,12C_5:aC_2,aC_2  & {(3,5:2,2)}   \\
\cI_a(\mI),\, \gcd(a,30)=30 &  30C_2,20C_3,12C_5:aC_2,aC_2  & {(2,3,5:2,2)}   \\ [0.1cm]
\hline
\end{tabular}
\end{center}
\end{table}

We now consider the imprimitive quaternionic reflection 
groups $G(n,a,b,r)$, which is more involved as there can be
$2$, $3$, $4$, or $5$ reflection orbits. To get an overview, see 
Table \ref{G(n,a,b,r)groups-table}, which gives the numbers of labels and of normal 
reflection subgroups. 

\begin{table}[!h]
\caption{The size of the lattice of labels, i.e., the divisors of $(k_1,\ldots,k_m)$,
and the size of the lattice of normal reflection subgroups of 
	$G(n,a,b,r)=H^{(k_1,\ldots,k_m)}$ for various indices $[n,a,b,r]\in\gL_n$.
	These differ if and only if there are subgroups with multiple labels.}
\label{G(n,a,b,r)groups-table}
\begin{center}
\begin{tabular}{ | >{$}l<{$}  >{$}l<{$}  >{$}l<{$} | >{$}l<{$} >{$}l<{$}  >{$}l<{$} | >{$}l<{$} >{$}l<{$} >{$}l<{$} |  >{$}l<{$} >{$}l<{$} >{$}l<{$}|}
\hline
\multicolumn{3}{|l|}{\small group / lattice sizes} &
\multicolumn{3}{l|}{\small group / lattice sizes} &
\multicolumn{3}{l|}{\small group / lattice sizes} &
\multicolumn{3}{l|}{\small group / lattice sizes} \\
\hline
&&&&&&&&&&&\\[-0.3cm]
[ 2, 1, 2, 1 ] & 8 & 8 & 
[ 3, 1, 1, 6 ] & 16 & 14 &
[ 4, 1, 2, 2 ] & 32 & 23 &
[ 5, 1, 1, 10 ] & 16 & 14 \cr
[ 2, 1, 1, 4 ] & 12 & 11 &
[ 3, 1, 3, 1 ] & 4 & 4 &
[ 4, 1, 1, 8 ] & 16 & 14 &
[ 5, 1, 5, 1 ] & 4 & 4 \cr
[ 2, 1, 1, 2 ] & 32 & 27 &
[ 3, 1, 1, 3 ] & 8 & 7 &
[ 4, 1, 4, 1 ] & 8 & 8 &
[ 5, 1, 5, 2 ] & 8 & 8 \cr
&&&
[ 3, 1, 3, 2 ] & 8 & 8 &
[ 4, 1, 1, 4 ] & 48 & 34 &
	[ 5, 1, 1, 5 ] & 8 & 7 \\[0.1cm]
\hline
&&&&&&&&&&&\\[-0.3cm]
[ 6, 1, 2, 3 ] & 16 & 13 &
[ 7, 1, 7, 2 ] & 8 & 8 & 
[ 8, 1, 2, 4 ] & 48 & 30 &
[ 9, 1, 1, 9 ] & 12 & 10 \cr
[ 6, 1, 1, 6 ] & 64 & 41 &
[ 7, 1, 1, 7 ] & 8 & 7 & 
[ 8, 1, 4, 2 ] & 32 & 23 &
[ 9, 1, 9, 2 ] & 8 & 8 \cr
[ 6, 1, 3, 2 ] & 32 & 27 &
[ 7, 1, 1, 14 ] & 16 & 14 &
[ 8, 1, 1, 16 ] & 20 & 17 &
[ 9, 1, 3, 3 ] & 8 & 7 \cr
[ 6, 2, 3, 1 ] & 8 & 8 &
[ 7, 1, 7, 1 ] & 4 & 4 &
[ 8, 1, 8, 1 ] & 8 & 8 &
[ 9, 1, 1, 18 ] & 24 & 20 \cr
[ 6, 1, 3, 4 ] & 12 & 11 &
	&&&
[ 8, 1, 1, 8 ] & 64 & 41 &
[ 9, 1, 3, 6 ] & 16 & 14 \cr
[ 6, 1, 6, 1 ] & 8 & 8 &
	&&& 
	&&& 
[ 9, 1, 9, 1 ] & 4 & 4 \cr
[ 6, 1, 1, 12 ] & 24 & 20 &
	&&&
	&&& && \\[0.1 cm]
\hline
\end{tabular}
\end{center}
\end{table}

It appears that
\begin{itemize}
\item The subgroups of $G(n,a,b,r)$ can have multiple labels,
	except for when $n=ab$.
\end{itemize}

To understand the mechanism behind this we consider the abelianisation, 
i.e., the surjective homomorphism
$$ G\mapsto {G\over[G,G]}, $$
where $G_{\rm ab}={G\over[G,G]}$ is the largest quotient of $G$ by 
a normal subgroup which is abelian (also called the abelianisation).
This maps each reflection orbit $n_j C_{k_j}$ to a cyclic group of order dividing
$k_j$. For $G$ a complex reflection group, the abelianisation maps
onto $C_{k_1}\times \cdots\times C_{k_\ell}$ (Lemma 7.1, \cite{W26c}), thereby
showing that the normal reflection subgroups of $G=G^{(k_1,\ldots,k_\ell)}$
are uniquely labelled by the divisors of $(k_1,\ldots,k_\ell)$.
The connection of the abelianisation with braid groups is explored in 
 \cite{BMR98}, \cite{B11}.

\begin{theorem}
\label{G(n,a,b,r)-abel-thm}
The abelianisation of $G=G(n,a,b,r)$, $[n,a,b,r]\in\gL_n$, has order
\begin{equation}
\label{G_ab-I}
|G(n,a,b,r)_{\rm ab}| = {16nr\over \lcm(2,{n\over a},r) \lcm(2,{n\over b},r) }
	\in \{4,8,16\},
\end{equation}
and is given by 
\begin{equation}
\label{G_ab-II}
G(n,a,b,r)_{\rm ab}=
\begin{cases}
C_2\times C_2, & \hbox{${n\over a},{n\over b},r$ are odd}; \cr
        C_2\times C_2\times C_2\times C_2, & \hbox{${n\over a},{n\over b},r$ are even, $r\ne {2n\over ab}$}; \cr
C_2\times C_2\times C_2, & \hbox{otherwise}. \cr
\end{cases}
\end{equation}
Therefore,
the number of factors of $C_2$ above is a lower bound on 
the number of reflections of even order from different reflection orbits that
	are required to generate $G(n,a,b,r)$.
\end{theorem}

\begin{proof}
By direct computation, the commutators for the five generating 
reflections for $G$ given in (\ref{G(n,a,b,r)formula}) are seen to generate
the following subgroup of diagonal matrices
\begin{equation}
\label{Gab-gens}
\inpro{\zeta_2,\zeta_{n\over b},\zeta_r,
\pmat{\zeta_{2}&0\cr0&\overline{\zeta_{2}}},
\pmat{\zeta_{n\over a}&0\cr0&\overline{\zeta_{n\over a}}},
\pmat{\zeta_{r}&0\cr0&\overline{\zeta_{r}}} 
	} \subset [G,G],
\end{equation}
which is normal in $G$, and hence is its commutator subgroup $[G,G]$.
For example
$$
[\pmat{0&\go^a\cr\overline{\go^a}&0},\pmat{\zeta_r&0\cr0&1}] 
= \pmat{{\zeta_r}&0\cr0&\overline{\zeta_r}}, \qquad
[\pmat{0&\go^b j \cr-\go^b j&0},\pmat{\zeta_r&0\cr0&1}]
= \pmat{{\zeta_r}&0\cr0&{\zeta_r}}.  $$	
The first three scalar matrices in (\ref{Gab-gens}) generate a cyclic
group of order $\lcm(2,{n\over b},r)$, and, similarly, the last three 
matrices (all of determinant $1$) generate a cyclic group of order
$\lcm(2,{n\over a},r)$. The intersection of these two cyclic groups
is the scalar matrices of determinant $1$, i.e., the group of order
two generated by $\zeta_2=-1$. Thus, we have
$$ |G(n,a,b,r)_{\rm ab}| = {|G|\over |[G,G]|}
= {8nr\over {1\over2}\lcm(2,{n\over a},r)\lcm(2,{n\over b},r)}
= {16nr\over \lcm(2,{n\over a},r)\lcm(2,{n\over b},r)}.  $$
The minimum and maximum values of $|G_{\rm ab}|$ are easily
	seen to be taken at ${n\over a},{n\over b},r$ all odd and all even
(respectively). Indeed, in these cases we have
$$ 
4={16nr\over 2{n\over a}\cdot 2{n\over b}} \hbox{ divides } |G_{\rm ab}|, \qquad
|G_{\rm ab}| \hbox{ divides } {16nr\over {n\over a}\cdot {n\over b}} =16,
$$
and so one obtains to orders of (\ref{G_ab-I}) and (\ref{G_ab-II}).

The abelianisation maps a reflection and its orbit to an element dividing 
the order of the reflection. Therefore the images of the reflections of
order $2$ have order $1$ or $2$. The reflection of order $r$ must have an
image which is a power of $2$, and so if $r$ is odd the image has order
$1$, and if $r$ is even, then the square of the reflection is
$$ \pmat{\zeta_r^2&0\cr0&1} = \pmat{{\zeta_r}&0\cr0&{\zeta_r}}
	\pmat{{\zeta_r}&0\cr0&\overline{\zeta_r}} \in [G,G], $$
so its image also has order $1$ or $2$.
Since the images of the generating reflection generates the 
abelianisation, we therefore conclude that it is a product cyclic
groups of order $2$.

The abelianisation has minimal generating sets of size $2,4,3$, for 
the cases of (\ref{G_ab-II}). Therefore, a generating set of reflections for
$G$, grouped  by reflection orbits, must map to at least $2,4,3$ nontrivial 
elements of the abelianisation (respectively).
\end{proof}


Similar calculations show that
$$ \cH_a(G_5)_{\rm ab} =
\begin{cases}
C_2\times C_3, & \hbox{$a$ odd}; \cr
C_2\times C_2\times C_3, & \hbox{$a$ even},
\end{cases}, $$
$$ \cH_a(G_8)_{\rm ab}=\cH_a(G_{12})_{\rm ab} =
\begin{cases}
C_2\times C_2, & \hbox{$a$ odd}; \cr
C_2\times C_2\times C_2, & \hbox{$a$ even}.
\end{cases} $$

The reflection orbits of 
$$ G(n,a,b,r)=\cH_a(G({2n\over a},{2n\over ar},2)) $$
can be calculated by 
\begin{itemize}
\item Taking the (type $A$) reflections of $G({2n\over a},{2n\over ar},2)$
	given in Lemma 8.1 of \cite{W26c},
	and checking if any they combine under the action of $s$
		(Lemma \ref{typeAorbits}).
\item 
Taking the one or two type $Bj$ orbits as given by Lemma \ref{typeBorbits} .
\end{itemize}
In this way, we obtain the following.

\begin{table}[H]
\caption{
\label{G(n,a,b,r)-orbitsTable}
The reflection types of the imprimitive reflection groups 
	$G(n,a,b,r)$.  }
\begin{center}
\begin{tabular}{ |  >{$}l<{$} | >{$}l<{$} | >{$}l<{$} | >{$}l<{$} | }
\hline
        &&&\\[-0.3cm]
	\hbox{orbits} & \hbox{condition} & \hbox{reflection type} & \hbox{example} \\[0.1cm]
\hline
  &&&\\[-0.3cm]
\hbox{five orbits} &  n=abr,\ \hbox{$r$ even ($n$ even)} & {(2,2,r:2,2)} 
	& [ 2, 1, 1, 2 ] \\
\hbox{four orbits} &  \hbox{$n=abr$ even,\ $r\ne1$ odd, ${n\over b}$ odd} & {(2,2,r:2)} & [ 6, 1, 2, 3 ] \\
 &  \hbox{$n=abr$ even,\ $r\ne1$ odd, ${n\over b}$ even} & {(2,r:2,2)} & [ 18, 2, 3, 3 ] \\
\hbox{three orbits} &  \hbox{$n=abr$ even,\ $r=1$, ${n\over b}$ odd} & {(2,2:2)} & [ 2, 1, 2, 1 ] \\
 &  \hbox{$n=abr$ even,\ $r=1$, ${n\over b}$ even} & {(2:2,2)} & [ 6, 2, 3, 1 ] \\
& \hbox{$n=abr$ odd,\ $r\ne 1$ ($r$ odd)} & {(2,r:2)} & [ 3, 1, 1, 3 ] \\
	& \hbox{$2n=abr$ (higher order group)} & {(2,r:2)} & [ 2, 1, 1, 4 ] \\
	\hbox{two orbits} &  \hbox{$n=abr$ odd,\ $r=1$} & {(2:2)} 
	& [ 3, 1, 3, 1 ] \\[0.1cm]
\hline
\end{tabular}
\end{center}
\end{table}

\begin{theorem}
The number of reflection orbits of the imprimitive quaternionic 
reflection groups
$$ G=G(n,a,b,r)=\cH_a(G({2n\over a},{2n\over ar},2)), 
\qquad [n,a,b,r]\in\gL_n, $$
can be two, three, four or five, as given by the reflection types in 
Table \ref{G(n,a,b,r)-orbitsTable}.
\end{theorem}

\begin{proof}
We consider the case of five orbits. 

\noindent
For $G({2n\over a},{2n\over ar},2)$ to have three orbits, we must have
	(see Lemma 8.1 of \cite{W26c}) 
\begin{equation}
\label{5orbeq1}
\hbox{${2n\over a}\ne{2n\over ar}$}, \quad \hbox{${2n\over ar}$ is even} 
\Iff r\ne1, \quad \hbox{${2n\over ar}$ is even}.
\end{equation}
For these 
	orbits to stay separate under the action of $s$,
by Lemma \ref{typeAorbits}, we must have
\begin{equation}
\label{5orbeq2}
\hbox{$\gcd({2n\over a},4)\ne2$ \quad or \quad $\gcd({2n\over ar},4)\ne2$}.
\end{equation}
Since $\gcd({2n\over a},4)\ge\gcd({2n\over ar},4)\ge2$
	(because ${2n\over ar}$ is even), this condition becomes
$$ \hbox{$\gcd({2n\over a},4)=4$} \Implies \hbox{${n\over a}$ is even}. $$
By Lemma \ref{typeBorbits},
the condition for there to be two $Bj$ orbits is
\begin{equation}
\label{5orbeq3}
\hbox{${n\over b}$ is even}, \quad
\hbox{$\nu_2(r)\ne\nu_2({2n\over b})$}
\Iff \hbox{${n\over b}$ is even}, \quad n=abr.
\end{equation}
The 
	equivalence above follows from the fact that if
	$2n=abr$ (so that $ab$ is odd), then
	$$ \hbox{$\nu_2({2n\over b})=\nu_2(ar)=\nu_2(a)+\nu_2(r)=\nu_2(r)$}. $$ 
Since $\gcd(a,b)=1$, the conditions ${n\over a}$
	and ${n\over b}$ above, imply that
	$r={n\over ab}$ is even.
Hence for $G(n,a,b,r)$ to have five reflection orbits,
we must have $n=abr$ and $r$ is even. Conversely, if
$n=abr$ and $r$ is even, then it is easily checked that
(\ref{5orbeq1}), (\ref{5orbeq2}), (\ref{5orbeq3}) hold, i.e.,
$$ r\ne1, \quad \hbox{${2n\over ar}=2b$ is even}, \qquad
\hbox{$\gcd({2n\over a},4)=\gcd({2br},4)=4\ne 2$}, \qquad
\hbox{${n\over b}=ar$ is even},
$$
so that there are five reflection orbits.

The other cases follow from similar elementary calculations.
\end{proof}

\begin{example} 
\label{Uniquelabellingexample}
	($n=ab$)
From Table \ref{G(n,a,b,r)groups-table}, 
we see that the labelling for
the normal reflection subgroups of $G(n,a,b,r)$ can be unique,
or there can be multiple labels. By considering the cases 
$n$ even ($8$ subgroups) and $n$ odd ($4$ subgroups), it can be
shown that 
$$ G(n,a,b,1)= \cH_1(G({2n\over a},{2n\over a},2))
=\inpro{\pmat{0&1\cr1&0},\pmat{0&\zeta_{2n\over a}\cr\zeta_{2n\over a}^{-1}&0},\zeta_{2n\over b},s} $$
has unique labelling,
	as do the higher order groups $G(n,a,b,2)$ for $n=ab$ odd.
Our calculations suggest these are the only cases of unique labelling
	(for all subgroups).
\end{example}

All of the reflections in the groups of Example \ref{Uniquelabellingexample} 
have order $2$ (and there can be one or two $Bj$ reflection orbits). 
This is not sufficient for unique labelling.

\begin{example}
\label{fivereflectionorbits}
(Five orbits) 
The imprimitive quaternionic reflection group
$$ K=G(2,1,1,2)=\cG(\{1,i,j,k\},\{\})=\inpro{G(4,2,2),\zeta_4,s} $$
has five reflection orbits, and is notable for many reasons:
\begin{itemize}
\item It is the unique quaternionic reflection group 
with the maximal number of five systems of imprimitivity \cite{W26}.
\item It has abelianisation $C_2\times C_2\times C_2\times C_2$
	-- 
the largest possible group (Theorem \ref{G(n,a,b,r)-abel-thm}).
\item It is intimately involved in the construction of 
	the primitive quaternionic reflection groups of type $P$
		and the five MUBs (mutually unbiased bases) in $\HH^2$
		\cite{BW26}.
\end{itemize}
It is generated by any four of its five reflection orbits, thereby 
giving six labels for $K$, and this accounts for all the 
multiple labelling, i.e., $K$ has $2^5-5=27$ normal reflection
subgroups. The other reflection groups with five reflection orbits
(with type $(2,2,r:2,2)$ and $2^4\cdot r$ labels) have
similar multiple labelling, e.g., we have the following numbers 
of labels and normal reflection subgroups
(see Table \ref{G(n,a,b,r)groups-table})
$$ G(4,1,2,2):\ 32,23, \quad G(4,1,1,4):\ 48,34, \quad
G(6,1,1,6):\ 64,41, \quad G(6,1,3,2):\ 32,27. $$
\end{example}

Another mechanism for there to be multiple labelling is the 
existence of a reflection of order $r\ne2$ (which is generated 
by the reflections of order $2$).

\begin{example} 
\label{genbyrefsordertwo}
The indices below give groups $G(n,a,b,r)$, $r\ne2$, 
which are generated by the reflections of order $2$, 
i.e., these reflections generate the reflections of order $r\ne2$
$$ [ 3, 1, 1, 3 ], \quad [ 4, 1, 1, 4 ], \quad [ 5, 1, 1, 5 ], \quad
[ 6, 1, 2, 3 ], \quad [ 6, 1, 1, 6 ], \quad [ 7, 1, 1, 7 ], \quad 
[ 8, 1, 2, 4 ], \quad \ldots.  $$
We consider the corresponding multiple labelling for the first
$$ G(3,1,1,3) =\cH_1(G(6,2,2))=\inpro{G(6,2,2),\zeta_6,s}, $$
which is a near miss for unique labelling.
The complex reflection group $G(6,2,2)$ has reflection type $3C_2,3C_2,2C_3$,
with the $3C_2,3C_2$ combining to a single orbit in $G(3,1,1,3)$,
which therefore has reflection type $6C_2,2C_3:6C_2$
Direct calculation shows that
        $$ H^{(2,3:2)}=G(3,1,1,3)= H^{(2,1:2)}, $$
i.e., $G(3,1,1,3)$ has two labels, and 
it is generated by its reflections of order $2$.
This behaviour does not happen for complex reflection groups.
\end{example}

It appears that the only case when $G=G(n,a,b,r)$ 
has a reflection of order $r\ne 2$ which is required 
in a generating set of reflections for $G$ is
$$ [n,a,b,r] = [n,1,b,{2n\over b}], \quad b\ne n \quad
\hbox{(higher order group, $b$ odd)}, $$
and these groups do have subgroups with multiple labels.


We have not looked exhaustively at the lattice of normal
reflection subgroups of the quaternionic reflection groups $\cH_a(G)$, 
but this could, in principle, be done using the techniques of this paper 
-- in particular to find a
simple condition which ensures unique labels (or not).
The mechanisms which are likely to lead to multiple labels include:
\begin{itemize}
\item Many reflection orbits, in particular two type $Bj$ orbits.
\item Reflections of order greater than $2$.
\item The number of factors of the abelianisation being less 
	than the number of reflection orbits.
\end{itemize}

We conclude with a general summary of our current understanding.

\begin{theorem} 
\label{HaGnormalrefgroups}
(Normal reflection subgroups) 
Let $\cH_a(G)=H^{(k_1,\ldots,k_m:\ell_1,\ldots,\ell_{m'})}$ be
one of the quaternionic reflection groups of 
Theorems \ref{TOIfamilyclass} and \ref{G(n,a,b,r)familyclass}.
Then
$$ (\ga_1,\ldots,\ga_m:\gb_1,\ldots,\gb_{m'}) \mapsto
 H^{(\ga_1,\ldots,\ga_m:\gb_1,\ldots,\gb_{m'})} $$
is a lattice homomorphism from the divisors of $(k_1,\ldots,k_m:\ell_1,\ldots,\ell_{m'})$ onto the normal reflection subgroups of $\cH_a(G)$.
This can be one-to-one (unique labelling) or not (multiple labelling),
as detailed in the Examples \ref{oneBjorbit}, \ref{TwoBjorbits}, 
\ref{Uniquelabellingexample}, \ref{fivereflectionorbits}, 
\ref{genbyrefsordertwo}.
\end{theorem}

\section{Appendix}

We constructed and calculated with the groups in this paper in the
{\tt magma} computer algebra system. We give some sample code (which may be
cut and pasted), to show how the quaternionic algebras were used
(there are some subtle technical issues). Code for the groups
$G(n,a,b,r)$ and the indices $\gL_n$ can be found in \cite{W25}.

\begin{verbatim}

// Given a matrix group H over a cyclotomic field F
// extend the field to K and define M_a(H) as group over K 
// K is chosen to contain i, so as to avoid issues when j is added
Ma := function(a,H);
  F:=BaseRing(H);
  K:=CyclotomicField(LCM(CyclotomicOrder(F)*a,4));
  G:=ChangeRing(H,K);
  zeta:=ScalarMatrix(K,2,RootOfUnity(Order(Centre(H))*a));
  return MatrixGroup<2,K|Generators(G) join {zeta}>;
end function;
\end{verbatim}

\begin{verbatim}

// Given a matrix group H over a cyclotomic field F
// extend the field to K and define G=M_a(H) as group over K
// extend K to a quaternionic algebra Q
// return the matrix group over Q generated by G and s
// _,Gp:=Ha(1,H) returns Gp the group H over Q
Ha := function(a,H);
  G:=Ma(a,H);
  F:=BaseRing(G);
  Q<i,j,k>:=QuaternionAlgebra<F|-1,-1>;
  s:=Matrix(Q,2,2,[0,j,-j,0]);
  Ggens:={Parent(s)![(z+ComplexConjugate(z))/2+(z-ComplexConjugate(z))
          /(2*RootOfUnity(4))*i:z in Eltseq(B)]: B in Generators(G)};
  return MatrixGroup<2,Q|Ggens join {s}>, MatrixGroup<2,Q|Ggens>;
end function;

// Check whether a matrix with entries in a quaternionic algebra
// over a cyclotomic field F is in the field F
// This can be used to see if a matrix is of type A or type Bj
IsComplex := function(g);
  quatparts:=&join[{ElementToSequence(z)[a]:a in [3..4]}:z in Eltseq(g)];
  return quatparts eq {0};
end function;





// The group H_r^A4, for r=0,1,2,...
// The choice r=-1 gives the binary tetrahedral group T
HA4 := function(r);
  F:=CyclotomicField(4*3^(r+1));
  i:=F!RootOfUnity(4); R:=Matrix(F,2,2,[1,0,0,i^3]);
  hZ:=1/2*Matrix(F,2,2,[-1+i,-1+i,1+i,-1-i]);
  return MatrixGroup<2,F|i*R^2,RootOfUnity(3^(r+1))*hZ>;
end function;

// The minimal group H_r^(S_4), for r=0,1,2,...
// The choice r=-1 gives the binary octahedral group O
HS4 := function(r);
  F:=CyclotomicField(LCM(8,2^(r+2)*3));
  i:=F!RootOfUnity(4); hR:=RootOfUnity(8)*Matrix(F,2,2,[1,0,0,i^3]);
  hZ:=1/2*Matrix(F,2,2,[-1+i,-1+i,1+i,-1-i]);
  return MatrixGroup<2,F|RootOfUnity(2^(r+2))*hR,hZ>;
end function;

// The minimal group I
bI := function();
  F:=CyclotomicField(120);
  i:=F!RootOfUnity(4); z5:=F!RootOfUnity(5); tau:=(z5+z5^-1+1);
  R:=Matrix(F,2,2,[1,0,0,i^3]);
  hM:=1/2*Matrix(F,2,2,[tau^-1-tau*i,-1,1,tau^-1+tau*i]);
  return MatrixGroup<2,F|i*R^2,hM>;
end function;

// The minimal groups H_r^(D_n) and K_r^(D_n), for r=0,1,2,...
// The choice r=-1 gives the binary dihedral group of order 4n
// The third argument gives the dihedral group D_n of order 2n (for n odd)
HKDn := function(n,r);
  F:=CyclotomicField(LCM(4*n,2^(r+2)));
  z:=F!RootOfUnity(2*n); i:=F!RootOfUnity(4);
  A:=Matrix(F,2,2,[z,0,0,z^-1]);
  B:=Matrix(F,2,2,[0,i,i,0]);  // B:=Matrix(F,2,2,[0,1,-1,0]);
  H:=MatrixGroup<2,F|A,RootOfUnity(2^(r+2))*B>;
  K:=MatrixGroup<2,F|RootOfUnity(2^(r+2))*A,RootOfUnity(2^(r+2))*B>;
  Dn:=MatrixGroup<2,F|-A,i*B>;
  return H,K,Dn;
end function;

\end{verbatim}

\begin{verbatim}
// Calculates the primitive subgroups G(a,H,G) of GL(2,C)
// Call up as PrimG(a:G,H) assumes two or three arguments with (a:G) implying H=G
// with H,G taken from {T,D2,O,I}, the order being unimportant
PrimG := function(a : T:=false, D2:=false, O:=false, I:=false);
  if T and D2 then r,c:=Valuation(a,3); return Ma(c,HA4(r)); end if;
  if O and T  then r,b:=Valuation(a,2); return Ma(b,HS4(r)); end if;
  if T then return Ma(a,HA4(-1)); end if;
  if O then return Ma(a,HS4(-1)); end if;
  if I then return Ma(a,bI()); end if;
  return "There is no group with these parameters";
end function;

// Calculates the imprimitive irreducible subgroups G(a,D_n,H) of GL(2,C)
// Call up as ImPrimG(a,n:H) with H taken from {Dn,C2n,Cn,Dn2}
ImPrimG := function(a,n : Dn:=false, C2n:=false, Cn:=false, Dn2:=false);
  error if n lt 2, "n must be >= 2";
  error if not Type(a) eq RngIntElt or a lt 1, "a must be a positive integer";
  if Dn then return Ma(a,HKDn(n,-1)); end if;
  if C2n then r,b:=Valuation(a,2); return Ma(b,HKDn(n,r)); end if;
  if Dn2 and IsEven(n) then r,b:=Valuation(a,2); _,K:=HKDn(n,r); 
                            return Ma(b,K); end if;
  if Cn then
    if IsOdd(n*a) then _,_,DnOdd:=HKDn(n,0); return Ma(a,DnOdd);
                  else error "There are no groups G(a,Dn,Cn) for na even";
    end if;
  end if;
  return "There is no group with these parameters";
end function;

// Calculate the orbits of the reflections in a group G
// This returns an orbit of an element which generates H_a
CalculateOrbits := function(G)
  d:=Ncols(Random(G)); STF:=BaseRing(G); I:=IdentityMatrix(STF,d);
  refs:={};
  for g in G do
    if Rank(g-I) eq 1 then Include(~refs,g); end if;
  end for;
  remainrefs:=refs;
  ords:=[]; orbs:=[];
  while #remainrefs gt 0 do
    // find a reflection h of maximal order
    maxord:=Maximum({Order(g):g in remainrefs});
    h:=Random({g : g in remainrefs | Order(g) eq maxord });
    horbit:={g*h*g^-1:g in G}; Include(~orbs,horbit);
    remainrefs:=remainrefs diff {g^j:g in horbit,j in [1..maxord-1]};
  end while;
  return orbs;
end function;
\end{verbatim}

\bibliographystyle{alpha}
\bibliography{references}
\nocite{*}

\end{document}